
%
%
%
%
%
%
\magnification=\magstephalf      
%
%
\vsize=7.5truein                 
\hsize=5.2truein                 
\newskip\stdskip                 
\stdskip=6pt plus3pt minus3pt    
\medskipamount=\stdskip          
\parindent=0pt                   
\parskip=\stdskip                
\abovedisplayskip=\stdskip       
\belowdisplayskip=\stdskip       
\mathsurround=0.75pt             
\overfullrule=0pt                
%
%
\def\ppar{\par\goodbreak\vskip 8pt plus 4pt minus 4pt}     
%
%
\def\stdspace{\hskip 0.75em plus 0.15em\ignorespaces}
\let\qua\stdspace 
%
%
%
%
%
%
%
\def\hexnumber#1{\ifcase#1 0\or 1\or 2\or 3\or 4\or 5\or 6\or 7\or 8\or
 9\or A\or B\or C\or D\or E\or F\fi}
%
%
\font\thirtnmsa=msam10 scaled 1315    
\font\tenmsa=msam10          \font\ninemsa=msam9
\font\sevenmsa=msam7         \font\sixmsa=msam6
\font\fivemsa=msam5
%
%
\newfam\msafam                  \textfont\msafam=\tenmsa
\scriptfont\msafam=\sevenmsa    \scriptscriptfont\msafam=\fivemsa
\edef\hexa{\hexnumber\msafam}        
\def\msa{\fam\msafam\tenmsa}         
%
%
\font\thirtnmsb=msbm10 scaled 1315   
\font\tenmsb=msbm10      \font\ninemsb=msbm9
\font\sevenmsb=msbm7     \font\sixmsb=msbm6
\font\fivemsb=msbm5
%
\newfam\msbfam                   \textfont\msbfam=\tenmsb       
\scriptfont\msbfam=\sevenmsb     \scriptscriptfont\msbfam=\fivemsb
\edef\hexb{\hexnumber\msbfam}    
\def\msb{\fam\msbfam\tenmsb}     
%
%
\font\thirtneufm=eufm10 scaled 1315   
\font\teneufm=eufm10                 \font\nineeufm=eufm9
\font\seveneufm=eufm7                \font\sixeufm=eufm6
\font\fiveeufm=eufm5
%
\newfam\eufmfam                    \textfont\eufmfam=\teneufm
\scriptfont\eufmfam=\seveneufm     \scriptscriptfont\eufmfam=\fiveeufm
\edef\hexf{\hexnumber\eufmfam}      
\def\frak{\fam\eufmfam\teneufm}     
%
%
%
\font\thirtnrm=cmr10 scaled 1315    
\font\ninerm=cmr9                   \font\sixrm=cmr6   
%
\font\thirtni=cmmi10 scaled 1315    
\font\ninei=cmmi9                   \font\sixi=cmmi6  
%
\font\thirtnsy=cmsy10 scaled 1315   
\font\ninesy=cmsy9                  \font\sixsy=cmsy6  
%
\font\thirtnbf=cmbx10 scaled 1315   
\font\ninebf=cmbx9                  \font\sixbf=cmbx6  
%
%
\font\thirtnex=cmex10 scaled 1315   
\font\nineex=cmex9                  
%
%
\font\thirtnit=cmti10 scaled 1315  
\font\nineit=cmti9                  
%
\font\thirtnsl=cmsl10 scaled 1315  
\font\ninesl=cmsl9                  
%
\font\thirtntt=cmtt10 scaled 1315  
\font\ninett=cmtt9                  
%
%
%
%
\def\small{%
%
%
\textfont0=\ninerm \scriptfont0=\sixrm \scriptscriptfont0=\fiverm
\def\rm{\fam0\ninerm}
%
%
\textfont1=\ninei \scriptfont1=\sixi \scriptscriptfont1=\fivei
%
%
\textfont2=\ninesy \scriptfont2=\sixsy \scriptscriptfont2=\fivesy
%
%
\textfont3=\nineex \scriptfont3=\nineex \scriptscriptfont3=\nineex
%
%
\textfont\bffam=\ninebf \scriptfont\bffam=\sixbf
\scriptscriptfont\bffam=\fivebf \def\bf{\fam\bffam\ninebf}%
%
%
\textfont\itfam=\nineit \def\it{\fam\itfam\nineit}%
\textfont\slfam=\ninesl \def\sl{\fam\slfam\ninesl}%
\textfont\ttfam=\ninett \def\tt{\fam\ttfam\ninett}%
%
%
%
\textfont\msafam=\ninemsa \scriptfont\msafam=\sixmsa
\scriptscriptfont\msafam=\fivemsa \def\msa{\fam\msafam\ninemsa}%
%
%
\textfont\msbfam=\ninemsb \scriptfont\msbfam=\sixmsb
\scriptscriptfont\msbfam=\fivemsb \def\msb{\fam\msbfam\ninemsb}%
%
%
\textfont\eufmfam=\nineeufm  \scriptfont\eufmfam=\sixeufm
\scriptscriptfont\eufmfam=\fiveeufm \def\frak{\fam\eufmfam\nineeufm}%
%
%
%
\normalbaselineskip=11pt%
\setbox\strutbox=\hbox{\vrule height8pt depth3pt width0pt}%
%
%
\normalbaselines\rm
%
%
\stdskip=4pt plus2pt minus2pt    
\medskipamount=\stdskip          
\parskip=\stdskip                
\abovedisplayskip=\stdskip       
\belowdisplayskip=\stdskip       
\def\ppar{\par\goodbreak\vskip 6pt plus 3pt minus 3pt}%
%
%
\def\section##1{\global\advance\sectionnumber by 1
\vskip-\lastskip\penalty-800\vskip 20pt plus10pt minus5pt 
\egroup{\bf\number\sectionnumber\quad##1}\bgroup\small         
\vskip 6pt plus3pt minus3pt
\nobreak\resultnumber=1}
}    
%
\def\beginsmall{\bgroup\small}
\let\endsmall\egroup
%
%
%
%
\def\large{%
\textfont0=\thirtnrm \scriptfont0=\ninerm \scriptscriptfont0=\sevenrm
\def\rm{\fam0\thirtnrm}%
\textfont1=\thirtni \scriptfont1=\ninei \scriptscriptfont1=\seveni
\textfont2=\thirtnsy \scriptfont2=\ninesy \scriptscriptfont2=\sevensy
\textfont3=\thirtnex \scriptfont3=\thirtnex \scriptscriptfont3=\thirtnex
\textfont\bffam=\thirtnbf \scriptfont\bffam=\ninebf
\scriptscriptfont\bffam=\sevenbf \def\bf{\fam\bffam\thirtnbf}%
\textfont\itfam=\thirtnit \def\it{\fam\itfam\thirtnit}%
\textfont\slfam=\thirtnsl \def\sl{\fam\slfam\thirtnsl}%
\textfont\ttfam=\thirtntt \def\tt{\fam\ttfam\thirtntt}%
\textfont\msafam=\thirtnmsa \scriptfont\msafam=\ninemsa
\scriptscriptfont\msafam=\sevenmsa \def\msa{\fam\msafam\thirtnmsa}%
\textfont\msbfam=\thirtnmsb \scriptfont\msbfam=\ninemsb
\scriptscriptfont\msbfam=\sevenmsb \def\msb{\fam\msbfam\thirtnmsb}%
\textfont\eufmfam=\thirtneufm  \scriptfont\eufmfam=\nineeufm
\scriptscriptfont\eufmfam=\seveneufm \def\frak{\fam\eufmfam\teneufm}%
\normalbaselineskip=16pt%
\setbox\strutbox=\hbox{\vrule height11.5pt depth4.5pt width0pt}%
\normalbaselines\rm}%
\let\Large\large   
%
\def\Bbb#1{{\msb#1}}

%

%
\mathchardef\plussquare="0\hexa01
\mathchardef\nge="3\hexb0B
\mathchardef\maltesecross="0\hexa7A
\mathchardef\del="0\hexf01
%
%
%
%
\font\sc=cmcsc10
%
%
%
%
\def\sqr#1#2{{\vcenter{\vbox{\hrule  height.#2truept
	\hbox{\vrule width.#2truept height#1truept 
	\kern#1truept \vrule width.#2truept}
	\hrule height.#2truept}}}}
\def\sq{\sqr55}    
%
%
%
%
\newcount\sectionnumber            
\newcount\resultnumber             
\sectionnumber=0\resultnumber=1    
%
%
%
\def\section#1{\global\advance\sectionnumber by 1
\xdef\nextkey{\number\sectionnumber}
\vskip-\lastskip\penalty-800\vskip 20pt plus10pt minus5pt 
{\large\bf\number\sectionnumber\quad#1}         
\vskip 8pt plus4pt minus4pt
\nobreak\resultnumber=1}      
%
%
%
%
%
\def\sh#1{\vskip-\lastskip\ppar{\bf #1}\par\nobreak\medskip}         
%
%
%
%

%
\def\proc#1{\xdef\nextkey{\number\sectionnumber.\number\resultnumber}%
\vskip-\lastskip\ppar\bf%
\noindent#1\ \number\sectionnumber.\number\resultnumber
\stdspace\sl\global\advance\resultnumber by 1\ignorespaces}
 
%
%
\def\qed{\hfill$\sq$\par\goodbreak\rm}   
%
%
%
%
%
%
%
%
\def\proclaim#1{\vskip-\lastskip\ppar\bf%
\noindent#1\stdspace\sl\ignorespaces} 
\let\endproclaim\endproc
%
%
%
%
\def\rk#1{\vskip-\lastskip\ppar{\bf #1}\stdspace\ignorespaces}                

%
%
%
%
%
%
\def\label{\xdef\nextkey{\number\sectionnumber.\number\resultnumber}%
\number\sectionnumber.\number\resultnumber
\global\advance\resultnumber by 1}
%
%
%
%
%
%
%
%
%
%
%
%
%
%
%
%
\newcount\refnumber              
\refnumber=1                     
\long\def\reflist#1\endreflist{%
\long\def\thereflist{#1}{\def\refkey##1##2\par{\xdef##1{\number\refnumber}%
\global\advance\refnumber by 1}%
\def\key##1##2\par{\expandafter\xdef%
\csname##1\endcsname{\number\refnumber}%
\global\advance\refnumber by 1}#1\par}}
\long\def\references{%
\penalty-800\vskip-\lastskip\vskip 15pt plus10pt minus5pt 
{\large\bf References}\ppar 
{\leftskip=25pt\frenchspacing    
\small\parskip=3pt plus2pt       
\def\refkey##1##2\par{\noindent  
\llap{[##1]\stdspace}\ignorespaces##2\par}         
\def\key##1##2\par{\noindent  
\llap{[\ref{##1}]\stdspace}\ignorespaces##2\par}  
\def\,{\thinspace}\thereflist\par}}
%
%
%
\newcount\footnotenumber         
\footnotenumber=1                
\def\fnote#1{\xdef\nextkey{\number\footnotenumber}%
{\small\ifnum\footnotenumber>9\parindent=14pt%
\else\parindent=10pt\fi\footnote{$^{\number\footnotenumber}$}%
{\hglue-5pt#1}\global\advance\footnotenumber by 1}}
%
%
%
%
%
%
%
\newcount\figurenumber          
\figurenumber=1                 
\def\caption#1{\xdef\nextkey{\number\figurenumber}%
\cl{\small Figure \number\figurenumber: #1}%
\global\advance\figurenumber by 1}
\def\figurelabel{\xdef\nextkey{\number\figurenumber}%
\cl{\small Figure \number\figurenumber}%
\global\advance\figurenumber by 1}
\long\def\figure#1\endfigure{{\xdef\nextkey{\number\figurenumber}%
\let\captiontext\relax\def\caption##1{\xdef\captiontext{##1}}%
\midinsert\cl{\ignorespaces#1\unskip\unskip\unskip\unskip}\vglue6pt\cl{\small 
Figure \number\figurenumber\ifx\captiontext\relax\else: \captiontext
\fi}\endinsert\global\advance\figurenumber by 1}}
%
%
%
%
%
%
%
\def\nextkey{??}   
%
\def\key#1{\expandafter\xdef\csname #1\endcsname{\nextkey}}
\def\ref#1{\expandafter\ifx\csname #1\endcsname\relax
\immediate\write16{Reference {#1} undefined}??\else
\csname #1\endcsname\fi}
%
%
%
%
%
%
%
\newread\gtinfile
\newwrite\gtreffile
\def\useforwardrefs{
\openin\gtinfile\jobname.ref
\ifeof\gtinfile
\closein\gtinfile
\immediate\write16{No file \jobname.ref}
\else
\closein\gtinfile
\input \jobname.ref
\fi
\immediate\openout\gtreffile \jobname.ref
%
%
\def\key##1{{\def\\{\noexpand}%
\expandafter\xdef\csname ##1\endcsname{\nextkey}%
\immediate\write\gtreffile{\\\expandafter\\\def\\\csname ##1\\\endcsname%
{\nextkey}}}}
%
%
\long\def\reflist##1\endreflist{%
\long\def\thereflist{##1}{\def\refkey####1####2\par{\xdef####1{%
\number\refnumber}{\def\\{\noexpand}\immediate\write\gtreffile
{\\\def\\####1{\number\refnumber}}}\global\advance\refnumber by 1}%
\def\key####1####2\par{\expandafter\xdef%
\csname####1\endcsname{\number\refnumber}%
{\def\\{\noexpand}\immediate\write\gtreffile
{\\\expandafter\\\def\\\csname ####1\\\endcsname{\number\refnumber}}}
\global\advance\refnumber by 1}##1\par}}
\long\def\biblio##1\endbiblio{\reflist##1\endreflist\references}%
%
%
\def\numkey##1{{\def\\{\noexpand}%
\xdef##1{\number\sectionnumber.\number\resultnumber}
\immediate\write\gtreffile{\\\def\\##1%
{\number\sectionnumber.\number\resultnumber}}}}
\def\seckey##1{{\def\\{\noexpand}\xdef##1{\number\sectionnumber}
\immediate\write\gtreffile{\\\def\\##1{\number\sectionnumber}}}}
\def\figkey##1{\xdef##1{\number\figurenumber}%
{\def\\{\noexpand}\immediate\write\gtreffile%
{\\\def\\##1{\number\figurenumber}}}
\number\figurenumber\global\advance\figurenumber by 1}
}   
%
%
%
%
\def\figkey#1{\xdef#1{\number\figurenumber}%
\number\figurenumber\global\advance\figurenumber by 1}
\def\fig#1#2\endfig{%
\midinsert\cl{#2}\vglue6pt\cl{\small Figure #1}\endinsert}
\def\newfig{\number\figurenumber\global\advance\figurenumber by 1}
\def\numkey#1{\xdef#1{\number\sectionnumber.\number\resultnumber}}
\def\seckey#1{\xdef#1{\number\sectionnumber}}
%
%
%
%
%
%
%
%
%
\def\verb{\catcode`\"=\active}       
\def\brev{\catcode`\"=12}            
\brev                                
\verb                                
{\obeyspaces\gdef {\ }}              
{\catcode`\`=\active\gdef`{\relax\lq}}
\def"{%
\begingroup\baselineskip=12pt\def\par{\leavevmode\endgraf}%
\tt\obeylines\obeyspaces\parskip=0pt\parindent=0pt%
\catcode`\$=12\catcode`\&=12\catcode`\^=12\catcode`\#=12%
\catcode`\_=12\catcode`\~=12%
\catcode`\{=12\catcode`\}=12\catcode`\%=12\catcode`\\=12%
\catcode`\`=\active\let"\endgroup}
\brev      
%
%
%
%
%
%
\def\item#1{\par\leavevmode\llap{#1\stdspace}%
\ignorespaces}                             
%
%

%
%
\def\np{\vfil\eject}         
\def\nl{\hfil\break}         
\def\cl{\centerline}         
\def\agt{{\mathsurround=0pt\it$\cal A\mskip-.7mu$lgebraic \&\ 
$\cal G\mskip-2mu$eometric $\cal T\!\!$opology}}  
%
%
%

%
%
%
%
%
\def\title#1{\def\thetitle{#1}}
\def\shorttitle#1{\def\theshorttitle{#1}}
\def\author#1{\edef\previousauthors{\theauthors}
 \ifx\theauthors\relax\def\theauthors{#1}\else
 \def\theauthors{\previousauthors\par#1}\fi}

%
\def\address#1{\edef\previousaddresses{\theaddress}
 \ifx\theaddress\relax\def\theaddress{#1}\else
 \def\theaddress{\previousaddresses\par\vskip 2pt\par#1}\fi}
\def\secondaddress#1{\edef\previousaddresses{\theaddress}
 \ifx\theaddress\relax\def\theaddress{#1}\else
 \def\theaddress{\previousaddresses\par{\rm and}\par#1}\fi}   

\def\email#1{\edef\previousemails{\theemail}
 \ifx\theemail\relax\def\theemail{#1}\else
 \def\theemail{\previousemails\hskip 0.75em\relax#1}\fi}
\def\secondemail#1{\edef\previousemails{\theemail}
 \ifx\theemail\relax\def\theemail{#1}\else
 \def\theemail{\previousemails\hskip 0.75em{\rm and}\hskip 0.75em
 \relax#1}\fi}
\def\url#1{\edef\previousurls{\theurl}
 \ifx\theurl\relax\def\theurl{#1}\else
 \def\theurl{\previousurls\hskip 0.75em\relax#1}\fi}
\def\secondurl#1{\edef\previousurls{\theurl}
 \ifx\theurl\relax\def\theurl{#1}\else
 \def\theurl{\previousurls\hskip 0.75em{\rm and}\hskip 0.75em
 \relax#1}\fi}
\long\def\abstract#1\endabstract{\long\def\theabstract{#1}}
\def\primaryclass#1{\def\theprimaryclass{#1}}
\let\subjclass\primaryclass                        
\def\secondaryclass#1{\def\thesecondaryclass{#1}}
\def\keywords#1{\def\thekeywords{#1}}
%
%
\let\\\par\let\thetitle\relax\let\theshorttitle\relax
\let\theauthors\relax\let\theshortauthors\relax
\let\theaddress\relax\let\theshortaddress\relax
\let\theemail\relax\let\theurl\relax
\let\theabstract\relax\let\theprimaryclass\relax
\let\thesecondaryclass\relax\let\thekeywords\relax
%
%
%
%
\long\def\maketitlepage{    

\vglue 0.2truein   

%
{\parskip=0pt\leftskip 0pt plus 1fil\def\\{\par\smallskip}{\large
\bf\thetitle}\par\medskip}   

\vglue 0.15truein 

%
{\parskip=0pt\leftskip 0pt plus 1fil\def\\{\par}{\sc\theauthors}
\par\medskip}%
 
\vglue 0.1truein 

%
{\small\parskip=0pt
{\leftskip 0pt plus 1fil\def\\{\par}{\sl\theaddress}\par}
\ifx\theemail\relax\else  
\vglue 5pt \def\\{\stdspace{\rm and}\stdspace} 
\cl{Email:\stdspace\tt\theemail}\fi
\ifx\theurl\relax\else    
\vglue 5pt \def\\{\stdspace{\rm and}\stdspace} 
\cl{URL:\stdspace\tt\theurl}\fi\par}

\vglue 7pt 

{\bf Abstract}

\vglue 5pt

\theabstract

\vglue 7pt 

{\bf AMS Classification numbers}\quad Primary:\quad \theprimaryclass\par

Secondary:\quad \thesecondaryclass

\vglue 5pt 

{\bf Keywords:}\quad \thekeywords

\np  

}    
%
%
\long\def\makeshorttitle{    


%
{\parskip=0pt\leftskip 0pt plus 1fil\def\\{\par\smallskip}{\large
\bf\thetitle}\par\medskip}   

\vglue 0.05truein 

%
{\parskip=0pt\leftskip 0pt plus 1fil\def\\{\par}{\sc\theauthors}
\par\medskip}%
 
\vglue 0.03truein 

%
{\small\parskip=0pt
{\leftskip 0pt plus 1fil\def\\{\par}{\sl\ifx\theshortaddress\relax
\theaddress\else\theshortaddress\fi}\par}
\ifx\theemail\relax\else  
\vglue 5pt \def\\{\stdspace{\rm and}\stdspace} 
\cl{Email:\stdspace\tt\theemail}\fi
\ifx\theurl\relax\else    
\vglue 5pt \def\\{\stdspace{\rm and}\stdspace} 
\cl{URL:\stdspace\tt\theurl}\fi\par}

\vglue 10pt 


{\small\leftskip 25pt\rightskip 25pt{\bf Abstract}\stdspace\theabstract

{\bf AMS Classification}\stdspace\theprimaryclass
\ifx\thesecondaryclass\relax\else; \thesecondaryclass\fi\par
{\bf Keywords}\stdspace \thekeywords\par}
\vglue 7pt
}    
\let\maketitle\makeshorttitle        
%
%

\def\volumenumber#1{\def\thevolumenumber{#1}}
\def\volumeyear#1{\def\thevolumeyear{#1}}
\def\pagenumbers#1#2{\def\startpage{#1}\def\finishpage{#2}}
\def\published#1{\def\publishdate{#1}}
\def\received#1{\def\receiveddate{#1}}
\def\revised#1{\def\reviseddate{#1}}
\let\reviseddate\relax
\volumenumber{X}
\volumeyear{20XX}
\pagenumbers{1}{XXX}
\published{XX Xxxember 20XX}

\long\def\makeagttitle{   
\agt\hfill      
\hbox to 60truept{\vbox to 0pt{\vglue -14truept{\bf [Logo here]}\vss}\hss}
\break
{\small Volume \thevolumenumber\ (\thevolumeyear)
\startpage--\finishpage\nl
Published: \publishdate}

\vglue .2truein

{\parskip=0pt\leftskip 0pt plus 1fil\def\\{\par\smallskip}{\large
\bf\thetitle}\par\medskip}   
\vglue 0.05truein 

%
{\parskip=0pt\leftskip 0pt plus 1fil\def\\{\par}{\sc\theauthors}
\par\medskip}%
 
\vglue 0.03truein 


{\small\leftskip 25truept\rightskip 25truept{\bf Abstract}\stdspace\theabstract

{\bf AMS Classification}\stdspace\theprimaryclass
\ifx\thesecondaryclass\relax\else; \thesecondaryclass\fi\par
{\bf Keywords}\stdspace \thekeywords\par}\vglue 7truept

}   


\def\Addresses{\bigskip
{\small \parskip 0pt \leftskip 0pt \rightskip 0pt plus 1fil \def\\{\par}
\sl\theaddress\par\medskip \rm Email:\stdspace\tt\theemail\par
\ifx\theurl\relax\else\smallskip \rm URL:\stdspace\tt\theurl\par\fi}}

\def\agtart{
\hoffset 14truemm
\voffset 31truemm
\font\phead=cmsl9 scaled 950
\font\pnum=cmbx10 scaled 913
\font\pfoot=cmsl9 scaled 950
\headline{\vbox to 0pt{\vskip -4.5mm\line{\small\phead\ifnum
\count0=\startpage ISSN numbers are printed here
\hfill {\pnum\folio}\else\ifodd\count0\def\\{ }%
\ifx\theshorttitle\relax\thetitle\else\theshorttitle\fi\hfill{\pnum\folio}
\else\def\\{ and }{\pnum\folio}\hfill\ifx\theshortauthors\relax\theauthors
\else\theshortauthors\fi\fi\fi}\vss}}
\footline{\vbox to 0pt{\vglue 0mm\line{\small\pfoot\ifnum\count0=\startpage
Copyright declaration is printed here\hfill\else
\agt, Volume \thevolumenumber\ (\thevolumeyear)\hfill\fi}\vss}}
\let\maketitle\makeagttitle\let\makeshorttitle\makeagttitle}


\def\ifplaintex{\expandafter\ifx\csname documentclass\endcsname\relax}

\def\gtp{{\mathsurround=0pt\it $\cal G\mskip-2mu$eometry \&\ 
$\cal T\!\!$opology $\cal P\!$ublications}}  

\def\recd{{\small Received:\qua\receiveddate\ifx\reviseddate\relax
\else\qquad Revised:\qua\reviseddate\fi\par}} 


\def\lognumber#1{\def\thelognumber{#1}}
\def\volumenumber#1{\def\thevolumenumber{#1}}
\def\volumeyear#1{\def\thevolumeyear{#1}}
\def\papernumber#1{\def\thepapernumber{#1}}
\def\pagenumbers#1#2{\def\startpage{#1}\def\finishpage{#2}}
\def\published#1{\def\publishdate{#1}}

\def\received#1{\def\receiveddate{#1}}
\def\revised#1{\def\reviseddate{#1}}
\def\accepted#1{\def\accepteddate{#1}}

\long\def\asciiabstract#1{\long\def\theasciiabstract{#1}}
\def\asciikeywords#1{\def\theasciikeywords{#1}}


\let\\\par\let\thelognumber\relax\let\thevolumenumber\relax
\let\thepapernumber\relax\let\thevolumeyear\relax\let\startpage\relax
\let\finishpage\relax\let\publishdate\relax\let\receiveddate\relax
\let\reviseddate\relax\let\accepteddate\relax\let\theasciititle\relax
\let\theasciiauthors\relax
\let\theasciiabstract\relax\let\theasciikeywords\relax

\let\theasciiemail\relax


\ifplaintex
\font\logobig=cmssbx10 scaled 3836
\font\logomed=cmssbx10 scaled 2557
\else
\font\logobig=cmssbx10 scaled 4200
\font\logomed=cmssbx10 scaled 2800
\fi

\long\def\makeagttitle{   
\count0=\startpage
\agt\hfill      
\hbox to 45truept{\vbox to 0pt{\vglue -13truept{\logomed A\kern -.37em{\logobig 
T}\kern -.38em G}\vss}\hss}
\break
{\small Volume \thevolumenumber\ (\thevolumeyear)
\startpage--\finishpage\nl
Published: \publishdate}

\vglue .25truein

{\parskip=0pt\leftskip 0pt plus
1fil\def\\{\par\smallskip}{\Large\bf\thetitle}\par\medskip} \vglue
0.05truein

%
{\parskip=0pt\leftskip 0pt plus 1fil\def\\{\par}{\sc\theauthors}
\par\medskip}%
 
\vglue 0.03truein 


{\small\leftskip 25truept\rightskip 25truept{\bf Abstract}\stdspace\theabstract

{\bf AMS Classification}\stdspace\theprimaryclass
\ifx\thesecondaryclass\relax\else; \thesecondaryclass\fi\par
{\bf Keywords}\stdspace \thekeywords\par}\vglue 7truept

}   

\ifplaintex
\hoffset 14truemm
\voffset 31truemm
\font\phead=cmsl9 scaled 950
\font\pnum=cmbx10 scaled 913
\font\pfoot=cmsl9 scaled 950
\headline{\vbox to 0pt{\vskip -4.5mm\line{\small\phead\ifnum
\count0=\startpage ISSN 1472-2739 (on-line) 1472-2747 (printed)
\hfill {\pnum\folio}\else\ifodd\count0\def\\{ }%
\ifx\theshorttitle\relax\thetitle\else\theshorttitle\fi\hfill{\pnum\folio}
\else\def\\{ and }{\pnum\folio}\hfill\ifx\theshortauthors\relax\theauthors
\else\theshortauthors\fi\fi\fi}\vss}}
\footline{\vbox to 0pt{\vglue 0mm\line{\small\pfoot\ifnum\count0=\startpage
\copyright\ \gtp\hfill\else
\agt, Volume \thevolumenumber\ (\thevolumeyear)\hfill\fi}\vss}}
\else
\headsep 23pt
\footskip 35pt
\hoffset -4truemm
\voffset 12.5truemm
\font\lhead=cmsl9 scaled 1050
\font\lnum=cmbx10 
\font\lfoot=cmsl9 scaled 1050
\makeatletter
\def\@oddhead{{\small\lhead\ifnum\count0=\startpage ISSN 1472-2739 
(on-line) 1472-2747 (printed)\hfill {\lnum\number\count0}\else\ifodd\count0
\def\\{ }\ifx\theshorttitle\relax \thetitle \else\theshorttitle\fi\hfill
{\lnum\number\count0}\else\def\\{ and }{\lnum\number\count0}
\hfill\ifx\theshortauthors\relax 
\theauthors\else\theshortauthors\fi\fi\fi}}\def\@evenhead{\@oddhead}
\def\@oddfoot{\small\lfoot\ifnum\count0=\startpage\copyright\ \gtp\hfill\else
\agt, Volume \thevolumenumber\ (\thevolumeyear)\hfill\fi}
\def\@evenfoot{\@oddfoot}
\makeatother
\fi
\let\maketitlepage\makeagttitle
\let\makeshorttitle\maketitlepage
\let\maketitle\maketitlepage


\newwrite\gtoutfile
\long\gdef\makeheadfile{  
{\def\\{, }\def\s{ }
\immediate\openout\gtoutfile head.xxx
\immediate\write\gtoutfile{To: math@arxiv.org}
\immediate\write\gtoutfile{Subject: put OR rep NNNNN:ppppp}
\immediate\write\gtoutfile{--text follows this line--}
\immediate\write\gtoutfile{Proxy-for: \ifx\theasciiauthors\relax
\theauthors\else\theasciiauthors\fi\s<\ifx\theasciiemail\relax\theemail\else\theasciiemail\fi>}
\immediate\write\gtoutfile{\noexpand\\}
\immediate\write\gtoutfile{Authors: \ifx\theasciiauthors\relax
\theauthors\else\theasciiauthors\fi}
{\def\\{ }\immediate\write\gtoutfile{Title: \ifx\theasciititle\relax
\thetitle\else\theasciititle\fi}}
\immediate\write\gtoutfile{Subj-class: GT or SG, GR etc}
\immediate\write\gtoutfile{MSC-class: \theprimaryclass\ifx\thesecondaryclass\relax\else, \thesecondaryclass\fi}
\immediate\write\gtoutfile{Journal-ref: Algebr. Geom. Topol. \thevolumenumber\s
(\thevolumeyear) \startpage-\finishpage}
\immediate\write\gtoutfile{Comments: Published by Algebraic and
Geometric Topology at}
\immediate\write\gtoutfile{\s\s\s  http://www.maths.warwick.ac.uk/agt/AGTVol\thevolumenumber/agt-\thevolumenumber-\thepapernumber.abs.html}
\immediate\write\gtoutfile{\noexpand\\}
\immediate\write\gtoutfile{}
\ifx\theasciiabstract\relax
\immediate\write\gtoutfile{\theabstract}\else
\immediate\write\gtoutfile{\theasciiabstract}\fi
\immediate\write\gtoutfile{}
\immediate\write\gtoutfile{\noexpand\\}
\immediate\write\gtoutfile{}
\immediate\closeout\gtoutfile}}  

\def\maketitlepage{\makeagttitle\makeheadfile}
\let\makeshorttitle\maketitlepage
\let\maketitle\maketitlepage

\lognumber{27}
\volumenumber{3}
\volumeyear{2003}
\papernumber{27}
\published{31 August 2003}
\pagenumbers{791}{856}
\received{11 February 2002}
\revised{31 July 2003}
\accepted{20 August 2003}

\input amsnames
\input amstex

\hoffset 14truemm
\voffset 31truemm

\let\cal\Cal          
\catcode`\@=12        

\let\\\par
\def\topmatter{\relax}
\def\endtopmatter{\maketitlepage}
\let\gttitle\title
\def\title#1\endtitle{\gttitle{#1}}
\let\gtauthor\author
\def\author#1\endauthor{\gtauthor{#1}}
\let\gtaddress\address
\def\address#1\endaddress{\gtaddress{#1}}
\let\gtemail\email
\def\email#1\endemail{\gtemail{#1}}
\def\subjclass#1\endsubjclass{\primaryclass{#1}}
\let\gtkeywords\keywords
\def\keywords#1\endkeywords{\gtkeywords{#1}}
\def\heading#1\endheading{{\def\S##1{\relax}\def\\{\relax\ignorespaces}
    \section{#1}}}
\def\head#1\endhead{\heading#1\endheading}

\def\subhead#1\endsubhead{\sh{#1}}
\def\subsubhead#1\endsubsubhead{\sh{#1}}
\def\specialhead#1\endspecialhead{\sh{#1}}
\def\demo#1{\rk{#1}\ignorespaces}
\def\enddemo{\ppar}
\let\remark\demo
\def\endremark{}

\def\qed{\ifmmode\quad\sq\else\hbox{}\hfill$\sq$\par\goodbreak\rm\fi}  
\def\proclaim#1{\rk{#1}\sl\ignorespaces}
\def\endproclaim{\rm\ppar}
\def\cite#1{[#1]}
\newcount\itemnumber

\let\itemold\item
\def\item{\itemold{{\rm(\number\itemnumber)}}%
\global\advance\itemnumber by 1\ignorespaces}
\def\S{section~\ignorespaces}  
\def\date#1\enddate{\relax}
\def\thanks#1\endthanks{\relax}   
\def\dedicatory#1\enddedicatory{\relax}  
\def\rom#1{{\rm #1}}  
\let\footnote\plainfootnote

\def\Refs{\ppar{\large\bf References}\ppar\bgroup\leftskip=25pt
\frenchspacing\parskip=3pt plus2pt\small}       
\def\endRefs{\egroup}
\def\widestnumber#1#2{\relax}
\def\endrefitem{}
\def\refdef#1#2#3{\def#1{\leavevmode\unskip\endrefitem#2\def\endrefitem{#3}}}
\def\ref{\par}
\def\endref{\endrefitem\par\def\endrefitem{}}
\refdef\key{\noindent\llap\bgroup[}{]\ \ \egroup}
\refdef\no{\noindent\llap\bgroup[}{]\ \ \egroup}
\refdef\by{\bf}{\rm, }
\refdef\manyby{\bf}{\rm, }
\refdef\paper{\it}{\rm, }
\refdef\book{\it}{\rm, }
\refdef\jour{}{ }
\refdef\vol{}{ }
\refdef\yr{$(}{)$ }
\refdef\ed{(}{ Editor) }
\refdef\publ{}{ }
\refdef\inbook{from: ``}{'', }
\refdef\pages{}{ }
\refdef\page{}{ }
\refdef\paperinfo{}{ }
\refdef\bookinfo{}{ }
\refdef\publaddr{}{ }
\refdef\eds{(}{ Editors)}
\refdef\bysame{\hbox to 3 em{\hrulefill}\thinspace,}{ }
\refdef\toappear{(to appear)}{ }
\refdef\issue{no.\ }{ }
\newcount\refnumber\refnumber=1
\def\refkey#1{\expandafter\xdef\csname cite#1\endcsname{\number\refnumber}%
\global\advance\refnumber by 1}
\def\cite#1{[\csname cite#1\endcsname]}
\def\Cite#1{\csname cite#1\endcsname}  
\def\key#1{\noindent\llap{[\csname cite#1\endcsname]\ \ }}

\refkey {Ba}
\refkey {BtD}
\refkey {Br}
\refkey {BF}
\refkey {DMVV}
\refkey {DHVW}
\refkey {HH}
\refkey {HKR}
\refkey {KW}
\refkey {M1}
\refkey {M2}
\refkey {Me}
\refkey {MP}
\refkey {Sh}
\refkey {St}
\refkey {Su}
\refkey {T}
\refkey {W}

\topmatter
\title
Generalized orbifold Euler characteristics of\\
symmetric orbifolds and covering spaces
\endtitle
\shorttitle{Generalized Orbifold Euler Characteristics}
\author
Hirotaka Tamanoi
\endauthor
\address
Department of Mathematics, University of California\\ 
Santa Cruz, CA 95064, USA
\endaddress
\email
tamanoi@math.ucsc.edu
\endemail
\keywords
Automorphism group, centralizer, combinatorial group theory, covering
space, equivariant principal bundle, free group, $\Gamma$-sets,
generating function, Klein bottle genus, (non)orientable surface
group, orbifold Euler characteristic, symmetric products, twisted
sector, wreath product
\endkeywords
\asciikeywords{Automorphism group, centralizer, combinatorial 
group theory, covering space, equivariant principal bundle, free
group, Gamma-sets, generating function, Klein bottle genus,
(non)orientable surface group, orbifold Euler characteristic,
symmetric products, twisted sector, wreath product}
\primaryclass{55N20, 55N91}
\secondaryclass{57S17, 57D15, 20E22, 37F20, 05A15}
\abstract
Let $G$ be a finite group and let $M$ be a $G$-manifold. We introduce
the concept of generalized orbifold invariants of $M/G$ associated to
an arbitrary group $\Gamma$, an arbitrary $\Gamma$-set, and an
arbitrary covering space of a connected manifold $\Sigma$ whose
fundamental group is $\Gamma$. Our orbifold invariants have a natural
and simple geometric origin in the context of locally constant
$G$-equivariant maps from $G$-principal bundles over covering spaces
of $\Sigma$ to the $G$-manifold $M$. We calculate generating functions
of orbifold Euler characteristic of symmetric products of orbifolds
associated to arbitrary surface groups (orientable or non-orientable,
compact or non-compact), in both an exponential form and in an
infinite product form. Geometrically, each factor of this infinite
product corresponds to an isomorphism class of a connected covering
space of a manifold $\Sigma$. The essential ingredient for the
calculation is a structure theorem of the centralizer of homomorphisms
into wreath products described in terms of automorphism groups of
$\Gamma$-equivariant $G$-principal bundles over finite
$\Gamma$-sets. As corollaries, we obtain many identities in
combinatorial group theory. As a byproduct, we prove a simple formula
which calculates the number of conjugacy classes of subgroups of given
index in any group. Our investigation is motivated by orbifold
conformal field theory.
\endabstract
\asciiabstract{Let G be a finite group and let M be a G-manifold. 
We introduce the concept of generalized orbifold invariants of M/G
associated to an arbitrary group Gamma, an arbitrary Gamma-set, and an
arbitrary covering space of a connected manifold Sigma whose
fundamental group is Gamma. Our orbifold invariants have a natural and
simple geometric origin in the context of locally constant
G-equivariant maps from G-principal bundles over covering spaces of
Sigma to the G-manifold M. We calculate generating functions of
orbifold Euler characteristic of symmetric products of orbifolds
associated to arbitrary surface groups (orientable or non-orientable,
compact or non-compact), in both an exponential form and in an
infinite product form. Geometrically, each factor of this infinite
product corresponds to an isomorphism class of a connected covering
space of a manifold Sigma. The essential ingredient for the
calculation is a structure theorem of the centralizer of homomorphisms
into wreath products described in terms of automorphism groups of
Gamma-equivariant G-principal bundles over finite Gamma-sets. As
corollaries, we obtain many identities in combinatorial group
theory. As a byproduct, we prove a simple formula which calculates the
number of conjugacy classes of subgroups of given index in any
group. Our investigation is motivated by orbifold conformal field
theory.}

\endtopmatter

\catcode`\@=\active

\document

\head
Introduction and summary of results
\endhead

Let $G$ be a finite group acting on a finite dimensional smooth closed
manifold $M$. The action of $G$ is not assumed to be free. The orbit
space $M/G$ is an example of an orbifold. In 1980s, string
physicists suggested that in the context of orbifold conformal field
theory the ordinary Euler characteristic $\chi(M/G)$ of the orbit
space is not the right invariant, but the correct invariant is the
{\it orbifold Euler characteristic} $e_{\sssize\text{orb}}(M;G)$
defined by
$$
e_{\sssize\text{orb}}(M;G)=
\frac1{|G|}\!\!\!\sum_{\ gh=hg}\!\!\!\chi(M^{\langle g,h\rangle}),
\tag1-1
$$
where the summation runs over all pairs of commuting elements $(g,h)$
in $G$ \cite{DHVW}.

We have started the investigation of generalized orbifold Euler
characteristic in \cite{T}. We continue our investigation: Our idea in
this paper is to study orbifold singularities of $M/G$ through the use
of a connected manifold $\Sigma$ as a ``probe.''  More precisely, we
examine an infinite dimensional mapping space $\text{Map}(\Sigma,M/G)$
which can be thought of as a ``thickening'' of the orbit space $M/G$.
We consider lifting maps $\Sigma @>>> M/G$ to $G$-equivariant maps $P
@>>> M$ from a $G$-principal bundle $P$ over $\Sigma$ to $M$. The set
of equivalence classes of these lifts is denoted by $\Bbb
L_{\Sigma}(M;G)$. We also consider a subspace of equivalence classes
of locally constant $G$-equivariant maps $P @>>> M$:
$$
\Bbb L_{\Sigma}(M;G)\overset\text{def}\to=
\!\!\!\!\!\coprod_{[P @>>> \Sigma]} 
\!\!\!\!\!\text{Map}_G(P,M)/\text{Aut}_G(P)
\supset \Bbb L_{\Sigma}^{}(M;G)_{\text{l.c.}}=
\!\!\!\!\!\!\!\!\!\!\!\!\!\!\!\!\!\!\!
\coprod_{\ \ \ \ [\phi]\in\text{Hom}(\Gamma,G)/G}
\!\!\!\!\!\!\!\!\!\!\!\!\!\!\!\!\!\!\!
\bigl[M^{\langle\phi\rangle}/C(\phi)\bigr].
\tag1-2
$$
Here $[P @>>> \Sigma]$ runs over the set of all isomorphism classes of
$G$-principal bundles over $\Sigma$ whose fundamental group is
$\pi_1(\Sigma)=\Gamma$. Recall that isomorphism classes of $G$-bundles
over $\Sigma$ are classified by the $G$-conjugacy classes
$\text{Hom}(\Gamma,G)/G$ of homomorphisms. For any $\phi:\Gamma @>>>
G$, the subgroup $C(\phi)\subset G$ is the centralizer of the image of
$\phi$ denoted by $\text{Im}\,\phi=\langle\phi\rangle$, and
$M^{\langle\phi\rangle}$ is the fixed point subset of $M$ under
$\langle\phi\rangle$. We have $\text{Aut}_G(P)\cong C(\phi)$. Since
the space $\Bbb L_{\Sigma}(M;G)$ is the space of equivalence classes
of lifts, there is a natural map $\Bbb L_{\Sigma}(M;G) @>>>
\text{Map}(\Sigma, M/G)$ which is a homeomorphism when $G$-action on
$M$ is free. In this case, there is no need to go over the above
construction, and thus we are primarily interested in non-free
$G$-actions on $M$. Note that we have replaced the mapping space
$\text{Map}(\Sigma,M/G)$ on an orbifold $M/G$ by the space $\Bbb
L_{\Sigma}(M;G)$ which is a union of orbifolds on mapping spaces. As
such, the space $\Bbb L_{\Sigma}(M;G)$ is, in a sense, a mild
desingularization of the mapping space $\text{Map}(\Sigma,M/G)$.
Intuitively, our manifold $\Sigma$ is used to probe the nature of
orbifold singularities of the orbit space $M/G$ by examining the
holonomy of all possible $G$-bundles on $\Sigma$ ``induced'' by a map
$\gamma:\Sigma @>>> M/G$.

The space $\Bbb L_{\Sigma}(M;G)$ is motivated by and is a
generalization of the notion of {\it twisted sectors} in orbifold
conformal field theory, which corresponds to the case $\Sigma=S^1$.
Intuitively speaking, twisted sectors are vector spaces obtained by
``quantizing'' the mapping spaces $\text{Map}_G(P,M)$, where $P$ runs
over isomorphism classes of $G$-bundles over $S^1$. Thus the study of
global topology and geometry of these equivariant mapping space will
shed light on the nature of twisted sectors in orbifold conformal
field theory. Mathematically, twisted sectors have been studied in the
algebraic framework of vertex operator algebras and their modules. The
study of geometry and topology of the space $\Bbb L_{\Sigma}(M;G)$
when $\Sigma$ is $S^1$ or Riemann surfaces may provide illuminating
topological and geometric insight into the nature of the algebraic
structure of twisted sectors.

In this paper, we study the subspace of locally constant equivariant
maps instead of directly studying the global topology of $\Bbb
L_{\Sigma}(M;G)$, and we define our generalized orbifold invariants as
follows. Let $\varphi(M;G)$ be any (multiplicative) invariant of
$(M;G)$. For example, as $\varphi(M;G)$ we may use the following
invariants coming from the usual Euler characteristic:
$$
\chi^{\text{orb}}(M;G)\overset\text{def}\to= 
\frac{\chi(M)}{|G|}\in\frac1{|G|}\Bbb Z,
\qquad \qquad \chi(M;G)\overset\text{def}\to= 
\chi(M/G)\in\Bbb Z. 
\tag1-3
$$
The first Euler characteristic $\chi^{\text{orb}}$ is a special case
of the equivariant Euler characteristic defined for general (not
necessarily finite) group $G$ [\Cite{Br},\,p.249].  The $\Gamma$-{\it
extension} $\varphi_{\sssize\Gamma}(M;G)$ of $\varphi(M;G)$ is defined
in terms of the subspace of locally constant equivariant maps given in
(1-2), viewed equivariantly:
$$
\varphi_{\sssize\Gamma}(M;G)\overset{\text{def}}\to=
\!\!\!\!\!\!\!\!\!\!\!\!\!\!\!\!\!\!\!\!\!\!
\sum_{\ \ \ \ \ \ \ [\phi]\in\text{Hom}(\Gamma,G)/G}
\!\!\!\!\!\!\!\!\!\!\!\!\!\!\!\!\!\!\!\!\!
\varphi\bigl(M^{\langle\phi\rangle};C(\phi)\bigr).
\tag1-4
$$
This is our generalization of physicists' orbifold Euler
characteristic $e_{\sssize\text{orb}}(M;G)$. In our previous paper
\cite{T}, we considered $\chi_{\sssize\Gamma}(M;G)$ for the case
$\Gamma=\Bbb Z^d$ and $\Gamma=(\Bbb Z_p)^d$ with $d\in\Bbb N$, where
$p$ is any prime and $\Bbb Z_p$ is the ring of $p$-adic integers. In
the same paper, we described an application of the result for the case
$\Gamma=(\Bbb Z_p)^d$ to the Euler characteristic of Morava $K$-theory
of classifying spaces of wreath products. The inductive method used
for these calculations cannot be extended to a general group
$\Gamma$. The purpose of this paper is to develop a theory and
techniques which allow us to handle orbifold invariants
$\varphi_{\sssize\Gamma}(M;G)$ for a more general family of groups
$\Gamma$, including the fundamental groups of two dimensional surfaces
which may be orientable or non-orientable, compact or non-compact.

For more discussion on general properties of the $\Gamma$-extension of
$\varphi(M;G)$, see \S2. For example, we will show that
$\Gamma$-extensions of the two orbifold Euler characteristics in (1-3)
are related by
$$
\chi_{\sssize\Gamma}(M;G)=
\chi_{\sssize\Gamma\times\Bbb Z}^{\text{orb}}(M;G)\in\Bbb Z,
\tag1-5
$$
for any group $\Gamma$. Thus, $\chi_{\sssize\Gamma\times
\Bbb Z}^{\text{orb}}(M;G)$ is always integer valued.  When
$\Gamma=\Bbb Z$, the invariant $\chi_{\sssize\Bbb Z^2}^
{\text{orb}}(M;G)$ coincides with the orbifold Euler characteristic
(1-1) considered by physicists \cite{DHVW}, \cite{HH}. 
See the formula (2-4) 
in \S2. The above formula says that $\chi^{\text{orb}}(M;G)$ is more
basic than $\chi(M;G)$.

We demonstrate that orbifold invariants (1-4) is well behaved by
calculating the generating function of these invariants of symmetric
products of global quotients. Recall that an $n$-fold symmetric
product of an orbit space $M/G$ is given by 
$$
SP^n(M/G)=(M/G)^n\!/\frak S_n
=M^n\!/(G\!\wr\! \frak S_n),
\tag1-6
$$
where $G\!\wr\!\frak S_n$, which we also denote by $G_n$, is the wreath
product of $G$ and the symmetric group $\frak S_n$. For more details
on wreath products, see \S 3. 

Let $\Sigma$ be a real 2 dimensional surface, orientable or
non-orientable, compact or non-compact. For example, let $\Sigma$ be a
genus $g+1$ orientable closed surface. Its fundamental group
$\Gamma_{g+1}$ is generated by $2g+2$ elements with one relation:
$$
\Gamma_{\!g+1}=\langle a_1,a_2,\dots,a_{g+1},
b_1,b_2,\dots,b_{g+1} \mid [a_1,b_1]\cdots[a_{g+1},b_{g+1}]=1 
\rangle,
\tag1-7
$$
where $g\ge0$. Consideration of the twisted space (1-2) associated to
the mapping space $\text{Map}\,\bigl(\Sigma,SP^n(M/G)\bigr)$ from the
surface $\Sigma$ to symmetric products is motivated by string theory
literature \cite{DMVV}. We calculate the orbifold Euler characteristic
$\chi^{\text{orb}}_{\sssize \Gamma_{\!g+1}}(M^n;G_n)$ associated to
the surface group $\Gamma_{\!g+1}$. It turns out that this orbifold
Euler characteristic of an $n$-th symmetric product can be expressed
by $\{\chi^{\text{orb}}_{\sssize H}(M;G)\}_H$, where $H$ runs over
subgroups of $\Gamma_{\!g+1}$ of index at most $n$. In fact, we have
the following formula of the generating function [Theorem 5-8].

\proclaim{Theorem A}{\rm (Higher genus orbifold Euler characteristics of
symmetric orbifolds)}\qua Let $g\ge0$. With the above notations, 
$$
\sum_{n\ge0}q^n
\chi^{\text{orb}}_{\sssize\Gamma_{\!g+1}}(M^n;G\!\wr\!\frak S_n)
=\exp\Bigl[\ \sum_{r\ge1}q^r\Bigl\{\frac{j_r(\Gamma_{\!g+1})}{r}
\chi^{\text{orb}}_{\sssize\Gamma_{rg+1}}(M;G)\Bigr\}\Bigr],
\tag1-8
$$
where $j_r(\Gamma_{\!g+1})$ is the number of index $r$ subgroups of
$\Gamma_{\!g+1}$. 
\endproclaim
In this formula the genus $rg+1$ surface group $\Gamma_{rg+1}$ appears
since any index $r$ subgroup of $\Gamma_{\!g+1}$ is isomorphic to it,
although they may not be conjugate to each other in $\Gamma_{g+1}$. To
calculate numbers $j_r(\Gamma_{\!g+1})$ for $r\ge1$ and $g\ge0$, see
the formula (5-10).

Furthermore, we can also consider non-orientable cases. Let $\Sigma$ be
a closed genus $h+2$ non-orientable surface with $h\ge0$. Its
fundamental group $\Lambda_{h+2}$ is described by
$$
\Lambda_{h+2}=\langle c_1,c_2,\dots,c_{h+2} \mid 
c_1^2c_2^2\cdots c_{h+2}^2=1\rangle, \qquad h\ge0.
\tag1-9
$$
Since any genus $1$ non-orientable closed surface is homeomorphic to
$\Bbb RP^2$ and $\Lambda_1=\Bbb Z/2\Bbb Z$ is abelian, this case was
discussed in our previous paper \cite{T}. Here, we only consider
non-orientable surfaces of genus $2$ or higher. We have a formula of
the generating function of the orbifold invariant
$\chi^{\text{orb}}_{\sssize
\Lambda_{h+2}}(M^n;G\!\wr\!\frak S_n)$ for non-orientable surface groups
similar to the one in Theorem A. See Theorem 5-10 for details. 

This formula simplifies for the case of a genus $2$ non-orientable
surface which is a Klein bottle, due to a fact that every finite
covering of a Klein bottle is either a torus or a Klein bottle. Here
we give a formula for Klein bottle orbifold Euler characteristic
[Theorem 5-11].

\proclaim{Theorem B}{\rm (Klein bottle orbifold Euler characteristic)}\qua With
the above notation, 
$$
\sum_{n\ge0}q^n\chi^{\text{orb}}_{\sssize\Lambda_2}
(M^n;G\!\wr\!\frak S_n)
=\bigl[\prod_{r\ge1}(1-q^{2r})\bigr]^
{\frac{-1}2\chi^{\text{orb}}_{\Gamma_1}(M;G)}
\Bigl[\prod_{r\ge1}\Bigl(\frac{1+q^r}{1-q^r}\Bigr)\Bigr]^
{\frac12\chi^{\text{orb}}_{\Lambda_2}(M;G)}.
\tag1-10
$$
\endproclaim
Notice the appearance of modular forms, and modular functions. Let
$\eta(q)=q^{\frac1{24}}\prod_{r\ge1}(1-q^r)$, where $q=e^{2\pi i\tau}$
with $\text{Im}\,\tau>0$, be the Dedekind eta function. Then, the
first infinite product in (1-10) is almost $\eta(q^2)$ and the second
infinite product is precisely $\eta(q^2)/\eta(q)^2$.

By specializing $M=\text{pt}$, we obtain a formula for
$|\text{Hom}(\Lambda_2,G\!\wr\!\frak S_n)|$ for all $n$. See (5-19). 

Recall that the geometry behind the above Klein bottle orbifold Euler
characteristic is the mapping space $\text{Map}\bigl(K, 
SP^n(M/G)\bigr)$ from a Klein bottle $K$ to $n$-fold symmetric
product of an orbifold $M/G$ for various $n$. The idea of using
nonorientable surfaces comes from unoriented strings
moving on manifolds. 

There is a curious relationship between orbifold Euler characteristics
associated to orientable surface groups and non-orientable surface
groups. For example, when $G$ is a trivial group, we have
$\chi_{\Gamma_{\!g+1}}^{\text{orb}}(M^n;\frak S_n)
=\chi_{\Lambda_{2g+2}}^{\text{orb}}(M^n;\frak S_n)$ for any $n\ge1$
and $g\ge0$. For this, see the end of \S5. 

Note that the formula in Theorem A is in an exponential form. If we
use $\chi(M;G)=\chi(M/G)$ instead of $\chi^{\text{orb}}(M;G)$, then
the corresponding generating function can be written in an infinite
product form. To describe this formula, we need to introduce an
extension of $\varphi(M;G)$ to an {\it orbifold invariant
$\varphi_{\sssize[X]}(M;G)$ associated to an isomorphism class $[X]$
of a finite $\Gamma$-set $X$}. When $X$ is a transitive $\Gamma$ set
of the form $\Gamma\!/\!H$ for some subgroup $H\subset \Gamma$, we let
$$
\varphi_{\sssize[\Gamma\!/\!H]}(M;G)\overset\text{def}\to=
\!\!\!\!\!\!\!\!\!\!\!\!\!\!\!\!\!\!\!\!\!\!\!\!\!\!\!\!\!\!\!\!\!\!\!
\sum_{\ \ \ \ \ \ \ \ \ [\rho]\in\text{Hom}(H,G)/(N_{\Gamma}(H)\times G)}
\!\!\!\!\!\!\!\!\!\!\!\!\!\!\!\!\!\!\!\!\!\!\!\!\!\!\!\!\!\!\!\!\!\!\!
\varphi\bigl(M^{\langle\rho\rangle};
\text{Aut}_{\Gamma\text{-}G}(P_{\rho})\bigr),
\tag1-11
$$
where $N_{\sssize\Gamma}(H)\subset\Gamma$ is the normalizer of $H$ in
$\Gamma$ acting on $H$ by conjugation, and
$\text{Aut}_{\Gamma\text{-}G}(P_{\rho})$ is the group of
$\Gamma$-equivariant $G$-bundle automorphisms of a $G$-bundle
$P_{\rho}=\Gamma\times_{\rho} G$. The group $G$ of course acts on
$\text{Hom}(H,G)$ by conjugation. For an explanation and the origin of
this definition, see Proposition 6-1 and its proof. The above formula
reduces to (1-4) when $H=\Gamma$. For an explicit nontrivial example
of $\varphi_{\sssize[\Gamma\!/\!H]}(M;G)$, see Lemma 6-2.

The invariant $\varphi_{\sssize [X]}$ can be defined for a general
$\Gamma$-set $X$ (see (6-13) for the definition) in such a way that it
is multiplicative in $[X]$ in the following sense. Let $X_1$ and $X_2$
be $\Gamma$-sets without common transitive $\Gamma$-sets in their
$\Gamma$-orbit decompositions, that is,
$\text{Map}_{\Gamma}^{}(X_1,X_2)=\emptyset$. Then, the orbifold invariant
of their disjoint union $[X]=[X_1\amalg X_2]$ is a product:
$$
\varphi_{\sssize[X]}(M;G)=\varphi_{\sssize[X_1]}(M;G)\cdot 
\varphi_{\sssize[X_2]}(M;G).
\tag1-12
$$
For details, see Proposition 6-9. 

This invariant $\varphi_{\sssize[X]}(M;G)$ has a very natural
geometric origin in terms of a twisted space similar to (1-2) using
$G$-bundles over a covering space
$\Sigma'=\widetilde{\Sigma}\times_{\Gamma}X$ over $\Sigma$, where
$\widetilde{\Sigma}$ is the universal cover of $\Sigma$. See Theorem G
below. In fact, we could start with this approach and deduce the
formula (1-11).

The next theorem describes the generating function of the orbifold
Euler characteristic of symmetric orbifolds associated to the genus
$g+1$ orientable surface group $\Gamma_{\!g+1}$, in terms of the
orbifold Euler characteristic of $(M;G)$ associated to finite {\it
transitive} $\Gamma_{\!g+1}$-sets.

\proclaim{Theorem A${}'$} {\rm (Higher genus orbifold Euler
characteristic: infinite product form)}\qua For any $g\ge0$, we have
$$
\sum_{n\ge0}q^n\chi_{\sssize\Gamma_{\!g+1}}^{}(M^n;G\!\wr\!\frak S_n)
=\prod_{[H]}\big(1-q^{|\Gamma_{g+1}/H|}\bigr)
^{-\chi_{[\Gamma_{g+1}/H]}(M;G)},
\comment
=\topsmash{
\prod_{r\ge1}\Bigl\{(1-q^r)^
{-\!\!\!\!\!\!\!\!\!\!\!\!\!\!\!\!
\sum\limits\Sb {\sssize[H]} \\ 
       \ \ \ \ \sssize{|\Gamma_{\!g+1}\!/\!H|=r} \endSb
\!\!\!\!\!\!\!\!\!\!\!\!\!\!\!\!\!
\chi_{[\Gamma_{\!g+1}\!/\!H]}(M;G)}\Bigr\}.
}\endcomment
\tag1-13
$$
where $[H]$ runs overthe set of conjugacy classes of $H$ in
$\Gamma_{g+1}$. 
\endproclaim
This is part of Theorem 6-6. As remarked above, the invariant
$\chi_{\sssize[\Gamma_{\!g+1}\!/\!H]}(M;G)$ with
$|\Gamma_{\!g+1}\!/\!H|=r$ arises from a twisted space (1-21) defined
in terms of $G$-principal bundles over an $r$-fold covering space
$\Sigma_{rg+1}$ over $\Sigma_{g+1}$ corresponding to the conjugacy
class $[H]$.

To contrast the formulae in Theorem A and Theorem $\text{A}'$, we
let $M=\text{pt}$. We get the following combinatorial
group theoretic formulae:
$$
\aligned
\sum_{n\ge0}q^n\frac{|\text{Hom}(\Gamma_{\!g+1},G_n)|}
{|G|^n\cdot n!}&=\exp
\biggl[\sum_{r\ge1}\frac{q^r}r j_r(\Gamma_{\!g+1})
\frac{|\text{Hom}(\Gamma_{rg+1},G)|}{|G|}\biggr] \\
\sum_{n\ge0}q^n|\text{Hom}(\Gamma_{\!g+1},G_n)/G_n|
&=\prod_{[H]}\bigl(1-q^{|\Gamma_{g+1}/H|}\bigr)
^{-|\text{\rm Hom}(H,G)/(N_{\Gamma_{g+1}}(H)\times G)|}
\endaligned
\comment
&=\prod_{r\ge1}(1-q^r)^{
-\!\!\!\!\!\!\!\!\!\!\!\!\!\!\!
\sum\limits\Sb {\sssize[H]} \\ 
      \ \ \ \ {\sssize|\Gamma_{\!g+1}\!/\!H|=r} \endSb
\!\!\!\!\!\!\!\!\!\!\!\!\!\!\!
|\text{Hom}(H,G)/(N_{\Gamma_{\!g+1}}(H)\times G)|}.
\endaligned
\endcomment
\tag1-14
$$
In the second formula, $[H]$ runs over all conjugacy classes of finite
index subgroups of $\Gamma_{\!g+1}$. Note that the size of the
normalizer $N_{\Gamma_{\!g+1}}(H)$ may depend on the conjugacy class
$[H]$, although all index $r$ subgroups of $\Gamma_{\!g+1}$ are
isomorphic to $\Gamma_{\!rg+1}$. See also remarks after Theorem 6-6.

Here we remark that the following formulae for the number of
(conjugacy classes of) homomorphisms are well known, and they can be
easily proved using characters of $G$. For any finite
group $G$, we have
$$
\aligned
\frac{|\text{Hom}\,(\Gamma_{\!g+1},G)|}{|G|}
&=\!\!\!\!\!\!\!
\sum_{[V]\in\text{Irred}(G)}\!\!\!\!\!\!\!
\Bigl(\frac{|G|}{\dim V}\Bigr)^{\!2g},\\
|\text{Hom}(\Gamma_{\!g+1},G)/G|
&=\sum_{[\sigma]}\!\!\!\!\!\!\sum_{\ \ \ [V]\in\text{Irred}(C(\sigma))}
\!\!\!\!\!\!\!\!
\Bigl(\frac{|C(\sigma)|}{\dim V}\Bigr)^{\!2g},
\endaligned
\tag1-15
$$
where $[V]\in\text{Irred}(G)$ means that $[V]$ runs over the set of
all isomorphism classes of irreducible representations of $G$, and
$[\sigma]$ runs over the set of all conjugacy classes of $G$. The
subgroup $C(\sigma)\subset G$ is the centralizer of $\sigma$. Similar
formulae for non-orientable surface groups $\Lambda_{h+2}$ are well
known, too. See (5-16) and (5-17). In view of (1-15), to understand
the left hand side of (1-14), we need the representation theory of
wreath product $G_n=G\!\wr\!\frak S_n$. For this topic, see for
example [\Cite{M2},\,Appendix B to Chapter I]. However, to prove
topological Theorems A and A${}'$, we do not need the representation
theory of wreath products.

Two formulae in (1-14) are closely related, although their right hand
sides look very different. See (6-11) and (6-12) for group theoretic
discussion on this connection.

Theorems A,B, and A${}'$ are all special cases of the following more
general formulae which are the ones we actually prove [Theorem 5-5,
Theorem 6-3]. These are part of our main results.

\proclaim{Theorem C} With the above notations, for any group $\Gamma$,
we have 
$$
\aligned
\sum_{n\ge0}q^n
\chi^{\text{orb}}_{\sssize\Gamma}(M^n;G\!\wr\!\frak S_n)
&=\exp\Bigl[\sum_{H}\frac{q^{|\Gamma/H|}}{|\Gamma/H|}
\chi_H^{\text{\rm orb}}(M;G)\Bigr]\\
\sum_{n\ge0}q^n\chi_{\sssize\Gamma}(M^n;G\!\wr\!\frak S_n)
&=\prod_{[H]}\bigl(1-q^{|\Gamma/H|}\bigr)^{-\chi_{[\Gamma/H]}(M;G)}.
\endaligned
\tag1-16
$$
Here, in the first formula, $H$ runs over all subgroups of $\Gamma$
of finite index. In the second formula, $[H]$ runs over conjugacy
classes of all finite index subgroups of $\Gamma$. 
\endproclaim

By specializing $\Gamma$, we obtain numerous corollaries. In addition
to the cases $\Gamma=\Gamma_{\!g+1}$, $\Lambda_{h+2}$, we consider free
groups on finitely many generators, abelian groups $\Bbb Z^d$ and
$\Bbb Z^d_p$ for $d\ge1$, and products of these groups. 

Special cases of the second formula in (1-16) have been known. When
$\Gamma$ is the trivial group, this is Macdonald's formula
\cite{M1}. When $\Gamma=\Bbb Z$ and $G=\{e\}$, this formula is due to
\cite{HH}. For $\Gamma=\Bbb Z$ and general finite group $G$, the
formula is due to \cite{W}. When $\Gamma=\Bbb Z^d$ for $d\ge1$ and
$G=\{e\}$, the formula was proved in \cite{BF}. Finally, for
$\Gamma=\Bbb Z^d$ and $\Bbb Z^d_p$ with an arbitrary finite group $G$,
the formula was proved in \cite{T}.

When $G$ is the trivial group, we have
$\chi_{\sssize[X]}(M;\{e\})=\chi(M)$ for any transitive $\Gamma$-set
$X$ (see Lemma 6-10, and a formula after Corollary 6-4). When $\Gamma$
is an abelian group $\Bbb Z^d$ or $\Bbb Z^d_p$, Lemma 6-2 says that
$\chi_{\sssize[X]}(M;G)=\chi_{\sssize\Gamma}(M;G)$ for any finite
transitive $\Gamma$-set $X$. For general abelian group $\Gamma$, the
second formula in (1-16) reduces to (6-5). Thus, in all the known
cases described in the above paragraph, the notion of orbifold
invariants {\it associated to $\Gamma$-sets} is not yet needed. Only
in the generality of the present paper, this notion plays a crucial
role.

The main ingredient of the proof of Theorem C is a structure theorem
of the centralizer $C(\theta)\subset G_n$ of an arbitrary homomorphism
$\theta:\Gamma @>>> G_n=G\!\wr\!\frak S_n$ into a wreath product.
When $\Gamma=\Bbb Z$, the subgroup $C(\theta)$ is nothing but the
centralizer of the element $\theta(1)\in G_n$. Detailed description of
the centralizer of an element in $G_n$ was given in section 3 of
\cite{T}. The method used there was purely group theoretic. It is
very complicated to extend this group theoretic approach to the
present context of centralizers of homomorphisms into wreath
products. We have a better approach in terms of geometry of
$\Gamma$-equivariant $G$-principal bundles. These are $G$-principal
bundles $P$ over $\Gamma$-sets such that $\Gamma$ acts on $P$ as
$G$-bundle automorphisms. This link between algebra and geometry is
supplied by the following theorem [Theorem 3-1]. For the next three
theorems, $G$ does not have to be a finite group.

\proclaim{Theorem D} Let $G$ and $\Gamma$ be any group. Then the
following bijective correspondence exists\rom{:}  
$$
\biggl\{\foldedtext\foldedwidth{3in}
{Isomorphism classes of $\Gamma$-equivariant $G$-principal bundles
over a $\Gamma$-set of order $n$}
\biggr\}
\underset{\text{\rm onto}}\to{\overset{1:1}\to\longleftrightarrow}
\text{\rm Hom}(\Gamma,G_n)/G_n. 
\tag1-17
$$
\endproclaim

This theorem allows us to apply geometric concepts and techniques
associated to principal bundles to the study of centralizers of
homomorphisms $\theta$ into wreath products. Note that this theorem is
a generalization of a well-known bijective correspondence between the
set $\text{Hom}(\Gamma, \frak S_n)/\frak S_n$ of conjugacy classes of
homomorphisms into the $n$-th symmetric group $\frak S_n$ and the set
of isomorphism classes of $\Gamma$-sets of order $n$. As an example of
the above correspondence, let $H\subset\Gamma$ be any subgroup of
index $n$. For any homomorphism $\rho:H @>>> G$, let
$P_{\rho}=\Gamma\times_{\rho}G @>>> \Gamma\!/\!H$ be a
$\Gamma$-equivariant $G$-principal bundle over a $\Gamma$-set
$\Gamma\!/\!H$ of order $n$. By choosing a bijection
$\Gamma\!/\!H\cong\{1,2,\dots,n\}$ and a section of $P_{\rho} @>>>
\Gamma\!/\!H$, we have an isomorphism $\text{Aut}_G(P_{\rho})\cong
G\!\wr\!\frak S_n$ [Lemma 3-3]. Thus, the wreath product appears
naturally in the geometric context of $G$-bundles. This is the reason
of our use of $G$-bundles in studying wreath products. Now the
ready-made action of $\Gamma$ on $P_{\rho}$ as $G$-bundle maps gives
rise to a homomorphism $\theta:\Gamma @>>> G\!\wr\!\frak S_n$ whose
conjugacy class is independent of choices made. For more details on
this and the construction of the converse correspondence, see the
proof of Theorem 3-1.

We call a $\Gamma$-equivariant $G$-principal bundle over a
$\Gamma$-transitive set a $\Gamma$-{\it irreduc\-ible} $G$-principal
bundle. The following theorem classifies $\Gamma$-irreducible
$G$-principal bundles [Theorem 3-6].

\proclaim{Theorem E}{\rm (Classification of $\Gamma$-$G$ bundles)}\qua Let $G$
and $\Gamma$ be any groups. Then there exists the following bijective
correspondence\rom{:}
$$
\biggl\{\foldedtext\foldedwidth{2.4in}{Isomorphism classes of
$\Gamma$-irreducible $G$-principal bundles over
$\Gamma$-sets of order $n$}\biggr\} 
\underset{\text{\rm onto}}\to{\overset{1:1}\to\longleftrightarrow}
\!\!\!\!\!\!\!\!\!\coprod\Sb[H] \\ \ \ \ |\Gamma\!/\!H|=n \endSb
\!\!\!\!\!\!\!\!\!\text{\rm Hom}(H,G)/(N_{\Gamma}(H)\times G),
\tag1-18
$$
where $[H]$ runs over the set of conjugacy classes of index $n$ 
subgroups of $\Gamma$.
\endproclaim

Note that the expression appearing on the right hand side of the above
correspondence has already appeared in (1-11). From the point of view
of a $\Gamma$-irreducible $G$-principal bundle $P @>>> Z$, the meaning
of quantities $[H]$, $N_{\Gamma}^{}(H)$-action, and $G$-action on
$\text{Hom}(H,G)$, appearing in the right hand side of (1-18), is as
follows: the isomorphism class of the transitive $\Gamma$-set $Z$ is
determined by the conjugacy class $[H]$ of an isotropy subgroup $H$,
different choices of a base point $z_0$ of $Z$ with the isotropy
subgroup $H$ correspond to the action of $N_{\Gamma}(H)$, and
different choices of a base point $p_0$ of $P$ over $z_0$ correspond
to the conjugation action of $G$ on $\text{Hom}(H,G)$.

For any homomorphism $\theta:\Gamma @>>> G_n$, the corresponding
$\Gamma$-equivariant $G$-bundle $\pi_{\theta}: P_{\theta} @>>>
Z_{\theta}$ is a disjoint union of $\Gamma$-irreducible $G$-bundles
each of which is of the form $P_{\rho} @>>>\Gamma\!/\!H$ for some
subgroup $H$ and some homomorphism $\rho:H @>>> G$. Given $\theta$,
let $r(H,\rho)$ denote the number of $\Gamma$-irreducible $G$-bundles
isomorphic to $[P_{\rho} @>>> \Gamma\!/\!H]$ appearing in the
irreducible decomposition of $P_{\theta} @>>> Z_{\theta}$. The number
$r(H,\rho)$ depends only on the conjugacy classes $[H]$ and
$[\rho]$. Now we can describe the centralizer of the image of $\theta$
in the wreath product $G\!\wr\!\frak S_n$ [Lemma 4-1, Theorem 4-2,
Theorem 4-4, Theorem 4-5].

\proclaim{Theorem F}{\rm (Structure of centralizer of homomorphisms into
wreath products)}\qua Let $\theta:\Gamma @>>> G_n$ be a homomorphism, and
let $\{r(H,\rho)\}_{[H],[\rho]}$ be the associated integers described 
above. Then, the centralizer of $\theta$ in the wreath product
$G\!\wr\!\frak S_n$ is a direct product of wreath products\rom{:}
$$
C_{G_n}(\theta)\cong \prod_{[H]}\prod_{[\rho]}
\bigl\{\text{\rm Aut}_{\Gamma\text{-}G}(P_{\rho})
\!\wr\!\frak S_{r(H,\rho)}\bigr\},
\tag1-19
$$
where $[H]$ runs over conjugacy classes of finite index subgroups, and
$[\rho]\in \text{\rm Hom}(H,G)/(N_{\Gamma}(H)\times G)$ for a given
conjugacy class $[H]$. The group $\text{\rm Aut}_{\Gamma\text{-}G}
(P_{\rho})$ of $\Gamma$-equivariant $G$-bundle automorphisms of
$P_{\rho}$ fits into the following exact sequences\rom{:}
$$
\gathered
1 @>>> C_G(\rho) @>>> \text{\rm Aut}_{\Gamma\text{-}G}(P_{\rho}) 
@>>> H\backslash N_{\Gamma}^{\rho}(H) @>>> 1, \\
1 @>>> H @>>> T_{\rho} @>>> \text{\rm Aut}_{\Gamma\text{-}G}(P_{\rho})
@>>> 1,
\endgathered
\tag1-20
$$
where $N_{\Gamma}^{\rho}(H)$ is the isotropy subgroup at the
$G$-conjugacy class $(\rho)$ of the $N_{\Gamma}(H)$-action on the set
$\text{\rm Hom}(H,G)/G$, and $T_{\rho}$ is given in \rom{(4-4)}. 
\endproclaim

Here, the centralizer $C_{G_n}(\theta)$ is isomorphic to
$\text{Aut}_{\sssize\Gamma-G}(P_{\theta})$ [Lemma 4-1]. 

This structure theorem can be applied to much wider context including,
for example, the calculation of generalized orbifold elliptic genus of
symmetric orbifolds.

Next, we describe the geometric meaning of the orbifold invariant (1-11)
associated to $\Gamma$-sets. Recall that the orbifold invariant
associated to the group $\Gamma$ is defined in terms of locally
constant $G$-equivariant maps (1-2) from $G$-principal bundles over
$\Sigma$ into $M$, where the fundamental group of $\Sigma$ is
$\Gamma$. We can generalize this idea and consider locally constant
$G$-equivariant maps from $G$-principal bundles over {\it covering
spaces} of $\Sigma$ into $M$. More precisely, let $\pi:\Sigma' @>>>
\Sigma$ be a not necessarily connected covering space. Let $\pi_i: P_i
@>>> \Sigma'$ for $i=1,2$ be $G$-bundles. A {\it $G$-bundle
isomorphism $\alpha: P_1 @>>> P_2$ over a covering space} $\Sigma'$
are defined to be a $G$-bundle map whose induced map on the base
$\Sigma'$ is a deck transformation of $\Sigma'$ over $\Sigma$. Next,
we introduce an equivalence relation among $G$-maps $P @>>> M$. Two
$G$-equivariant maps $\gamma_1:P_1 @>>> M$ and $\gamma_2:P_2 @>>> M$
are said to be equivalent if there exists a $G$-bundle isomorphism
$\alpha: P_1 @>>> P_2$ over the covering space $\Sigma'$ such that
$\gamma_1=\gamma_2\circ\alpha$. We then consider the space of
equivalence classes of $G$-equivariant maps from $G$-bundles $P$ over
$\Sigma'$ into $M$. This turns out to be
$$
\Bbb L_{\Sigma'\!\!/\Sigma}^{}(M;G)
=\!\!\!\!\!\!\!\!\coprod_{[P @>>> \Sigma'\!\!/\Sigma]}
\!\!\!\!\!\!\!
\bigl[\text{Map}_G^{}(P,M)/\text{Aut}_G^{}
(P)_{\Sigma'\!\!/\Sigma}\bigr],
\tag1-21
$$
where $[P @>>> \Sigma'\!\!/\Sigma]$ denotes an isomorphism classes of
$G$-principal bundles over the covering space $\Sigma'$ over $\Sigma$. We
then take the subset $\Bbb L_{\Sigma'\!\!/\Sigma}^{}
(M;G)_{\text{l.c.}}$ of locally constant $G$-equivariant maps and
define the {\it orbifold invariant associated to a covering space}
$\Sigma' @>>> \Sigma$ by
$$
\varphi_{\sssize[\Sigma'\!\!/\Sigma]}^{}(M;G)
=
\!\!\!\!\!\!\!
\sum_{[P @>>> \Sigma'\!\!/\Sigma]}
\!\!\!\!\!\!\!
\varphi\bigl(\text{Map}_G^{}(P,M)_{\text{l.c.}}^{};
\text{Aut}_G^{}(P)_{\Sigma'\!\!/\Sigma}\bigr). 
\tag1-22
$$
The next result [Theorem 7-1] clarifies the geometric meaning of orbifold
Euler characteristic associated to $\Gamma$-sets. See (6-13) for the 
precise definition. 

\proclaim{Theorem G} Let $\pi : \Sigma' @>>> \Sigma$ be a covering
space, and let $X=\pi^{-1}(x_0)$ be a fibre over a base point
$x_0\in\Sigma$. Let $\Gamma=\pi_1(\Sigma,x_0)$. Then $X$ is a
$\Gamma$-set and for any $G$-manifold $M$, we have 
$$
\varphi_{\sssize[\Sigma'\!\!/\Sigma]}^{}(M;G)=
\varphi_{\sssize[X]}^{}(M;G).
\tag1-23
$$
Consequently, $\Gamma$-orbifold Euler characteristic of symmetric
products of $G$-manifold $M$ can be expressed in terms of orbifold
Euler characteristic associated to connected covering spaces\rom{:}
$$
\sum_{n\ge0}q^n\chi_{\sssize\Gamma}(M^n;G\!\wr\!\frak S_n)
=\!\!\!\!\!\!\!\!\!\!\!\!\!\!
\prod_{\ \ \ [\Sigma' @>>> \Sigma]_{\text{conn.}}}
\!\!\!\!\!\!\!\!\!\!\!\!\!
(1-q^{\sssize|\Sigma'\!\!/\Sigma|})
^{-\chi_{\sssize[\Sigma'\!\!/\Sigma]}(M;G)},
\tag1-24
$$
where the product range over all isomorphism classes of finite
connected covering spaces $\Sigma' @>>> \Sigma$ over $\Sigma$, and
$|\Sigma'\!\!/\Sigma|$ denotes the order of the covering.
\endproclaim
Note that the definition (1-22) is geometrically very natural and
conceptually very simple, compared with the practical definition
(1-11) of orbifold invariants associated to transitive $\Gamma$-sets
which is suitable for calculation.

The above formula is very interesting because it explains the
geometric origin of the infinite product (one factor for each
connected covering space) and describes each factor geometrically in
terms of covering spaces.

Intuitive explanation of this identity may be given as follows. There
are two ways to lift maps in $\text{Map}\bigl(\Sigma, SP^n(M/G)\bigr)$
to equivariant maps from certain principal bundles. One way is to lift
a map $\overline{\gamma}: \Sigma @>>> SP^n(M\!/\!G)$ to a $G\!\wr\!\frak
S_n$-equivariant map $\gamma:P_{\!\sssize G\wr\frak S_n} @>>> M^n$,
where $P_{\!\sssize G\wr\frak S_n}$ is a $G\!\wr\!\frak S_n$-principal
bundle over $\Sigma$. The orbifold invariant using this lifting gives
rise to the left hand side of (1-24).

To describe another way of lifting, note that any map $f:\Sigma @>>>
SP^n(M)$ into an $n$-th symmetric product corresponds to a (branched)
$n$-sheeted covering space $\pi: \Sigma' @>>> \Sigma$ whose fibre over
$x\in\Sigma$ is $f(x)\subset M$. That is, $\Sigma'$ is given by 
$\Sigma'=\{(x,y)\in\Sigma\times M
\mid y\in f(x)\}$. This is not always a covering space because for
some $x\in\Sigma$ the corresponding un-ordered set of $n$ points
$f(x)$ may not consist of $n$ distinct points. The (branched) covering
space $\Sigma'$ comes equipped with a natural map $f':\Sigma' @>>> M$
which is unique up to ``deck transformations'' of $\Sigma'$. Namely,
if we have a homeomorphism $h:\Sigma' @>>> \Sigma'$ commuting with
$\pi$, then both maps $f'$ and $f'\circ h$ induce the same map
$f:\Sigma @>>> SP^n(M)$ by associating to each $x\in\Sigma$ the fibres
of $f'$ and $f'\circ h$. 

Now in our context, to any map $\overline{\gamma}:\Sigma @>>>
SP^n(M\!/\!G)$ we cab associate a map $\overline{\gamma}':\Sigma' @>>>
M\!/\!G$ from an $n$-sheeted (branched) covering space $\Sigma'$ over
$\Sigma$ described above. This map $\overline{\gamma}'$ is determined
up to ``deck transformations'' of $\Sigma'$. We then consider lifting
this map to a $G$-equivariant map $\gamma': P_{\sssize G} @>>> M$,
where $P_{\sssize G}$ is a $G$-principal bundle over
$\Sigma'$. Orbifold invariants corresponding to these lifts are the
orbifold invariants associated to covering spaces in the right hand
side of (1-24). We expect the covering space formalism in (1-24) would
be valid for other orbifold invariants and also for the global
topology of twisted sectors.

We suspect that the formula (1-24) will continue to be valid for
general orbifolds which are not necessarily global quotients. We hope
to come back to this question later. 

Finally, as an application to combinatorial group theory, we calculate
the number $u_r(\Gamma)$ of conjugacy classes of index $r$ subgroups
of $\Gamma$. Our combinatorial formulae obtained by specializing our
topological formulae to the case $M=\text{pt}$ allow us to prove the
following formula on $u_r(\Gamma)$ [Theorem 8-1]. This simple formula
does not seem to be known before. 

\proclaim{Theorem H} The number $u_r(\Gamma)$ of conjugacy classes of
index $r$ subgroups of $\Gamma$ satisfies the following recursive
relation in terms of subgroups of $\Gamma$\rom{:}
$$
j_m(\Gamma\times\Bbb Z)=\sum_{r|m}r\cdot u_r(\Gamma)
=\botsmash{
\sum_{r|m}\!\!\!\!\!\!\!\!\!\!\!\!\!
\sum\Sb H \\\ \ \ \ \ \ \ |\Gamma\!/\!H|=m/r \endSb
\!\!\!\!\!\!\!\!\!\!\!\!\!
|\text{\rm Hom}(H_{\text{\rm ab}},\Bbb Z_r)|,
}
\tag1-25
$$
where $H_{\text{\rm ab}}$ denotes the abelianization of $H$. 
\endproclaim
If we like, instead of having a recursive formula as above, we can
write down the formula for $u_r(\Gamma)$ using M\"obius inversion
formula. See the paragraph after Theorem 8-1. The above formula allows
us to compute numbers $u_r(\Gamma)$ for surface groups
$\Gamma_{\!g+1}$, $F_{s+1}$, and $\Lambda_{h+2}$ for $g,s,h\ge0$ very
quickly [Corollary 8-2]. Geometrically these numbers $u_r(\Gamma)$ are
important because they count the number of non-isomorphic connected
covering spaces over surfaces. These numbers were calculated before in
\cite{Me}, \cite{MP} by different methods. However, our formula (1-25) is
valid for {\it any} group $\Gamma$, not just for surface groups. In
our formula, all we have to do is to count the number of homomorphisms
from abelianization of subgroups to cyclic groups. This simplicity
makes our formula useful.

The organization of this paper is as follows. In \S2, we discuss
general properties of $\Gamma$-extended orbifold invariants and
explain their geometric origin. In \S3, we study homomorphisms into
wreath products in terms of $\Gamma$-equivariant $G$-bundles and prove
Theorems D and E. In \S4, we study the structure of centralizers of
homomorphisms into wreath products in terms of $\Gamma$-equivariant
$G$-automorphisms and prove Theorem F. In \S5, we prove Theorems A and
B, together with the first part of Theorem C. We also discuss the
cases when $\Gamma$ is a free abelian group, the fundamental group of
a non-orientable surface, and a free group. In \S6, we further extend
orbifold invariants to ones associated to isomorphism classes of
$\Gamma$-sets. This will be used to prove Theorem A${}'$ and Theorem
C.  Our main results in this paper deal with much more general cases
than \cite{T}, and results in \cite{T} follow quickly as special
cases. In \S7, we describe orbifold invariants associated to covering
spaces and prove Theorem G. Finally, in \S8, we apply our
combinatorial formula to compute the numbers $u_r(\Gamma)$ and prove
Theorem H.

\head
Orbifold invariants associated to a group $\Gamma$
\endhead

In this section, we describe general properties of the
$\Gamma$-extended orbifold invariant $\varphi_{\sssize\Gamma}(M;G)$
defined in (1-4). In particular, properties of the $\Gamma$-extensions
of orbifold Euler characteristics defined in (1-3) are discussed in
detail. Later in this section, we explain twisted spaces in the
context of mapping spaces into orbifolds and motivate the definition
of our generalized orbifold invariant (1-4).

Let $G$ be a finite group and let $\varphi(M;G)$ be a multiplicative
invariant for a $G$-manifold $M$. Namely, for any $G_1$-manifold $M_1$
and $G_2$-manifold $M_2$, the invariant $\varphi$ satisfies 
$$
\varphi(M_1\times M_2;G_1\times G_2)
=\varphi(M_1;G_1)\varphi(M_2;G_2). 
$$
This invariant $\varphi(M;G)$ may depend only on the topology and
geometry of the orbit space $M/G$. For example, Euler characteristic
$\chi(M/G)$, signature $\text{Sgn}(M/G)$, and spin index
$\text{Spin}(M/G)$, in appropriate geometric settings, are such
examples. On the other hand, the invariant $\varphi(M;G)$ may actually
depend on the way $G$ acts on $M$. For example, for any group
$\Gamma$, $\Gamma$-extensions of the above geometric invariants are
such examples. For convenience, we repeat the definition of
$\Gamma$-extension here. For any such invariant $\varphi(\ \cdot\ ;\
\cdot\ )$, its $\Gamma$-extension is defined by
$$
\topsmash{
\varphi_{\sssize\Gamma}(M;G)=
\!\!\!\!\!\!\!\!\!\!\!\!\!\!\!\!\!\!\!\!\!\!
\sum_{\ \ \ \ \ \ \ [\phi]\in\text{Hom}(\Gamma,G)/G}
\!\!\!\!\!\!\!\!\!\!\!\!\!\!\!\!\!\!\!\!\!
\varphi\bigl(M^{\langle\phi\rangle};C(\phi)\bigr),
}
\tag2-1
$$
where $M^{\langle\phi\rangle}$ is the
$\langle\phi\rangle=\text{Im}\,\phi$ fixed point subset, and $C(\phi)$
is the centralizer of $\langle\phi\rangle$ in $G$.  Note that it can
happen that two non-conjugate homomorphisms $\phi_1,\phi_2:\Gamma @>>>
G$ have the same image subgroups
$\langle\phi_1\rangle=\langle\phi_2\rangle$. Also note that the
quantity $\varphi\bigl(M^{\langle\phi\rangle};C(\phi)\bigr)$ depends
only on the $G$-conjugacy class of $\phi$.

Basic properties of $\Gamma$-extended orbifold invariants are the
following ones. 

\proclaim{Proposition 2-1} \!Let $\varphi(M;G)$ be a multiplicative
invariant as above. 

\rom{(1)}\qua For any group $\Gamma$, the $\Gamma$-extended orbifold
invariant $\varphi_{\sssize \Gamma}(M;G)$ is multiplicative\rom{:} for
any $(M_1;G_1)$ and $(M_2;G_2)$, we have
$$
\varphi_{\sssize\Gamma}(M_1\times M_2;G_1\times G_2)
=\varphi_{\sssize\Gamma}(M_1;G_1)
\varphi_{\sssize\Gamma}(M_2;G_2).
$$

\rom{(2)}\qua For any groups $K,L$, 
$$
\varphi_{\sssize K\times L}(M;G)
=(\varphi_{\sssize L})_{\sssize K}(M;G)
=\smash{
\!\!\!\!\!\!\!\!\!\!\!\!\!\!\!\!\!\!\!
\sum_{\ \ \ \ [\phi]\in\text{\rm Hom}(K,G)/G}
\!\!\!\!\!\!\!\!\!\!\!\!\!\!\!\!\!\!\!\varphi_{\sssize L}
\bigl(M^{\langle\phi\rangle};C(\phi)\bigr). 
}
\tag2-2
$$
\endproclaim
\demo{Proof} The formula in (1) is straightforward using the
multiplicativity of $\varphi$. For (2), by definition,
$$
\align
\text{(R.H.S.)}
&=\!\!\!\!\sum_{\phi_1:K @>>> G}\!\!\frac1{\#[\phi_1]}
\Bigl\{\!\!\!\sum_{\phi_2:L @>>> C(\phi_1)}
\!\!\frac1{\#[\phi_2]_{\sssize C(\phi_1)}}\varphi
\bigl((M^{\langle \phi_1\rangle})^{\langle\phi_2\rangle};
C_{C(\phi_1)}(\phi_2)\bigr)\Bigr\} \\
\intertext{Here, $\#[\phi_1]$ denotes the number of elements in the
conjugacy class $[\phi_1]$, and $[\phi_2]_{C(\phi_1)}$ is the
$C(\phi_1)$-conjugacy class of $\phi_2$. Since the images of $\phi_1$
and $\phi_2$ commute, using the notation $\phi_1\times\phi_2:K\times L
@>>> G$, we have} 
&=\!\!\!\!\!\!\!\!\!\!\!\!\!
\sum_{\phi_1:K @>>> G\ \ \ \ \ }\!\!\!\!\!\!\!\!\!\!\!
\frac1{\#[\phi_1]}\!\!\!\!\!\!\!\!
\sum_{\ \ \ \phi_2:L @>>> C(\phi_1)}
\!\!\!\!\!\!\!\!\!\frac{|C(\phi_1\times\phi_2)|}{|C(\phi_1)|}
\varphi\bigl(M^{\langle\phi_1\times\phi_2\rangle};
C(\phi_1\times\phi_2)\bigr) \\
&=\!\!\!\!\!\!\!\!\!\!\!\!\!\!\!\!\!\!\!\!
\sum_{\phi_1\times\phi_2:K\times L@>>> G\ \ \ \ \ }
\!\!\!\!\!\!\!\!\!\!\!\!\!\!\!\!\!\!\!
\frac{|C(\phi_1\times\phi_2)|}{|G|}
\varphi\bigl(M^{\langle\phi_1\times\phi_2\rangle};
C(\phi_1\times\phi_2)\bigr) \\
&=\!\!\!\!\!\sum_{[\phi_1\times\phi_2]}
\!\!\!\!\!\varphi\bigl(M^{\langle\phi_1\times\phi_2\rangle};
C(\phi_1\times\phi_2)\bigr)
=\varphi_{\sssize K\times L}(M;G).
\endalign
$$
In the last line, $[\phi_1\times\phi_2]$ runs over all the conjugacy
classes $\text{Hom}(K\times L,G)/G$, and we used a fact that the
number of elements in the conjugacy class of the image of
$\phi_1\times\phi_2$ is $|G|/|C(\phi_1\times\phi_2)|$. In the above,
we used a formula $|G|=\#[\phi]\cdot |C(\phi)|$ a few times. This
completes the proof.
\qed
\enddemo
 
Recall that in section 1, we defined two kinds of Euler
characteristics for $(M;G)$ in (1-3). One of them is
$\chi^{\text{orb}}(M;G)=\chi(M)/|G|$, and the other is
$\chi(M;G)=\chi(M/G)$. It is well known that the Euler characteristic
of the orbit space $M/G$ of a $G$-manifold $M$ can be calculated as
the average of the Euler characteristic of $g$-fixed point subsets (see for
example [\Cite{Sh},\,p.127]):
$$
\chi(M/G)=\frac1{|G|}\sum_{g\in G}\chi(M^{\langle g\rangle}).
\tag2-3
$$
Here are some useful formulae for $\Gamma$-extended orbifold Euler
characteristics.

\proclaim{Proposition 2-2} \rom{(1)}\qua The invariant
$\chi^{\text{orb}}_{\sssize\Gamma}(M;G)$ can be calculated as an
average over homomorphisms $\phi:\Gamma @>>> G$\rom{:}
$$
\chi^{\text{orb}}_{\sssize\Gamma}(M;G)
=\frac1{|G|}\!\!\!\!\!\!\!\sum_{\ \ \ \phi:\Gamma @>>> G}
\!\!\!\!\!\!\!\chi\bigl(M^{\langle\phi\rangle}\bigr).
\tag2-4
$$
In particular, $\chi^{\text{\rm orb}}_{\sssize\Bbb Z}(M;G)
=\chi(M/G)$. 

\rom{(2)}\qua The $\Gamma$-extensions of $\chi^{\text{\rm orb}}(M;G)$ and
$\chi(M;G)$ are related as follows\rom{:}
$$
\chi^{\text{\rm orb}}_{\sssize\Gamma\times\Bbb Z}(M;G)
=\chi_{\sssize\Gamma}(M;G)
=\!\!\!\!\!\!\!\!\!\!\!\!\!\!\!\!\!\!\!\!\!
\sum_{\ \ \ \ \ \ [\phi]\in\text{\rm Hom}(\Gamma,G)/G}
\!\!\!\!\!\!\!\!\!\!\!\!\!\!\!\!\!\!\!\!\!
\chi\bigl(M^{\langle\phi\rangle}/C(\phi)\bigr). 
\tag2-5
$$
\endproclaim
\demo{Proof} For part (1), we unravel definitions. 
$$
\chi^{\text{orb}}_{\sssize\Gamma}(M;G)
=\sum_{[\phi]}\chi^{\text{orb}}
\bigl(M^{\langle\phi\rangle};C(\phi)\bigr)
=\sum_{\phi}\frac1{\#[\phi]}\frac{\chi(M^{\langle\phi\rangle})}
{|C(\phi)|}
=\frac1{|G|}\sum_{\phi}\chi(M^{\langle\phi\rangle}), 
$$
where $[\phi]\in\text{Hom}(\Gamma,G)/G$ and
$\phi\in\text{Hom}(\Gamma,G)$. This proves formula (2-4), 
In particular, when $\Gamma=\Bbb Z$, we have 
$\chi^{\text{orb}}_{\sssize\Bbb Z}(M;G)=(1/|G|)\sum_{g\in G}
\chi(M^{\langle g\rangle})$, which is equal to $\chi(M/G)$ by (2-3). 

For (2), by the product formula (2-2), we have 
$$
\chi^{\text{orb}}_{\sssize\Gamma\times\Bbb Z}(M;G)
=\!\!\!\!\!\!\!\!\!\!\!\!\!\!\!\!\!\!\!
\sum_{\ \ \ \ \ \ [\phi]\in\text{Hom}(\Gamma,G)/G}
\!\!\!\!\!\!\!\!\!\!\!\!\!\!\!\!\!\!\!
\chi^{\text{orb}}_{\sssize\Bbb Z}
\bigl(M^{\langle \phi\rangle};C(\phi)\bigr)
=\sum_{[\phi]}\chi\bigl(M^{\langle\phi\rangle}/C(\phi)\bigr)
=\chi_{\sssize\Gamma}(M;G).
$$
This completes the proof. 
\qed
\enddemo

For example, when $\Gamma=\Bbb Z^2$, we have 
$$
\chi^{\text{orb}}_{\sssize\Bbb Z^2}(M;G)
=\frac1{|G|}\!\!\!\!\!\sum_{\ \ gh=hg}\!\!\!\!\!
\chi(M^{\langle g,h\rangle})
=\sum_{[g]}\chi\bigl(M^{\langle g\rangle}/C(g)\bigr),
$$
which is the physicist's orbifold Euler characteristic (1-1). The
following corollary is obvious. 

\proclaim{Corollary 2-3} When $M=\text{\rm pt}$, we have 
$$
\chi^{\text{\rm orb}}_{\sssize\Gamma}(\text{\rm pt};G)
=\frac{|\text{\rm Hom}(\Gamma,G)|}{|G|},\quad
\chi_{\sssize\Gamma}(\text{\rm pt};G)=|\text{\rm Hom}(\Gamma,G)/G|
=\frac{|\text{\rm Hom}(\Gamma\times\Bbb Z,G)|}{|G|}.
\tag2-6
$$
When $G$ is the trivial group $\{e\}$, we have 
$$
\chi^{\text{\rm orb}}_{\sssize\Gamma}(M;\{e\})
=\chi_{\sssize\Gamma}(M;\{e\})=\chi(M). 
\tag2-7
$$
\endproclaim

Here, the identity $|\text{Hom}(\Gamma,G)/G|=
|\text{Hom}(\Gamma\times\Bbb Z,G)|/|G|$ comes from 
$\chi_{\sssize\Gamma}(\text{pt};G)$ $=\,
\chi_{\sssize\Gamma\times\Bbb Z}^{\text{orb}}(\text{pt};G)$ given in
(2-5). This is a well-known and useful formula.

When $\chi^{\text{orb}}_{\sssize\Gamma}(\ \cdot\ ;G)$ is regarded as a
function on $G$-CW complexes with values in $\frac1{|G|}\Bbb Z$, it is
an additive function. If either $G$ or $\Gamma$ is abelian, it is a
complex oriented additive function in the sense of \cite{HKR}. Since
an additive function defined on $G$-CW complexes is completely
determined by its values on transitive $G$-sets, we calculate these
values. We also do the same for $\chi_{\Gamma}^{}(\ \cdot\ ;G)$.

\proclaim{Proposition 2-4} Let $H$ be a subgroup of $G$. Then 
$$
\aligned
\chi^{\text{\rm orb}}_{\sssize\Gamma}(G/H;G)
&=\chi_{\sssize\Gamma}^{\text{\rm orb}}(\text{\rm pt};H)
=\topsmash{
\frac{|\text{\rm Hom}(\Gamma,H)|}{|H|}
} \\
\chi_{\sssize\Gamma}(G/H;G)
&=\chi_{\sssize\Gamma}(\text{\rm pt};H)
=|\text{\rm Hom}(\Gamma,H)/H|.
\endaligned
\tag2-8
$$
\endproclaim
\demo{Proof} Using (2-4), we have 
$$
\chi^{\text{orb}}_{\sssize\Gamma}(G/H;G)
=\frac1{|G|}\!\!\!\!\!\!\sum_{\ \ \phi:\Gamma @>>> G}
\!\!\!\!\!\!\chi\bigl((G/H)^{\langle\phi\rangle}\bigr) 
=\frac1{|G|}\sum_{\phi}|(G/H)^{\langle\phi\rangle}|.
$$
Here $|(G/H)^{\langle\phi\rangle}|$ counts the number of elements in
$G/H$ whose isotropy subgroups contain $\phi(\Gamma)$. Since there are
$|N_G(H)/H|$ elements in $G/H$ with the same isotropy subgroups, we
have 
$$
|(G/H)^{\langle\phi\rangle}|
=\!\!\!\!\sum\Sb\phi(\Gamma)\subset H' \\ H'\sim H \endSb
\!\!\!\!|N_G(H)/H|.
$$
Here $H'\sim H$ means that $H'$ and $H$ are conjugate in $G$. 
Continuing our calculation, we have 
$$
\align
\frac1{|G|}\sum_{\phi}|(G/H)^{\langle\phi\rangle}|
&=\frac{|N_G(H)/H|}{|G|}
\sum_{\phi}\!\!
\sum\Sb \phi(\Gamma)\subset H' \\ H'\sim H \endSb \!\!\!\!\!1
=\frac{|N_G(H)|}{|G|\,|H|}\sum_{H'\sim H}
\sum_{\phi:\Gamma@>>> H'}\!\!\!\!\!1\\
&=\frac{|N_G(H)|}{|G|\,|H|}|\text{Hom}(\Gamma,H)|
\sum_{H'\sim H}\!\!1. 
\endalign
$$
since the number of subgroups $H'$ conjugate to $H$ in $G$ is given by
$|G|/|N_G(H)|$, the above quantity is equal to
$|\text{Hom}(\Gamma,H)|/|H|$.  

For $\chi_{\sssize\Gamma}(\ \cdot\ ;G)$, using (2-5), (2-6), and what
we have already proved, we have
$$
\chi_{\sssize\Gamma}(G/H;G)
=\chi_{\sssize\Gamma\times\Bbb Z}^{\text{orb}}(G/H;G)
=\chi_{\sssize\Gamma\times\Bbb Z}^{\text{orb}}(\text{pt};H)
=\chi_{\sssize\Gamma}(\text{pt};H)=|\text{Hom}(\Gamma,H)/H|.
$$ 
This completes the proof.
\qed
\enddemo

Next, we turn our attention to orbifold Euler characteristics of
symmetric products without $\Gamma$-extension. Macdonald's formula
\cite{M1} says that 
$$
\sum_{n\ge0}q^n\chi(M^n;G\!\wr\!\frak S_n)=
\sum_{n\ge0}q^n\chi\bigl(SP^n(M\!/\!G)\bigr)
=\frac1{(1-q)^{\chi(M/G)}}.
\tag2-9
$$
The corresponding formula for
$\chi^{\text{orb}}(\ \cdot\ ;\ \cdot\ )$ is given by
$$
\sum_{n\ge0}q^n\chi^{\text{orb}}(M^n;G\!\wr\!\frak S_n)
=\exp\Bigl\{q\frac{\chi(M)}{|G|}\Bigr\}
=\exp\{q\,\chi^{\text{orb}}(M;G)\}.
\tag2-10
$$
This formula follows immediately since
$\chi^{\text{orb}}(M_n;G\!\wr\!\frak S_n) =\chi(M)^n/(|G|^nn!)$ by
definition. The objective of this paper is to develop a theory and a
technique to calculate orbifold invariants {\it with}
$\Gamma$-extensions of symmetric products.

For the remaining part of this section, we explain the notion of
twisted sectors (or spaces) from which follows the geometric origin of
our definition of the orbifold invariant (2-1) associated to a group
$\Gamma$. The notion of twisted sectors arises in orbifold conformal
field theory \cite{DHVW} in which strings moving in an orbifold $M/G$
are considered. This situation is analyzed through their lifts to
$M$. When lifted to $M$, a string may not close, but its end points
are related by the action of an element $g\in G$. Thus, we are led to
consider the $g$-twisted free loop space $L_gM$ for $g\in G$ given by
$$
L_gM=\{\gamma: \Bbb R @>>> M \mid \gamma(t+1)=g\cdot\gamma(t) 
\text{ for all }t\in\Bbb R \}.
$$
On this space the centralizer $C(g)$ acts.  The $g$-twisted sector in
orbifold conformal field theory is, intuitively, a vector space
obtained by ``quantizing'' $L_gM$ (read: by taking a certain function
space on $L_gM$), and as such $C(g)$ acts projectively on this vector
space.

The same point of view can be applied to more general context. Let
$\Sigma$ be any connected manifold. We consider $\Sigma$ moving in an
orbifold $M/G$. The basic configuration space is a mapping space
$\text{Map}(\Sigma,M/G)$. To analyze this situation, we consider their
lifts to $G$-equivariant maps from $G$-principal bundles $P$ over
$\Sigma$ to $M$. Since the action of $G$ on $M$ is not assumed to be
free, not all maps $\overline{\gamma} : \Sigma @>>> M/G$ lift to
$G$-equivariant maps from principal bundles. Even when they do, there
can be lifts from several non-isomorphic principal bundles to $M$.  We
consider the set of all possible lifts $\gamma$ of
$\overline{\gamma}\in\text{Map}(\Sigma,M/G)$ for all possible
$G$-principal bundles $P$ over $\Sigma$, and we introduce an
equivalence relation among lifts: two lifts $\gamma_1:P_1 @>>> M$ and
$\gamma_2:P_2@>>> M$ of $\overline{\gamma}$ are equivalent if there
exists a $G$-bundle isomorphism $\alpha: P_1 @>{\cong}>> P_2$ inducing
the identity map on $\Sigma$ such that
$\gamma_2\circ\alpha=\gamma_1$. We denote the set of all equivalence
classes by $\Bbb L_{\Sigma}(M;G)$. There is a natural map from $\Bbb
L_{\Sigma}(M;G)$ to the original mapping space. Thus we have
$$
\pi_{\sssize \Sigma}:\Bbb L_{\Sigma}(M;G)=\coprod_{[P]}
\text{Map}_G(P,M)/\text{Aut}_G^{}(P) @>>>
\text{Map}(\Sigma,M\!/\!G).
\tag2-11
$$
Here the union runs over all the isomorphism classes of $G$-principal
bundles over $\Sigma$. When the $G$-action on $M$ is free, the above
map $\pi_{\sssize \Sigma}$ is a homeomorphism. In general, the map
$\pi_{\sssize\Sigma}$ is a kind of mild resolution of orbifold
singularities. Note that $\Bbb L_{\Sigma}(M;G)$ is still an orbifold,
although the group $\text{Aut}_G(P)$ is in general a proper subgroup of
$G$ for any $P$. The part of $\Bbb L_{\Sigma}(M;G)$ corresponding to
$[P]$ is called the $[P]$-twisted sector. Note that the space
$\text{Map}_G(P,M)$ can be described as the space of sections of a
fibre bundle $P_M=P\times_GM @>>> \Sigma$. 

The same space can be described in a different way. Since the
isomorphism class of a $G$-principal bundle $P @>>> \Sigma$ is
determined by the $G$-conjugacy class of an associated holonomy
homomorphism $\phi:\pi_1(\Sigma) @>>> G$ and since the group of
$G$-bundle automorphisms $\text{Aut}_G^{}(P)$ of a $G$-bundle $P @>>>
\Sigma$ associated to $\phi$ is given by the centralizer
$C(\phi)\subset G$, the space of total twisted sectors $\Bbb
L_{\Sigma}(M;G)$ can be described alternately as follows:
$$
\Bbb L_{\Sigma}(M;G)=\coprod_{[\phi]}
\text{Map}_{\phi}(\widetilde{\Sigma},M)/C(\phi).
\tag2-12
$$
Here the disjoint union runs over all $G$-conjugacy classes of
homomorphisms $\phi: \pi_1(\Sigma) @>>> G$, the space
$\widetilde{\Sigma}$ is the universal cover of $\Sigma$, and the
$\phi$-twisted mapping space $\text{Map}_{\phi}$ is defined by
$$
\text{Map}_{\phi}(\widetilde{\Sigma},M)=\{\gamma:\widetilde{\Sigma}
@>>> M \mid \gamma(ph)=\phi(h)^{-1}\gamma(p) \text{ for all }
p\in\widetilde{\Sigma}, h\in\pi_1(\Sigma)\}.
\tag2-13
$$
We call $\text{Map}_{\phi}(\widetilde{\Sigma},M)$ the $\phi$-twisted
sector. 

Our future goal is to investigate global topology of the infinite
dimensional twisted space $\Bbb L_{\Sigma}(M;G)$. In this paper, we
examine the subspace of (2-11) consisting of locally constant
$G$-equivariant maps. The same space can be described as the space of
constant $\phi$-equivariant maps in (2-12). Although this is a very
small portion of the space $\Bbb L_{\Sigma}(M;G)$, it already contains
very interesting information. This space of (locally) constant
equivariant maps comes equipped with a map to $M\!/\!G$:
$$
\overline{\pi}_{\Sigma}:\coprod_{[\phi]}
M^{\langle\phi\rangle}/C(\phi)
@>>> M\!/\!G,
\tag2-14
$$
where $[\phi]\in\text{Hom}(\Gamma, G)/G$, with $\Gamma=\pi_1(\Sigma)$.
Note that the space corresponding to the trivial $\phi$ is $M/G$, and
the map $M^{\langle\phi\rangle}/C(\phi) @>>> M\!/\!G$ is the obvious map:
mapping $C(\phi)$-orbits to $G$-orbits in $M$. By considering the
Euler characteristic of the left hand side, we get our definition
(2-1) of the orbifold Euler characteristic associated to the group
$\Gamma$. In this form, the group $\Gamma$ can be any group, it does
not have to be the fundamental group of a manifold. This gives our
definition (2-1) for any group $\Gamma$. However, interesting cases
are for those groups $\Gamma$ arising as fundamental groups of
manifolds.

\head
Geometry of homomorphisms into wreath products 
\endhead

For results in this section and the next, the group $G$ does not have
to be finite. We assume finiteness only in the orbifold context. 

Let $M$ be a $G$-manifold. On the $n$-fold Cartesian product $M^n$,
the direct product $G^n$ acts coordinate-wise, and also the $n$-th
symmetric group $\frak S_n$ acts by permuting factors. Combining these
two actions, we have an action of the wreath product
$G_n=G\!\wr\!\frak S_n$ on $M^n$, which is formally defined as
follows. Let $\bold n=\{1,2,\dots,n\}$ for $n\ge1$. The $n$-th
symmetric group $\frak S_n$ is, by definition, the totality of
bijections of $\bold{n}$ to itself. The wreath product $G\!\wr\!\frak
S_n$ is defined as a semi-direct product
$$
G\!\wr\!\frak S_n=\text{Map}(\bold{n},G)\rtimes\frak S_n,
\tag3-1
$$
where $\sigma\in\frak S_n$ acts on $f\in\text{Map}(\bold{n},G)$ by
$(\sigma\cdot f)(k)=f\bigl(\sigma^{-1}(k)\bigr)$ for $k\in\bold{n}$.
The product and inverse in $G\!\wr\!\frak S_n$ is given by
$(f_1,\sigma_1)(f_2,\sigma_2)=(f_1\cdot\sigma_1f_2,\sigma_1\sigma_2)$
and $(f,\sigma)^{-1}=(\sigma^{-1}f^{-1},\sigma^{-1})$ for any
$f_1,f_2,f\in \text{Map}(\bold{n},G)$ and
$\sigma_1,\sigma_2,\sigma\in\frak S_n$. When $n=0$, we set $G_n$ to be
the trivial group. 

The left action of $(f,\sigma)\in G\!\wr\!\frak S_n$ on
$(x_1,x_2,\dots,x_n)\in M^n$ is given by
$$
(f,\sigma)(x_1,x_2,\dots,x_n)
=\bigl(f(1)x_{\sigma^{-1}(1)},f(2)x_{\sigma^{-1}(2)},\dots,
f(n)x_{\sigma^{-1}(n)}\bigr)\in M^n.
\tag3-2
$$
In the calculation of the $\Gamma$-extended orbifold invariant
$\varphi_{\sssize\Gamma}(M^n;G\wr\frak S_n)$ of the $n$-th
symmetric product of an orbifold (see (5-1) for an explicit
expression), we have to deal with homomorphisms into wreath
products. Let $\theta: \Gamma @>>> G\!\wr\!\frak S_n$ be such a
homomorphism. Let $\theta(u)=\bigl(f(u,\
\cdot\ ),\phi(u)\bigr)$, where $f(u,\ \cdot\
)\in\text{Map}(\bold{n},G)$, and $\phi:\Gamma @>>>
\frak S_n$ is the composition of $\theta$ and the canonical projection
map $G\!\wr\!\frak S_n @>>> \frak S_n$. Through $\phi$, any $\theta$
induces a $\Gamma$-set structure on $\bold{n}$. That $\theta$ is a
homomorphism immediately implies that the two variable function $f(\
\cdot\ ,\ \cdot\ ):\Gamma\times\bold{n} @>>> G$ satisfies the
following identities:
$$
\aligned
f(u_1u_2,\ell)
&=f(u_1,\ell)f\bigl(u_2,\phi(u_1)^{-1}(\ell)\bigr) \\
f(1,\ell)&=1 \\
f(u,\ell)^{-1}&=f\bigl(u^{-1},\phi(u)^{-1}(\ell)\bigr)
\endaligned
\tag3-3
$$
for any $u_1,u_2\in\Gamma$, and $\ell\in\bold{n}$.

Recall that there is a bijective correspondence between the conjugacy
classes of homomorphisms $\text{Hom}(\Gamma,\frak S_n)/\frak S_n$ and
the set of isomorphism classes of $\Gamma$-sets of order $n$. We prove
a similar correspondence for homomorphisms into wreath
products. Namely,

\proclaim{Theorem 3-1}  Let $G$ and $\Gamma$ be any group. Then 
the following bijective correspondence exists\rom{:} 
$$
\biggl\{\foldedtext\foldedwidth{3in}
{Isomorphism classes of $\Gamma$-equivariant $G$-principal bundles
over a $\Gamma$-set of order $n$}
\biggr\}
\overset{1:1}\to\longleftrightarrow
\text{\rm Hom}(\Gamma,G_n)/G_n. 
\tag3-4
$$
\endproclaim

Before proving this theorem, we discuss some generalities on
$\Gamma$-equivariant $G$-principal bundles. These objects arise
naturally in the context of $G$-bundles over finite covering
spaces. See \S7 for details. Let $\pi: P @>>> Z$ be such a
bundle. Then $P @>>> Z$ is a (right) $G$-principal bundle over a left
$\Gamma$-set $Z$ and $\Gamma$ acts $G$-equivariantly on $P$ from the
left, inducing the $\Gamma$-action on $Z$. Thus, the left $\Gamma$
action and the right $G$ action on $P$ commute. We can also let
$\Gamma\times G$ act on $P$ from the right by setting
$p\cdot(u,g)=u^{-1}pg\in P$ for any $p\in P$, $u\in\Gamma$ and $g\in
G$. We use both of these two view points on group actions on $P$.

For $p\in P$, let $x=\pi(p)\in Z$ and let $H_x$ be the isotropy
subgroup at $x\in Z$ of the $\Gamma$ action on $Z$. Since any $h\in
H_x$ preserves the fibre $\pi^{-1}(x)$, for any point
$p\in\pi^{-1}(x)$ we may write $h\cdot p=p\cdot\rho_p(h)$ for some
unique $\rho_p(h)\in G$. Since left $\Gamma$ action and the right $G$
action commute, it can be easily checked that $\rho_p: H_x @>>> G$ is
a homomorphism. Thus, to each point $p\in P$, there is an associated
homomorphism $\rho_p$. At any two points $p$ and $p'$ of $P$ over the
same $\Gamma$-orbit in $Z$, we can easily check that the corresponding
homomorphisms $\rho_p$ and $\rho_{p'}$ are related by
$$
\rho_{p'}(uhu^{-1})=g^{-1}\rho_p(h)g, \qquad p'=upg ,
\tag3-5
$$
for some $u\in \Gamma$ and $g\in G$. Here $h\in H_{\pi(p)}$ so that
$uhu^{-1}\in H_{\pi(p')}$.

The basic structure of a $\Gamma$-equivariant $G$-principal bundle over
a $\Gamma$-transitive set is simple.

\proclaim{Lemma 3-2} Let $\pi:P @>>> Z$ be a $\Gamma$-equivariant
$G$-principal bundle over a transitive $\Gamma$-set $Z$. For any point
$p_0\in P$, let $\rho : H @>>> G$ be the homomorphism associated to
$p_0\in P$, where $H\subset\Gamma$ is the isotropy subgroup at
$\pi(p)\in Z$ so that $Z\cong \Gamma\!/\!H$ as $\Gamma$-sets. Then, we
have a $\Gamma$-equivariant $G$-principal bundle isomorphism
$\Gamma\times_{\rho}G \cong P$, where $[u,g]\in\Gamma\times_{\rho}G$
corresponds to $up_0g\in P$. Here, the equivalence relation in
$P_{\rho}=\Gamma\times_{\rho}G$ is given by $[u,g]=[uh,\rho(h)^{-1}g]$
for all $h\in H$.
\endproclaim

The proof is straightforward. One notes that $\Gamma\times G$ acts
transitively on $P$ from the right whose isotropy subgroup at $p_0\in
P$ is $\{(h,\rho(h))\mid h\in H\}\cong H$.

The reason of our use of $G$-bundles in our study of wreath products
is that wreath products arise naturally as the full automorphisms
groups of $G$-bundles over finite sets. 

\proclaim{Lemma 3-3} The group of $G$-bundle automorphisms of the
trivial $G$-bundle over the set $\bold{n}$ is the wreath product
$G\!\wr\!\frak S_n$. That is, $\text{\rm Aut}_{G}(\bold{n}\times G)
=G\!\wr\!\frak S_n$.  Here, the action of $(f,\sigma)\in G\!\wr\!\frak
S_n$ given in {\rm (3-1)} on $(\ell,g)\in \bold{n}\times G$ is given
by $(f,\sigma)\cdot(\ell,g)=
\bigl(\sigma(\ell),f(\sigma(\ell))g\bigr)$. 
\endproclaim
One can easily check that this defines an action. This action is
compatible with the action (3-2) in the following way. Through the
identification $M^n\cong\text{Map}_G(\bold{n}\times G,M)$ given by the
evaluation of $G$-equivariant maps at $1\in G$, the action of
$(f,\sigma)$ on $M^n$ given in (3-2) corresponds exactly to the action
$\bigl((f,\sigma)^{-1}\bigr)^*$ on $\text{Map}_G(\bold{n}\times G,M)$.

Lemma 3-3 gives a geometric characterization of wreath products. This
is the very reason why $G$-principal bundles are relevant to our study
of wreath products. 

We are ready to prove Theorem 3-1. 

\demo{Proof of Theorem 3-1}  First we describe the correspondence
(3-4) and show that it is well defined. 

For a given $\Gamma$-equivariant $G$-principal bundle $\pi:P @>>> Z$
with $|Z|=n$, we choose a section $s:Z @>>> P$ and a bijection
$t:\bold{n} @>{\cong}>> Z$. We then have a $G$-bundle isomorphism
$\bold{n}\times G @>{\cong}>{t}> Z\times G @>{\cong}>{s}> P$. We
transfer the $\Gamma$ action on $P$ to an action on $\bold{n}\times G$
through this isomorphism. Since the action of $\gamma\in\Gamma$ on $P$
is $G$-equivariant, the corresponding action $\theta(\gamma)$ on
$\bold{n}\times G$ is also $G$-equivariant. Since
$\text{Aut}_G(\bold{n}\times G)=G\!\wr\!\frak S_n$ by Lemma 3-3, we
have a homomorphism $\theta:\Gamma @>>> G\!\wr\!\frak S_n$.

We show that the conjugacy class $[\theta]$ is independent of the
choices $(s,t)$ and depends only on the isomorphism class of $\pi:P
@>>> Z$ as a $\Gamma$-equivariant $G$-bundle. Let $\pi_1:P_1 @>>> Z_1$
and $\pi_2: P_2 @>>> Z_2$ be two isomorphic $\Gamma$-equivariant
$G$-bundles. Let $\varphi:P_1 @>{\cong}>> P_2$ be an isomorphism. Let
$(s_1,t_1)$ and $(s_2,t_2)$ be the choices as above for $P_1$ and
$P_2$. We have the following commutative diagram. 
$$
\CD 
\bold{n}\times G @>{\cong}>{(s_1,t_1)}> P_1 @>{\cong}>{\varphi}> P_2
@<{\cong}<{(s_2,t_2)}< \bold{n}\times G \\
@V{\cong}V{\theta_1(\gamma)}V  @V{\cong}V{\gamma}V @V{\cong}V{\gamma}V 
@V{\cong}V{\theta_2(\gamma)}V \\
\bold{n}\times G @>{\cong}>{(s_1,t_1)}> P_1 @>{\cong}>{\varphi}> P_2
@<{\cong}<{(s_2,t_2)}< \bold{n}\times G 
\endCD
$$
The composition of horizontal maps is a $G$-bundle isomorphism $\eta$
of $\bold{n}\times G$, and hence $\eta\in G\!\wr\!\frak S_n$, and the
commutativity of the diagram shows that
$\theta_2(\gamma)=\eta\theta_1(\gamma)\eta^{-1}$.  Thus,
$\theta_1,\theta_2:\Gamma @>>> G\!\wr\!\frak S_n$ are conjugate in
$G\!\wr\!\frak S_n$. This proves that the correspondence (3-4) is well
defined.

Next we show that the correspondence (3-4) is onto. Let $\theta:\Gamma
@>>> G\!\wr\!\frak S_n$ be a homomorphism. Let $P_{\theta}$ be
$\bold{n}\times G$ with $\Gamma$-action induced by $\theta:
\Gamma @>>> \text{Aut}_G(P_{\theta})=G\!\wr\!\frak S_n$. Since the
left $\Gamma$-action commutes with the right $G$-action on
$P_{\theta}$, this $\Gamma$-action induces a $\Gamma$-action on
$\bold{n}$. Hence $P_{\theta} @>>> \bold{n}$ is a $\Gamma$-equivariant
$G$-bundle. With the obvious choice of $(s,t)$, the homomorphism
corresponding to $P_{\theta} @>>> \bold{n}$ is clearly $\theta :
\Gamma @>>> G\!\wr\!\frak S_n$. Hence the correspondence (3-4) is
onto. 

We give another proof of surjectivity. In this proof, a bundle $P @>>>
Z$ corresponding to $\theta$ is constructed from a different point of
view as a union of $\Gamma$-irreducible $G$-bundles.  This description
clarifies the detailed internal structures of the bundle $P @>>> Z$,
and it will be useful for later sections. Let $\theta:\Gamma @>>>
G\wr\frak S_n$ be a homomorphism. Let $\phi:\Gamma @>{\theta}>>
G\wr\frak S_n @>{\text{proj}}>> \frak S_n$ be a composition. Let
$\theta(u)=\bigl(f(u,\ \cdot\ ),\phi(u)\bigr)$ for $u\in\Gamma$. Let
$\bold{n}=\coprod_{\xi}X_{\xi}$ be the orbit decomposition of
$\bold{n}$ under the $\Gamma$-action on $\bold{n}$ defined by
$\phi$. We choose a base point $x_{\xi}\in X_{\xi}$ from each orbit,
and let $H_{\xi}$ be the isotropy subgroup at $x_{\xi}$. Now $\theta$
decomposes as $\theta=\prod_{\xi}\theta_{\xi}:
\Gamma @>>> \prod_{\xi}G\!\wr\!\frak S(X_{\xi})$, where 
$\theta_{\xi}:\Gamma @>>> G\!\wr\!\frak S(X_{\xi})$, and $\frak
S(X_{\xi})$ is the group of permutations of $X_{\xi}$. For each $\xi$,
let $\rho_{\xi}:H_{\xi} @>>> G$ be defined by
$\rho_{\xi}(u)=f(u,x_{\xi})$ for $u\in H_{\xi}$. Since $H_{\xi}$ fixes
$x_{\xi}$, the first formula in (3-3) shows that $\rho_{\xi}$ is a
homomorphism. Now set $P=\coprod_{\xi}(\Gamma\times_{\rho_{\xi}}G @>>>
X_{\xi})$. This $P$ is a $\Gamma$-equivariant $G$-principal bundle
over $\bold{n}$. We show that by choosing sections for these bundles
appropriately, the homomorphism corresponding to $P$ is precisely
$\theta$ we started with. Now a section $s_{\xi}$ for
$P_{\xi}=\Gamma\times_{\rho_{\xi}}G @>>> X_{\xi}$ is the same as a
$\rho_{\xi}$-equivariant map $s_{\xi}:\Gamma @>>> G$. Let
$s_{\xi}(u)=f(u^{-1},x_{\xi})$ for $u\in\Gamma$. By (3-3), for $h\in
H_{\xi}$ we have $s_{\xi}(uh)=f(h^{-1}u^{-1},x_{\xi})
=f(h^{-1},x_{\xi})f(u^{-1},\phi(h)x_{\xi})
=\rho_{\xi}(h)^{-1}s_{\xi}(u)$, where $\phi(h)x_{\xi}=x_{\xi}$ since
$h\in H_{\xi}$. This shows that $s_{\xi}$ is
$\rho_{\xi}$-equivariant. Thus, a map $\alpha_{\xi}: X_{\xi} @>>>
\Gamma\times_{\rho_{\xi}}G$ given by
$\alpha_{\xi}(\phi(u)x_{\xi})=[u,s_{\xi}(u)]$ for $u\in\Gamma$ is a
well defined section of $P_{\xi} @>>> X_{\xi}$. Using this section, we
have a $G$-bundle isomorphism $X_{\xi}\times G @>>>
\Gamma\times_{\rho_{\xi}}G$ mapping $(\phi(u)x_{\xi},g)$ to
$\alpha_{\xi}(\phi(u)x_{\xi})g=[u,s_{\xi}(u)g]$. We transfer the
$\Gamma$ action on $\Gamma\times_{\rho_{\xi}}G$ to $X_{\xi}\times G$
using this isomorphism. Let $\theta_{\xi}(\gamma)$ be the action on
$X_{\xi}\times G$ corresponding to the action of $\gamma$ on
$\Gamma\times_{\rho_{\xi}}G$, in which the action of $\Gamma$ is given
by $\gamma[u,s_{\xi}(u)g]=[\gamma u,s_{\xi}(u)g]$. Let
$\theta_{\xi}(\gamma)(\phi(u)x_{\xi},g)=(\phi(\gamma u)x_{\xi},g')$
for some $g'\in G$. The corresponding element in
$\Gamma\times_{\rho_{\xi}}G$ is $[\gamma u,s_{\xi}(\gamma u)g']$ which
must be equal to $[\gamma u, s_{\xi}(u)g]$. Using (3-3), we have
$s_{\xi}(\gamma u)=s_{\xi}(u)f(\gamma,\phi(\gamma u)x_{\xi})^{-1}$,
hence $g'=f(\gamma,\phi(\gamma u)x_{\xi})g$. This means that
$\theta_{\xi}(\gamma)\bigl(\phi(u)x_{\xi},g\bigr) =\bigl(\phi(\gamma
u)x_{\xi}, f(\gamma,\phi(\gamma u)x_{\xi})g\bigr) =\bigl(f(\gamma,\
),\phi(\gamma)\bigr)\cdot(\phi(u)x_{\xi},g)$ by the description of the
action of $G_n$ on $\bold{n}\times G$ given in Lemma 3-3. This shows
that the $\Gamma$-action $\theta_{\xi}'$ on $X_{\xi}\times G$ is given
by $\theta_{\xi}'(\gamma)=\bigl(f(\gamma,\ ),\phi(\gamma)\bigr)$,
which is the original $\theta$ restricted to the bundle over
$X_{\xi}$, namely $\theta_{\xi}:\Gamma @>>> G\!\wr\!\frak S(X_{\xi})$.
Hence the homomorphism $\prod\theta_{\xi}'$ corresponding to the
$\Gamma$-equivariant $G$-bundle
$P=\coprod_{\xi}(\Gamma\times_{\rho_{\xi}}G @>>> X_{\xi})$ with chosen
sections $\{\alpha_{\xi}\}_{\xi}$ is precisely $\theta :\Gamma @>>>
G\!\wr\!\frak S_n$. This proves that the correspondence (3-4) is onto.

Finally, we show that the correspondence (3-4) is injective.  Let
$\pi_i:P_i @>>> Z_i$ for $i=1,2$ be two $\Gamma$-equivariant
$G$-bundles such that corresponding homomorphisms $\theta_1$ and
$\theta_2$, after making some choices $(s_i,t_i)$ for $i=1,2$, are
conjugate by an element $\eta\in G\!\wr\!\frak S_n$. We then have the
following commutative diagram.
$$
\CD 
P_1 @<{\cong}<{(s_1,t_1)}< \bold{n}\times G @>{\cong}>{\eta}> 
\bold{n}\times G @>{\cong}>{(s_2,t_2)}> P_2 \\
@V{\cong}V{\gamma}V @V{\cong}V{\theta_1(\gamma)}V 
@V{\cong}V{\theta_2(\gamma)}V @V{\cong}V{\gamma}V \\
P_1 @<{\cong}<{(s_1,t_1)}< \bold{n}\times G @>{\cong}>{\eta}> 
\bold{n}\times G @>{\cong}>{(s_2,t_2)}> P_2
\endCD
$$
The composition $\varphi$ of the horizontal maps is a $G$-bundle
isomorphism. The commutativity of the diagram shows that $\varphi:P_1
@>>> P_2$ is a $G$-bundle map commuting with the action of
$\Gamma$. Hence $\varphi$ is a $\Gamma$-equivariant $G$-bundle
isomorphism. This shows that the correspondence in (3-4) is
injective. This completes the proof.
\qed
\enddemo

Next we discuss isomorphisms between $\Gamma$-equivariant
$G$-principal bundles. Recall that for any point $p\in P$, we can
associate a homomorphism $\rho_p:H_{\pi(p)} @>>> G$. This homomorphism
completely characterize the $\Gamma\text{-}G$ bundle $P @>>> Z$ when
the base is a transitive $\Gamma$-set.

\proclaim{Proposition 3-4}  \rom{(1)}\qua Let $\pi_i:P_i @>>> Z_i$ be
$\Gamma$-equivariant $G$-bundles for $i=1,2$. Let $\varphi:P_1
@>{\cong}>> P_2$ be $\Gamma$-equivariant $G$-bundle isomorphism. Then
for every $p_1\in P_1$, letting $p_2=\varphi(p_1)$, we have 
$$
\rho_{p_1}=\rho_{p_2}: H @>>> G,
\tag3-6
$$
where $H=H_{\pi_1(p_1)}=H_{\pi_2(p_2)}\subset\Gamma$ is the isotropy
subgroup at $\pi_1(p_1)\in Z_1$ and $\pi_2(p_2)\in Z_2$. 

\rom{(2)}\qua Conversely, if there exists $p_1\in P_1$ and $p_2\in P_2$
such that $\rho_{p_1}=\rho_{p_2}$ and if $Z_1$ and $Z_2$ are
transitive $\Gamma$-sets, then there exists a unique
$\Gamma$-equivariant $G$-bundle isomorphism $\varphi:P_1 @>{\cong}>>
P_2$ such that $\varphi(p_1)=p_2$.
\endproclaim
\demo{Proof} The $\Gamma$-equivariant $G$-bundle isomorphism $\varphi$
induces a $\Gamma$-isomorphism $\overline{\varphi}:Z_1 @>>> Z_2$ on
the base sets. Hence isotropy subgroups at $x_1=\pi_1(p_1)$ and
$x_2=\pi_2(p_2)=\overline{\varphi}(x_1)$ must be the same subgroup of
$\Gamma$ which we call $H$. Let the homomorphisms associated to points
$p_1$ and $p_2$ be $\rho_1,\rho_2:H @>>> G$. For any $h\in H$,
$\Gamma$-equivariance of $\varphi$ and the definition of $\rho_2$
imply that $\varphi(hp_1)=h\varphi(p_1)=hp_2=p_2\rho_2(h)$. On the
other hand, $G$-equivariance implies that $\varphi(hp_1)=
\varphi\bigl(p_1\rho_1(h)\bigr)=\varphi(p_1)\rho_1(h)=p_2\rho_1(h)$.
Since the $G$-action on $P_2$ is free, we have $\rho_1(h)=\rho_2(h)$
for all $h\in H$. This proves the first part.

For the second part, we first show uniqueness of $\varphi$. Suppose a
$\Gamma$-equivariant $G$-bundle isomorphism $\varphi: P_1 @>>> P_2$
such that $\varphi(p_1)=p_2$ exists. Then it must have the property
$\varphi(up_1g)=up_2g$ for all $u\in\Gamma$ and $g\in G$. Since $Z_1$
and $Z_2$ are transitive $\Gamma$-sets, the actions of $\Gamma\times G$
on $P_1$ and $P_2$ are transitive. Thus the above identity uniquely
determines $\varphi$. This proves uniqueness of $\varphi$.

For existence, since $\rho_{p_1}=\rho_{p_2}=\rho$, say, we have $P_1
@<{\cong}<{\varphi_1}<
\Gamma\times_{\rho}G @>{\cong}>{\varphi_2}> P_2$, where $\varphi_i$
are $\Gamma$-equivariant $G$-bundle isomorphisms in Lemma 3-2 such
that $\varphi_i([1,1])=up_ig$ for $i=1,2$. Then
$\varphi=\varphi_2\circ\varphi_1^{-1}$ is the desired $\Gamma$-$G$
isomorphism.  This completes the proof.
\qed
\enddemo

Next, we classify $\Gamma$-equivariant $G$-bundles over transitive
$\Gamma$-sets. Let $\pi:P @>>> Z$ be such a bundle. For $z\in Z$, let
$H_z$ be the isotropy subgroup at $z$. Then the collection of
homomorphisms $\rho_p:H_z @>>> G$ for all $p\in\pi^{-1}(z)$ forms a
complete $G$-conjugacy class, that is, an element in
$\text{Hom}(H_z,G)/G$ in view of (3-5). The normalizer
$N_{\Gamma}(H_z)$ acts on this set of conjugacy classes by conjugating
$H_z$. For each point $z\in Z$, we consider the set of
$N_{\Gamma}(H_z)$-orbits in $\text{Hom}(H_z,G)/G$. This gives us a
bundle of sets: 
$$
\omega: \coprod_{z\in Z}\text{Hom}(H_z,G)/(N_{\Gamma}(H_z)\times G)
@>>> Z.
\tag3-7
$$
This bundle $\omega$ depends only on the $\Gamma$-set $Z$. The above
argument shows that for any $\Gamma$-equivariant $G$-bundle $P$ over
$Z$, we have an associated section $\eta_{\sssize P}$ of $\omega$
whose value at $z\in Z$ is the $N_{\Gamma}(H_z)$-orbit of the
$G$-conjugacy class $[\rho_p]\in\text{Hom}(H_z,G)/G$ for any
$p\in\pi^{-1}(z)$. Note that for any $P$ over a $\Gamma$-transitive
set $Z$, this section $\eta_{\sssize P}$ is determined by its value at
any single point $z\in Z$, since $\rho_p$ for $p\in\pi^{-1}(z)$
determines $\rho_p'$ for any other point $p'\in P$ when $Z$ is
$\Gamma$-transitive by (3-5). We show that this section $\eta_{\sssize
P}$ of the bundle $\omega$, or equivalently, its value at any point of
$Z$, determines the $\Gamma$-$G$ isomorphism class of the bundle $\pi:
P @>>> Z$ for a $\Gamma$-transitive set $Z$.

\proclaim{Proposition 3-5} Let $\pi_i: P_i @>>> Z$ be
$\Gamma$-equivariant $G$ principal bundles over the same transitive
$\Gamma$-set $Z$ for $i=1,2$. Then the following statements are
equivalent\rom{:}

\rom{(1)}\qua Two bundles $P_1$ and $P_2$ are isomorphic as
$\Gamma$-equivariant $G$-principal bundles. 

\rom{(2)}\qua For some $z_0\in Z$, we have $\eta_{\sssize P_1}(z_0)
=\eta_{\sssize P_2}(z_0)\in
\text{\rm Hom}(H_{z_0},G)/(N_{\Gamma}(H_{z_0})\times G)$. 

\rom{(3)}\qua For all $z\in Z$, we have  $\eta_{\sssize P_1}(z)
=\eta_{\sssize P_2}(z)\in
\text{\rm Hom}(H_z,G)/(N_{\Gamma}(H_z)\times G)$.
\endproclaim
\demo{Proof} $(1)\Rightarrow(3).$ Let $\varphi:P_1 @>>> P_2$ be a 
$\Gamma$-equivariant $G$-bundle isomorphism. For any $p_1\in P_1$, let
$p_2=\varphi(p_1)$, and $\pi_1(p_1)=z_1$, $\pi_2(p_2)=z_2$. Since the
induced $\Gamma$-isomorphism $\overline{\varphi}$ on $Z$ is such that
$\overline{\varphi}(z_1)=(z_2)$, the isotropy subgroups at these
points are the same $H_{z_1}=H_{z_2}$, and since $\Gamma$ acts
transitively on $Z$, we must have $z_2=uz_1$ for some $u\in
N_{\Gamma}(H_{z_1})$. To compare $\eta_{\sssize P_1}$ and
$\eta_{\sssize P_2}$ over the same point $z_1$, let
$p_2'=u^{-1}p_2g^{-1}\in P_2$ for any $g\in G$. Then both $p_2'$ and
$p_1$ are points above $z_1$. Let $\rho_i=\rho_{p_i}$ for $i=1,2$, and
$\rho_2'=\rho_{p_2'}$. Then by (3-5), we have
$\rho_2(uhu^{-1})=g^{-1}\rho_2'(h)g$ for all $h\in H_{z_1}$. Since
$\varphi(p_1)=p_2$, by Proposition 3-6 we have $\rho_1=\rho_2$. Thus,
we have $\rho_1(uhu^{-1})=g^{-1}\rho_2'(h)g$ for all $h\in
H_{z_1}$. This means that $\eta_{\sssize P_1}(z_1)=[\rho_1]
=[\rho_2']=\eta_{\sssize P_2}(z_1)$ in
$\text{Hom}(H_{z_1},G)/(N_{\Gamma}(H_{z_1})\times G)$. Since $p_1\in
P_1$ can be an arbitrary point, we have $\eta_{\sssize
P_1}=\eta_{\sssize P_2}$ as sections of the bundle $\omega$. This
proves (3). 

$(3)\Rightarrow (2).$ This is obvious. 

$(2)\Rightarrow(1).$ Suppose for some $z_0\in Z$, we have $\eta_{\sssize
P_1}(z_0)=\eta_{\sssize P_2}(z_0)$. this means that for any choices of
points $p_1\in\pi_1^{-1}(z_0)$ and $p_2\in\pi_2^{-1}(z_0)$, there exists
$u\in N_{\Gamma}(H_{z_0})$ and $g\in G$ such that
$\rho_1(uhu^{-1})=g^{-1}\rho_2(h)g$ for all $h\in H_{z_0}$. Here we put
$\rho_i=\rho_{p_i}$ for $i=1,2$. Let $p_2'=up_2g\in P_2$. Then,
letting $\rho_2'=\rho_{p_2'}$, we have
$\rho_2'(uhu^{-1})=g^{-1}\rho_2(h)g$ for all $h\in H_{z_0}$ by
(3-5). Thus, we have $\rho_1=\rho_2':H_{z_0} @>>> G$. By the second part
of Proposition 3-4, there exists a unique $\Gamma$-equivariant
$G$-bundle map $\varphi:P_1 @>>> P_2$ such that
$\varphi(p_1)=p_2'$. Note that the induced $\Gamma$-map
$\overline{\varphi}$ on $Z$ is such that
$\overline{\varphi}(z_0)=uz_0$, and $z_0$ and $uz_0$ have the same
isotropy subgroups. 
This completes the proof.
\qed
\enddemo

\remark{Remark} We can also think of the above situation in the
following ({\it better}) way. When $Z$ is a transitive $\Gamma$-set,
we consider a bundle over $Z$:
$$
\coprod_{z\in Z}\text{Hom}(H_{z},G)/G @>>> Z.
\tag3-8
$$
Any $\Gamma$-equivariant $G$ bundle $\pi:P@>>> Z$ gives rise to a
section $s_{\sssize P}$ whose value at $z\in Z$ is the $G$-conjugacy
class $\{\rho_p\}_{p\in\pi^{-1}(z)}$. Now this section $s_{\sssize P}$
is $\Gamma$-equivariant. Here the result of an action of $u\in\Gamma$
on $f\in\text{Hom}(H,G)$ is given by $uf\in\text{Hom}(uHu^{-1},G)$
defined by $(uf)(uhu^{-1})=f(h)$ for any $h\in H$. Thus dividing by
the action of $\Gamma$, each bundle $P$ determines a unique element
$[s_{\sssize P}]$ whose representative is given by $\eta_{\sssize
P}(z)\in
\text{Hom}(H_z,G)/(N_{\Gamma}(H_z)\times G)$ for any $z\in
Z$. Proposition 3-5 says that this element classifies the
$\Gamma$-equivariant $G$-bundle isomorphism class of $\pi:P@>>> Z$.
\endremark

Reformulating this proposition, we have the classification theorem of
$\Gamma$-equivar\-iant $G$-principal bundles over transitive
$\Gamma$-sets. 

\proclaim{Theorem 3-6}{\rm (Classification of $\Gamma$-$G$ bundles)}\qua Let
$G$ and $\Gamma$ be any groups. Then there exists the following
bijective correspondence\rom{:}
$$
\biggl\{\foldedtext\foldedwidth{2.4in}{Isomorphism classes of
$\Gamma$-irreducible $G$-principal bundles over 
$\Gamma$-sets of order $n$}\biggr\} \overset{1:1}\to\longleftrightarrow
\!\!\coprod_{[H]_n}
\!\!\text{\rm Hom}(H,G)/(N_{\Gamma}(H)\times G),
\tag3-9
$$
where $[H]_n$ runs over the set of conjugacy classes of index $n$ 
subgroups of $\Gamma$. 

If we allow $\Gamma$-irreducible $G$-bundles over arbitrary not
necessarily finite $\Gamma$-sets, then the bijective correspondence is
still valid without the finiteness restriction on $|\Gamma\!/\!H|$ on
the right hand side.
\endproclaim

\head
Centralizers of homomorphisms into wreath products 
\endhead

Let $\theta: \Gamma @>>> G\!\wr\!\frak S_n=G_n$ be a homomorphism into a
wreath product. In this section, we determine the structure of the
centralizer of the image of $\theta$ in $G_n$. Of course this is a
purely group theoretic problem, but we found it illuminating and 
simpler to consider this problem from a geometric point of view.

From the proof of Theorem 3-1, given $\theta$ as above, there exists a
$\Gamma$-equivariant $G$-principal bundle $\pi:P @>>> Z$ over a
$\Gamma$-set $Z$ of order $n$ such that for appropriate choices of a
section $s:Z @>>> P$ and a bijection $t:\bold{n} @>{\cong}>> Z$, we
have the following commutative diagram of groups and homomorphisms.
$$
\CD
\Gamma @= \Gamma \\
@V{\theta_P}VV @VV{\theta}V  \\
\text{Aut}_G(P) @>{\cong}>{(s,t)_*}> G\!\wr\!\frak S_n
\endCD
$$
Here the bottom arrow is the one induced by a trivialization
$(s,t):\bold{n}\times G \cong P$. Now any element in the centralizer
$C(\theta)\subset G_n$ corresponds, under $(s,t)_*$, to a $G$-bundle
automorphism of $P$ which commutes with the action $\theta_P$ of
$\Gamma$ on $P$. These are $\Gamma$-equivariant $G$-bundle
automorphisms of $P$ which we denote by
$\text{Aut}_{\Gamma\text{-}G}(P)$.  This proves the following lemma.

\proclaim{Lemma 4-1} Let $\theta : \Gamma @>>> G\!\wr\!\frak S_n$ be a
homomorphism. Let $\pi: P@>>> Z$ be any $\Gamma$-equivariant
$G$-bundle whose isomorphism class corresponds to the $G_n$-conjugacy
class $[\,\theta\,]$. Then we have $C(\theta)\cong
\text{\rm Aut}_{\Gamma\text{-}G}(P)$. 
\endproclaim

When the base $Z$ is a transitive $\Gamma$-set, we have called $\pi: P
@>>> Z$ a $\Gamma$-irreducible $G$-bundle. In this case, the group
$\Gamma\times G$ acts transitively on the total space $P$. Our next
task is to describe $\Gamma\text{-}G$ automorphisms of a
$\Gamma$-irreducible $G$-bundle.  By Lemma 3-2, such a bundle is of
the form $\pi_{\rho}: P_{\rho}=\Gamma\times_{\rho}G @>>>
Z=\Gamma\!/\!H$ for some subgroup $H\subset\Gamma$ and some
homomorphism $\rho:H @>>> G$. We have the following obvious exact
sequence:
$$
1 @>>> \text{Aut}_{\Gamma\text{-}G}(P_{\rho})_Z 
@>>> \text{Aut}_{\Gamma\text{-}G}(P_{\rho})
@>>> \text{Aut}_{\Gamma}^P(Z) @>>> 1.
$$
Here $\text{Aut}_{\Gamma\text{-}G}(P_{\rho})_Z$ is a subgroup of
$\text{Aut}_{\Gamma\text{-}G}(P_{\rho})$ inducing the identity map on
$Z$. The group $\text{Aut}_{\Gamma}^P(Z)$ is a subgroup of the group
$\text{Aut}_{\Gamma}(Z)\cong H\backslash N_{\Gamma}(H)$ of
$\Gamma$-automorphisms of $Z$ which extend to $\Gamma\text{-}G$
automorphisms of $P_{\rho}$.

\remark{Remark} Since $H$ is a normal subgroup of $N_{\Gamma}^{}(H)$,
it does not matter whether we use right cosets $H\backslash
N_{\Gamma}^{}(H)$ or left cosets $N_{\Gamma}^{}(H)/H$: they are the
same group. We use right cosets notation here only because it is more
appropriate due to the formula of the isomorphism $\Phi: H\backslash
N_{\Gamma}(H)\cong \text{Aut}_{\Gamma}(Z)$ in the proof of Theorem
4-2.
\endremark

For $\rho:H @>>> G$ as above and $u\in N_{\Gamma}(H)$, let $\rho^u:H
@>>> G$ be a homomorphism defined by conjugation by $u$, that is, 
$\rho^u(h)=\rho(uhu^{-1})$ for $h\in H$. Let 
$$
N_{\Gamma}^{\rho}(H)=\{u\in N_{\Gamma}(H) \mid \rho^u \text{ and }\rho
\text{ are $G$-conjugate} \}.
\tag4-1
$$
This subset of $N_{\Gamma}(H)$ is easily seen to be a subgroup of
$N_{\Gamma}^{}(H)$. Another way to think about
$N_{\Gamma}^{\rho}(H)$ is as the isotropy subgroup at $[\rho]$
of the $N_{\Gamma}^{}(H)$-action on $\text{Hom}(H,G)/G$. 

\proclaim{Theorem 4-2} For $\rho: H @>>> G$, let $\pi_{\rho}: P_{\rho}
@>>> Z(=\Gamma\!/\!H)$ be the associated $\Gamma$-irreducible $G$-bundle.
Then $\text{\rm Aut}_{\Gamma\text{-}G}(P_{\rho})_Z\cong C_G(\rho)$ and
$\text{\rm Aut}_{\Gamma}^P(Z)\cong H\backslash
N_{\Gamma}^{\rho}(H)$. Thus, the group of $\Gamma$-$G$ automorphisms
of $P_{\rho}$ fits into the following exact sequence\rom{:}
$$
1 @>>> C_G(\rho) @>>> \text{\rm Aut}_{\Gamma\text{-}G}(P_{\rho})
@>>> H\backslash N_{\Gamma}^{\rho}(H) @>>> 1.
\tag4-2
$$
Let $p_0=[1,1]\in P_{\rho}$. For any $g'\in C_G(\rho)$, the
corresponding $\Gamma$-$G$ automorphism $\varphi_{g'}$ of $P_{\rho}$
is given by $\varphi_{g'}(up_0g)=up_0{g'}^{-1}g$ for any $u\in\Gamma$
and $g\in G$.
\endproclaim
\demo{Proof} First note that when $g'\in C(\rho)$, we have
$\rho_{p_0{g'}^{-1}}=g'\rho_{p_0}{g'}^{-1}=\rho_{p_0}$, where
$\rho_{p_0}=\rho$. Hence by Proposition 3-4, the map $\varphi_{g'}$
given above is indeed a $\Gamma$-equivariant $G$-bundle isomorphism
inducing the identity map on the base $Z$. This gives us a map
$C(\rho) @>>> \text{\rm Aut}_{\Gamma\text{-}G}(P_{\rho})_Z$. This can
be easily checked to be an injective homomorphism directly from the
definition of $\varphi_{g'}$. To see that this is a surjective map,
let $\varphi:P_{\rho} @>>> P_{\rho}$ be a $\Gamma$-equivariant
$G$-bundle automorphism of $P_{\rho}$ inducing the identity map on
$Z$. Then $\varphi(p_0)$ is a point in the same fibre. So
$\varphi(p_0)=p_0{g'}^{-1}$ for some unique $g'\in G$. By Proposition
3-4, we must have $\rho_{p_0}=\rho_{p_0{g'}^{-1}}$. By (3-5), this
means that ${g'}\rho_{p_0}(h){g'}^{-1}=\rho_{p_0}(h)$ for any $h\in
H$. Thus, $g'\in C(\rho)$. This proves that $C(\rho)\cong \text{\rm
Aut}_{\Gamma\text{-}G}(P_{\rho})_Z$.

Let $z_0=[H]\in\Gamma\!/\!H$. Recall that the isomorphism
$\Phi:H\backslash N_{\Gamma}(H) @>{\cong}>> \text{Aut}_{\Gamma}(Z)$ is
given by $\Phi(Hu')(uz_0)=u{u'}^{-1}z_0$. If $u'\in
N_{\Gamma}^{\rho}(H)$, then by definition (4-1) there exists $g'\in G$
such that $\rho({u'}^{-1}hu')={g'}^{-1}\rho(h)g'$ for any $h\in
H$. This means $\rho_{p_0}=\rho_{{u'}^{-1}p_0g'}$ by (3-5), where
$p_0=[1,1]\in P_{\rho}$. Then by
Proposition 3-4, there exists a unique $\Gamma$-equivariant $G$-bundle
isomorphism $\varphi:P @>>> P$ such that
$\varphi(up_0g)=u({u'}^{-1}p_0g')g$ for $u\in\Gamma$ and $g\in G$. The
induced $\Gamma$-map on the base $Z$ is
$\overline{\varphi}(uz_0)=u{u'}^{-1}z_0$. Thus, any element in the
image of the subgroup $H\backslash N_{\Gamma}^{\rho}(H)$ under $\Phi$
lifts to a $\Gamma$-$G$-map of $P$. There are $|C(\rho)|$ many choices
of such lifts. Thus by restriction, we have a well defined injective
homomorphism $\Phi': H\backslash N_{\Gamma}^{\rho}(H) @>>>
\text{Aut}_{\Gamma}^P(Z)$. To see that this map is surjective, let
$\overline{\varphi}: Z @>>> Z$ be a $\Gamma$-equivariant isomorphism
induced from a $\Gamma$-$G$ map $\varphi :P_{\rho} @>>>
P_{\rho}$. Suppose $\varphi(p_0)={u'}^{-1}p_0g'$ for some
$u'\in\Gamma$ and $g'\in G$. Then the $\Gamma$-map
$\overline{\varphi}$ on $Z$ is given by
$\overline{\varphi}(uz_0)=u{u'}^{-1}z_0$. By Proposition 3-4,
$\rho_{p_0}=\rho_{{u'}^{-1}p_0g'}$. This means
$\rho({u'}^{-1}hu')={g'}^{-1}\rho(h)g'$ for all $h\in H$ by
(3-5). This in turn means that $u'\in N_{\Gamma}^{\rho}(H)$ by
(4-1). Hence we have $\overline{\varphi}=\Phi'(Hu')$, and this shows
that $\Phi'$ is surjective. Thus $\Phi'$ is an isomorphism of groups. 
This completes the proof.
\qed
\enddemo

We give another description of $\text{Aut}_{\Gamma}^P(Z)$. Recall that
given a $\Gamma$-equivariant $G$-bundle $\pi:P@>>> Z$, we have a
section $s_{\sssize P}$ of the bundle (3-8) over $Z$. Let $D_P^{}(Z)$
be the set of all $\Gamma$-automorphisms of $Z$ which leaves
$s_{\sssize P}$ invariant:
$$
D_P^{}(Z)=\{\tau:Z @>>> Z \mid \text{$\tau$ is
$\Gamma$-equivariant bijection and $s_{\sssize P}\circ\tau=
s_{\sssize P}$}\}.
\tag4-3
$$
We show that any element in $D_P^{}(Z)$ extends to a
$\Gamma$-equivariant $G$-bundle map of $P$. For the next proposition,
we do not have to assume that $Z$ is a transitive $\Gamma$-set. 

\proclaim{Proposition 4-3} Let $\pi:P @>>> Z$ be a
$\Gamma$-equivariant $G$-principal bundle over a \rom{(}not
necessarily transitive\rom{)} $\Gamma$-set $Z$. With the above
notations, we have $\text{\rm Aut}_{\Gamma}^P(Z)=D_P^{}(Z)$.
\endproclaim
\demo{Proof} Let $\varphi: P@>{\cong}>> P$ be a $\Gamma$-equivariant
$G$-bundle automorphism. We show that the induced map
$\overline{\varphi}$ on the base is in $D_P^{}(Z)$. Since
$\overline{\varphi}$ is a $\Gamma$-map, for any $p\in P$, isotropy
subgroups at $\pi(\varphi(p))$ and $\pi(p)$ are the same. By
Proposition 3-4 we have $\rho_{\varphi(p)}=\rho_p$. Now, $s_{\sssize
P}\circ\overline{\varphi}(\pi(p))=[\rho_{\varphi(p)}]
=[\rho_p]=s_{\sssize P}(\pi(p))$ for any $p\in P$. Thus $s_{\sssize
P}\circ\overline{\varphi}=s_{\sssize P}$, and $\overline{\varphi}\in
D_P^{}(Z)$. This proves $\text{Aut}_{\Gamma}^P(Z)\subseteq D_P^{}(Z)$.

For the other direction of inclusion, let $\tau: Z @>>> Z$ be a
$\Gamma$-map such that $s_{\sssize P}\circ\tau=s_{\sssize P}$. Let
$Z=\coprod_{\xi}Z_{\xi}$ be the $\Gamma$-orbit decomposition of $Z$,
and let the corresponding decomposition of $\Gamma$-$G$ bundles be $(P
@>>> Z)=\coprod_{\xi}(P_{\xi} @>>> Z_{\xi})$. The $\Gamma$-map $\tau$
permutes orbits of the same type. We pick any point $z\in Z$, say,
$z\in Z_{\xi}$. Let $\tau(z)=z'\in Z_{\xi'}$. Thus isotropy subgroups
at $z$ and $z'$ are the same, say $H$. Since $s_{\sssize P}\circ
\tau(z)= s_{\sssize P}(z)$, for any points $p\in\pi^{-1}(z)$ and
$p'\in\pi^{-1}(z')$ we have $[\rho_p]=[\rho_{p'}]$ in
$\text{Hom}(H,G)/G$. Thus there exists $g\in G$ such that
$\rho_{p'}=g\rho_pg^{-1}$. By (3-5), letting $p''=p'g$ we have
$\rho_{p''}=\rho_p$. Hence by Proposition 3-4, there exists a unique
$\Gamma$-equivariant $G$-bundle map $\varphi_{\xi}:P_{\xi} @>>>
P_{\xi'}$ such that $\varphi_{\xi}(p)=p''$. Since the induced
$\Gamma$-map on the base is such that
$\overline{\varphi}_{\xi}(z)=\pi(p'')=z'=\tau(z)$ and both $Z_{\xi}$
and $Z_{\xi'}$ are transitive $\Gamma$-sets, we have
$\overline{\varphi}_{\xi}=\tau|_{Z_{\xi}}:Z_{\xi} @>>> Z_{\xi'}$. We
repeat this construction for the remaining $Z_{\xi}$'s, and we get a
$\Gamma$-equivariant $G$-bundle map $\varphi:P@>>> P$ lifting
$\tau$. Hence $\tau\in\text{Aut}_{\Gamma}^P(Z)$. This completes the
proof.
\qed
\enddemo

Next we give an alternate description of the structure of the group
$\text{Aut}_{\Gamma\text{-}G}(P_{\rho})$ when the base set of $P$ is a
transitive $\Gamma$-set. This description is more direct than the one
in Theorem 4-2. We take the point of view that $(u,g)\in\Gamma\times
G$ acts on $p\in P_{\rho}$ from the {\it right} by
$p\cdot(u,g)=u^{-1}pg$. The isotropy subgroup at $p_0=[1,1]\in
P_{\rho}$ is $\{(h,\rho(h))\mid h\in H\}\cong H$, so that as a right
$\Gamma\times G$-set, $P_{\rho}\cong H\backslash(\Gamma\times G)$
using right cosets.

By Proposition 3-4, we have a bijection $P_{\rho}\supset\{p\in
P_{\rho}\mid \rho_p=\rho\} @>{\cong}>>
\text{Aut}_{\Gamma\text{-}G}(P_{\rho})$ as sets. When $p=u^{-1}p_0g$,
by (3-5) we have $\rho_p(h)=g^{-1}\rho(uhu^{-1})g$ for all $h\in
H$. We consider a subset $T_{\rho}$ of $N_{\Gamma}^{}(H)\times
G\subset\Gamma\times G$ consisting of pairs $(u,g)$ such that the
corresponding point $p=u^{-1}p_0g$ has the property that
$\rho_p=\rho$. Namely, let
$$
T_{\rho}=\{(u,g)\in N_{\Gamma}^{}(H)\times G\mid 
g^{-1}\rho(uhu^{-1})g=\rho(h), \text{ for all }h\in H\}.
\tag4-4
$$
From this definition, we see that if $(u,g)\in T_{\rho}$, then it
follows that $u\in N_{\Gamma}^{\rho}(H)$, in view of its definition
(4-1). We can easily check that $T_{\rho}$ is a subgroup of 
$N_{\Gamma}^{}(H)\times G$. 

Another way to look at $T_{\rho}$ is as follows. Let $(u,g)\in
N_{\Gamma}(H)\times G$ act from the right on an element $f\in
\text{Hom}(H,G)$ by $[f\cdot(u,g)](h)=g^{-1}f(uhu^{-1})g$ for $h\in
H$. Then, $T_{\rho}\subset N_{\Gamma}^{}(H)\times G$ is the isotropy
subgroup at $\rho$. This makes it clear that $T_{\rho}$ is indeed a
subgroup.

We have the following diagram:
$$
\CD
\Gamma\times G @>{\text{onto}}>{H\backslash}> P_{\rho} @.  \\
@AAA  @AAA @. \\
T_{\rho} @>{\text{onto}}>{H\backslash}> \{p\in P\mid \rho_p=\rho\} 
@>{\cong}>{\text{as sets}}> \text{Aut}_{\Gamma\text{-}G}(P_{\rho})
\endCD
$$
where all vertical maps are inclusions. This shows that $H\backslash
T_{\rho}\cong \text{Aut}_{\Gamma\text{-}G}(P_{\rho})$ as sets. The
correspondence is the following one. For $(u,g)\in T_{\rho}$, it
corresponds to a point $u^{-1}p_0g\in P_{\rho}$. This point
corresponds to a $\Gamma$-equivariant $G$-bundle automorphism
$\varphi_{(u,g)}$ of $P_{\rho}$ characterized by
$\varphi_{(u,g)}(p_0)=u^{-1}p_0g$ by Proposition 3-4. In fact, this
bijective correspondence is as groups.

\proclaim{Theorem 4-4} Let an injective  homomorphism $\iota: H @>>>
T_{\rho}$ be given by $\iota(h)=\bigl(h,\rho(h)\bigr)$ for $h\in
H$. Then, $\text{\rm Im}\,\iota$ is a normal subgroup of $T_{\rho}$,
and the group of $\Gamma$-equivariant $G$-bundle automorphism of
$P_{\rho}$ fits into the following exact sequence of groups\rom{:}
$$
1 @>>> H @>{\iota}>> T_{\rho} @>{j}>> 
\text{\rm Aut}_{\Gamma\text{-}G}(P_{\rho}) @>>> 1.
\tag4-5
$$
Here, for $(u,g)\in T_{\rho}$, its image $j\bigl((u,g)\bigr)
=\varphi_{(u,g)}$ is a $\Gamma$-equivariant $G$-bundle automorphism of
$P_{\rho}$ given by $\varphi_{(u,g)}([v,k])=[vu^{-1},gk]$ for any
$v\in\Gamma$ and $k\in G$. In particular, $H\backslash T_{\rho}\cong
\text{\rm Aut}_{\Gamma\text{-}G}(P_{\rho})$ as groups. 
\endproclaim
\demo{Proof} For the first part, for any $h\in H$ and any $(u,g)\in
T_{\rho}$, we have 
$$
(u,g)\bigl(h,\rho(h)\bigr)(u,g)^{-1}
=(uhu^{-1}, g\rho(h)g^{-1})
=\bigl(uhu^{-1}, \rho(uhu^{-1})\bigr), 
$$
where the second equality holds because $(u,g)\in T_{\rho}$. Since $u\in
N_{\Gamma}(H)$, we have $uhu^{-1}\in H$. This proves that the
isomorphic image of $H$ under $\iota$ in $T_{\rho}$ is a normal
subgroup.

For the second part, it is straightforward to check that the map $j$
described above is a homomorphism with kernel $\iota(H)$. Since the
bijective correspondence $H\backslash T_{\rho}\cong
\text{\rm Aut}_{\Gamma\text{-}G}(P_{\rho})$ described earlier is
precisely the one induced by $j$ above, this correspondence is an
isomorphism of groups. This completes the proof.
\qed
\enddemo

We also have a related exact sequence involving $T_{\rho}$. Almost by
definition (4-4), we have $T_{\rho}\subset N_{\Gamma}^{\rho}(H)\times
G$. Projection to $N_{\Gamma}^{\rho}(H)$ is onto by definition (4-1),
and we have the following exact sequence of groups:
$$
1 @>>> C_G(\rho) @>>> T_{\rho} @>>> N_{\Gamma}^{\rho}(H) @>>> 1.
\tag4-6
$$
Dividing by the {\it left} action by $H$, we recover (4-2). Here, note
that $H\cap C(\rho)$ is trivial in $T_{\rho}$. 

Theorem 4-2 and Theorem 4-4 deal with the group of
$\Gamma$-equivariant $G$-bundle automorphisms of a
$\Gamma$-irreducible $G$-bundle $P_{\rho}$. Next we describe the
$\Gamma$-$G$-automorphism group of a general $\Gamma$-equivariant
$G$-bundle $\pi:P @>>> Z$, where the base set $Z$ is not necessarily a
transitive $\Gamma$-set, nor a finite $\Gamma$-set. Any such bundle
decomposes into a disjoint union of $\Gamma$-irreducible $G$-bundles,
and these irreducible bundles are classified by Theorem 3-6. For each
conjugacy class $[H]$ of subgroups of $\Gamma$, and for each conjugacy
class $[\rho]\in\text{Hom}(H,G)/(N_{\Gamma}(H)\times G)$, let
$r(H,\rho)$ be the number of isomorphic copies of
$P_{\rho}=\Gamma\times_{\rho}G @>>> \Gamma\!/\!H$ appearing in the
decomposition of $\pi:P @>>> Z$. Let $P[H,\rho]$ be the sub-bundle of
$P$ of type $([H],[\rho])$. Thus $P[H,\rho]$ is isomorphic to a
disjoint union of $r(H,\rho)$ copies of $\pi_{\rho}:P_{\rho} @>>>
\Gamma\!/\!H$. If $([H],[\rho])\not=([H'],[\rho'])$, then there are no
$\Gamma$-$G$ isomorphisms between $P[H,\rho]$ and
$P[H',\rho']$. Furthermore, it is elementary that 
$$
\text{Aut}_{\Gamma\text{-}G}(P[H,\rho])\cong 
\text{Aut}_{\Gamma\text{-}G}(P_{\rho})\!\wr\!\frak S_{r(H,\rho)}.
$$
Combining these results, we finally obtain the structure theorem for
$\text{Aut}_{\Gamma\text{-}G}(P)$. 

\proclaim{Theorem 4-5}{\rm (Structure of $\Gamma$-$G$ automorphism groups
and centralizers)}\qua Let $\pi: P @>>> Z$ be a $\Gamma$-equivariant
$G$-principal bundle over a $\Gamma$-set $Z$. Let $r(H,\rho)$ be the
number of isomorphic copies of $P_{\rho}=\Gamma\times_{\rho}G @>>>
\Gamma\!/\!H$ appearing in the irreducible decomposition of $P @>>>
Z$. Then,
$$
\text{\rm Aut}_{\Gamma\text{-}G}(P)\cong
\prod_{[H]}\prod_{[\rho]}
\text{\rm Aut}_{\Gamma\text{-}G}(P_{\rho})\!\wr\!\frak S_{r(H,\rho)},
\tag4-7
$$
where $[H]$ runs over all conjugacy classes of subgroups of $\Gamma$,
and for a given $[H]$, $[\rho]$ runs over the set $\text{\rm
Hom}(H,G)/(N_{\Gamma}(H)\times G)$. The structure of $\text{\rm
Aut}_{\Gamma\text{-}G}(P_{\rho})$ is described in Theorem {\rm 4-2} and
Theorem {\rm 4-4}. 

For a given homomorphism $\theta:\Gamma @>>> G\!\wr\!\frak S_n$, let
$\pi: P @>>> Z$ be a $\Gamma$-$G$ bundle whose isomorphism class
corresponds to $[\theta]$ under the correspondence \rom{(3-4)}. Then
the centralizer $C_{G_n}(\theta)$ is isomorphic to $\text{\rm
Aut}_{\Gamma\text{-}G}(P)$ given in \rom{(4-7)}.
\endproclaim

\head
Generating functions of orbifold Euler characteristics of 
symmetric products:\strut\ exponential formulae
\endhead

In this section, we let $G$ be a finite group and let $M$ be a
$G$-manifold. Let $\varphi(M;G)$ be an arbitrary multiplicative
orbifold invariant of $(M;G)$. Let $\Gamma$ be an arbitrary group. We
are interested in calculating the generating function
$\sum_{n\ge0}q^n\varphi_{\sssize\Gamma}(M^n;G_n)$ of the
$\Gamma$-extended orbifold invariant of symmetric products of an
orbifold. Here $G_n=G\!\wr\!\frak S_n$ is a wreath product and
$$
\varphi_{\sssize \Gamma}(M^n;G_n)
=\!\!\!\!\!\!\!\!\!\!\!\!\!\!\!\!\!\!\!\!\!\!\!\!
\sum_{\ \ \ \ \ \ \ [\theta]\in\text{Hom}(\Gamma,G_n)/G_n}
\!\!\!\!\!\!\!\!\!\!\!\!\!\!\!\!\!\!\!\!\!\!\!\!
\varphi\bigl((M^n)^{\langle\theta\rangle};C_{G_n}(\theta)\bigr).
\tag5-1
$$
Let $\theta:\Gamma @>>> G_n=G\!\wr\!\frak S_n$ be a
homomorphism. Calculating the above invariant involves identifying
fixed point subset $(M^n)^{\langle\theta\rangle}$, and calculating the
centralizer $C_{G_n}(\theta)$ together with its action on the fixed
point subset. These tasks can be done geometrically using
$\Gamma$-equivariant $G$-principal bundles. 

Given $\theta$ as above, let $\pi_{\theta}: P @>>> \bold{n}$ be a
$\Gamma$-equivariant $G$-bundle over a set $\bold{n}=\{1,2,\dots,n\}$
equipped with a section $s_{\theta}:
\bold{n} @>>> P$ such that the natural action $\theta_{\sssize P}$ of
$\Gamma$ on $P$ corresponds to $\theta$ through the $G$-bundle
trivialization $s: \bold{n}\times G @>{\cong}>> P$ induced by
$s_{\theta}$. Note that the $\Gamma$-action on $\bold{n}$ comes from
$\theta$ via the composition with a projection $\Gamma @>>>
G\!\wr\!\frak S_n @>>> \frak S_n=\text{Aut}(\bold{n})$. Such a
$\Gamma$-$G$ bundle $\pi_{\theta} :P @>>> \bold{n}$ together with
section $s_{\theta}$ was constructed in the proof of Theorem 3-1 in
two different ways. Let $\Cal S[P\times_GM]$ be the set of sections of
the fibre bundle $\pi_M: P\times_GM @>>>
\bold{n}$. This set can be identified with the set of $G$-equivariant
maps $\alpha\in\text{Map}_G(P,M)$ satisfying
$\alpha(pg)=g^{-1}\alpha(p)$ for all $g\in G$ and $p\in P$. A
$G$-automorphism $\varphi\in\text{Aut}_G(P)$ acts on
$\alpha\in\text{Map}_G(P,M)$ by
$\varphi(\alpha)=(\varphi^{-1})^*(\alpha)=\alpha\circ\varphi^{-1}$. Thus,
through $\theta_{\sssize P}$, $u\in\Gamma$ acts on $\text{Map}_G(P,M)=
\Cal S[P\times_GM]$ as $\theta_{\sssize P}(u^{-1})^*$. On the other
hand, we have an identification $\text{Map}_G(\bold{n}\times G,M)
@>{\cong}>> M^n$ by evaluation at $1\in G$. Since $u\in \Gamma$ acts
on $\bold{n}\times G$ as $\theta(u)$, the action of $u\in \Gamma$ on
$\text{Map}_G(\bold{n}\times G,M)\cong M^n$ is given by
$\theta(u^{-1})^*$. As we remarked earlier after Lemma 3-3, this
action of $\Gamma$ on $M^n$ is precisely the one given in (3-2).
Since $s:P @>>>\bold{n}\times G$ is $\Gamma$-equivariant with respect
to $\theta_{\sssize P}$ and $\theta$, we see that the $\Gamma$ action
on $\Cal S[P\times_GM]$ corresponds precisely to the $\Gamma$ action
on $M^n$.  The situation can be summarized by the following diagrams.
$$
\CD
\Cal S[P\times_GM]\cong \text{Map}_G(P,M)
@>{\cong}>{\theta_{\sssize P}(u^{-1})^*}> 
\text{Map}_G(P,M)\cong \Cal S[P\times_GM] \\
@V{\cong}V{s^*}V @V{\cong}V{s^*}V  \\
M^n\cong \text{Map}_G(\bold{n}\times G,M)
@>{\cong}>{\theta(u^{-1})^*=u\cdot}>
\text{Map}_G(\bold{n}\times G,M)\cong M^n
\endCD
$$
Here, horizontal arrows denote actions of $u\in\Gamma$.
From this, the first part of the next proposition is clear.

\proclaim{Proposition 5-1} Let $\theta:\Gamma @>>> G\!\wr\!\frak S_n$ be a
homomorphism. Let $\pi:P @>>> Z$ be a $\Gamma$-equivariant
$G$-principal bundle whose isomorphism class corresponds to the
conjugacy class $[\theta]$. Then we have
$(M^n)^{\langle\theta\rangle}\cong \Cal S[P\times_GM]^{\Gamma}$, where
the right hand side is the set of $\Gamma$-fixed sections. 

Let $\pi_{\rho}:P_{\rho}=\Gamma\times_{\rho}G @>>> \Gamma\!/\!H$ be an
irreducible $\Gamma$-$G$-bundle. Then we have 
$\Cal S[P_{\rho}\times_GM]^{\Gamma}\cong M^{\langle\rho\rangle}$. 
\endproclaim
\demo{Proof} We only have to prove the second part. Any section of
$P_{\rho}\times_G M @>>> \Gamma\!/\!H$ is represented by a
$G$-equivariant map $\alpha: P_{\rho} @>>> M$ satisfying
$\alpha(pg)=g^{-1}\alpha(p)$ for all $p\in P_{\rho}$ and $g\in
G$. This section is $\Gamma$-invariant if $\alpha(u^{-1}p)=\alpha(p)$
for all $u\in\Gamma$ and $p\in P_{\rho}$. Since $\Gamma$ acts
transitively on the base $\Gamma\!/\!H$, such a $\Gamma$-invariant
$G$-equivariant map $\alpha$ is uniquely determined by its value at
$p_0=[1,1]\in P_{\rho}$. This value $\alpha(p_0)\in M$ cannot be any
point in $M$. For any $h\in H\subset\Gamma$, using $\Gamma$-invariance
of $\alpha$, we must have $\alpha(p_0)=\alpha(hp_0)
=\alpha\bigl(p_0\rho(h)\bigr)=\rho(h)^{-1}\alpha(p_0)$. Thus,
$\alpha(p_0)$ must belong to $\rho(H)$-invariant subset of $M$. This
gives us an injective map $\Cal S[P_{\rho}\times_{G}M]^{\Gamma} @>>>
M^{\langle\rho\rangle}$.

To see that this map is surjective, let $x\in M^{\langle\rho\rangle}$
be any $\langle\rho\rangle$-fixed point in $M$. Let a
$\Gamma$-invariant $G$-equivariant map $\alpha_x:P_{\rho} @>>> M$ be
defined by a formula $\alpha_x(u^{-1}p_0g)=g^{-1}x$ for all
$u\in\Gamma$ and $g\in G$. To see that this well defined, let
$u_1^{-1}p_0g_1=u_2^{-1}p_0g_2$. Then there exists a unique $h\in H$
such that $u_2=hu_1$ and $g_2=\rho(h)g_1$. Then
$g_2^{-1}x=g_1^{-1}\rho(h)^{-1}x=g_1^{-1}x$, since $x\in
M^{\langle\rho\rangle}$. This proves that $\alpha_x$ is well defined
and the above correspondence is a bijection. This completes the proof.
\qed
\enddemo

From this proposition, the next corollary is straightforward.

\proclaim{Corollary 5-2} Let $r(H,\rho)$ be the number of isomorphic
copies of $\pi_{\rho}:P_{\rho} @>>> \Gamma\!/\!H$ appearing in the
irreducible decomposition of a $\Gamma$-equivariant $G$-principal
bundle $\pi: P @>>> Z$. Then
$$
\Cal S[P\times_GM]^{\Gamma}\cong
\prod_{[H]}\prod_{[\rho]}
\bigl(M^{\langle\rho\rangle}\bigr)^{r(H,\rho)}.
\tag5-2
$$
Furthermore, the action of $\text{\rm Aut}_{\Gamma\text{-}G}(P)$ on
the above set respects $[H]$-$[\rho]$ product decomposition described
in \rom{(4-7)}. 
\endproclaim
\demo{Proof} Let $[P @>>> Z]=\coprod\limits_{[H]}\coprod\limits_{[\rho]}
\!\!\coprod\limits^{\sssize r(H,\rho)}\!\![P_{\rho} @>>> \Gamma/H]$ 
be the irreducible decomposition of $P @>>> Z$.  Since the
$\Gamma$-action respects the decomposition of $P@>>> Z$ into
irreducible bundles, we have
$$
\Cal S[P\times_GM]^{\Gamma}\cong
\prod_{[H]}\prod_{[\rho]}\!\!\prod^{r(H,\rho)}\!\!
\Cal S[P_{\rho}\times_GM]^{\Gamma}.
$$
From the second part of Proposition 5-1, the isomorphism (5-2)
follows. 

Since there are no $\Gamma$-$G$ isomorphisms between $P_{\rho} @>>>
\Gamma\!/\!H$ and $P_{\rho'} @>>> \Gamma\!/\!H'$ if
$([H],[\rho])\not=([H'],[\rho'])$ by Classification Theorem 3-6, the
action of the $\Gamma$-$G$ automorphism group
$\text{Aut}_{\Gamma\text{-}G}(P)$ respects $[H]$-$[\rho]$
decomposition. This completes the proof.
\qed
\enddemo
 
Next we describe the action of
$\text{Aut}_{\Gamma\text{-}G}(P_{\rho})$ on $M^{\langle\rho\rangle}$.
By Theorem 4-4, an arbitrary element in
$\text{Aut}_{\Gamma\text{-}G}(P_{\rho})$ is of the form
$\varphi_{(u,g)}$ for some $(u,g)\in T_{\rho}$. This pair $(u,g)$ is
unique up to the left action by $H$, in view of Theorem 4-4.

\proclaim{Proposition 5-3} The left action of 
$\varphi_{(u,g)}\in\text{\rm Aut}_{\Gamma\text{-}G}(P_{\rho})$ on 
$x\in M^{\langle\rho\rangle}\cong \Cal S[P_{\rho}\times_GM]^{\Gamma}$
is given by $\varphi_{(u,g)}(x)=gx$. 
\endproclaim
\demo{Proof} By definition, the left action of
$\varphi\in\text{Aut}_{\Gamma\text{-}G}(P_{\rho})$ on
$\alpha\in\text{Map}_G(P_{\rho},M)\cong \Cal S[P_{\rho}\times_GM]$ is
given by the composition
$\varphi(\alpha)=\alpha\circ\varphi^{-1}$. For any $x\in
M^{\langle\rho\rangle}$, let $\alpha_x:P_{\rho} @>>> M$ be the
$\Gamma$-invariant $G$-equivariant map such that
$\alpha_x(p_0)=x$. When $(u,g)\in T_{\rho}$, we have
$g^{-1}\rho(uhu^{-1})g=\rho(h)$ for any $h\in H$. Since $u\in
N_{\Gamma}(H)$, we have $uhu^{-1}\in H$. Since the correspondence
$\text{Map}_G(P_{\rho},M)\cong M^{\langle\rho\rangle}$ is given by
evaluation at $p_0\in P_{\rho}$, we evaluate
$\varphi_{(u,g)}(\alpha_x) =\alpha_x\circ \varphi_{(u,g)}^{-1}$ at
$p_0$. Here $\Gamma$-$G$ map $\varphi_{(u,g)}$ is characterized by
$\varphi_{(u,g)}(p_0)=u^{-1}p_0g$. Thus
$\varphi_{(u,g)}^{-1}(p_0)=up_0g^{-1}$. Now,
$\alpha_x\circ\varphi_{(u,g)}^{-1}(p_0)=\alpha_x(up_0g^{-1})
=g\alpha_x(p_0)=gx$. Hence the action of $\varphi_{(u,g)}$ on $x$ is
given by $\varphi_{(u,g)}(x)=gx$. Indeed, $gx\in
M^{\langle\rho\rangle}$ because for any $h\in H$,
$\rho(h)(gx)=g\bigl(g^{-1}\rho(h)g\bigr)x=g\rho(u^{-1}hu)x=gx$, since
$x\in M^{\langle\rho\rangle}$.  This completes the proof.
\qed
\enddemo

Note that when $(u,g)\in T_{\rho}$, the element $g\in G$ is not
necessarily in $C(\rho)$, rather it is in the normalizer
$N_{G}^{}\bigl(\rho(H)\bigr)$, which can be seen from the definition
(4-4) of $T_{\rho}$. In general, the subgroup $\{g\in G \mid (u,g)\in
T_{\rho}\}$ is a proper subgroup of $N_G^{}\bigl(\rho(H)\bigr)$. The
last part of the above proof says that $N_{G}^{}\bigl(\rho(H)\bigr)$
acts on the $\langle\rho\rangle$-fixed point subset
$M^{\langle\rho\rangle}$.

We now calculate the generating function of the
invariants $\varphi_{\sssize\Gamma}(M^n;G_n)$ given in (5-1) for
$n\ge0$. An actual explicit concise formula depends on individual
properties of $\varphi$. However, the next formula applies to any
$\varphi$.

\proclaim{Proposition 5-4} Let $\varphi(M;G)$ be a multiplicative
orbifold invariant of a $G$-manifold $M$. Then 
$$
\sum_{n\ge0}q^n\varphi_{\sssize\Gamma}(M^n;G\!\wr\!\frak S_n)
=\prod_{[H]}\prod_{[\rho]}\Big[\sum_{r\ge0}q^{|\Gamma\!/\!H|r}
\varphi\bigl((M^{\langle\rho\rangle})^r;
\text{\rm Aut}_{\Gamma\text{-}G}(P_{\rho})\!\wr\!\frak S_r\bigr)\Bigr], 
\tag5-3
$$
where $[H]$ runs over the set of all the conjugacy classes of
subgroups of $\Gamma$, and for each given $[H]$, $[\rho]$ runs over
the set $\text{\rm Hom}(H,G)/(N_{\Gamma}(H)\times G)$.
\endproclaim
\demo{Proof} First, we rewrite the right hand side of (5-1) in terms
of $\Gamma$-equivariant $G$-bundles. By Theorem 3-1, each conjugacy
class $[\theta]\in\text{Hom}(\Gamma,G_n)/G_n$ corresponds to an
isomorphism class of $\Gamma$-equivariant $G$-principal bundle $P$
over a $\Gamma$-set $Z$ of order $n$. Furthermore, using Lemma 4-1,
Theorem 4-5, Proposition 5-1, and Corollary 5-2, we have
$$
\align
\varphi_{\sssize \Gamma}(M^n;G_n)
&=\!\!\!\!\!\!\sum\Sb [P] \\\ \ \ |Z|=n\endSb 
\!\!\!\!\!\!
\varphi\bigl(\Cal S[P\times_GM]^{\Gamma};
\text{Aut}_{\Gamma\text{-}G}(P)\bigr) \\ 
&=\!\!\!\!\!\!\sum_{\{r(H,\rho)\}\ \ }\!\!\!\!
\prod_{[H]}\prod_{[\rho]}
\varphi\bigl((M^{\langle\rho\rangle})^{r(H,\rho)};
\text{Aut}_{\Gamma\text{-}G}(P_{\rho})\!\wr\!\frak S_{r(H,\rho)}\bigr),
\endalign
$$
where the second summation is over all possible sets of non-negative
integers $\{r(H,\rho)\}_{[H],[\rho]}$ such that
$\sum_{[H],[\rho]}|\Gamma\!/\!H|r(H,\rho)=n$.  Taking the summation
over all $n\ge0$, we have
$$
\align
\sum_{n\ge0}q^n\varphi_{\sssize\Gamma}(M^n;G_n)
&=\!\!\!\!\!\!\!\!\!\sum_{r(H,\rho)\ge0\ \ }\!\!\!\!
\prod_{[H]}\prod_{[\rho]}
q^{|\Gamma\!/\!H|r(H,\rho)}
\varphi\bigl((M^{\langle\rho\rangle})^{r(H,\rho)};
\text{Aut}_{\Gamma\text{-}G}(P_{\rho})\!\wr\!\frak S_{r(H,\rho)}\bigr) \\
&=\prod_{[H]}\prod_{[\rho]}
\Big[\sum_{r\ge0}q^{|\Gamma\!/\!H|r}
\varphi\bigl((M^{\langle\rho\rangle})^r;
\text{\rm Aut}_{\Gamma\text{-}G}(P_{\rho})\!\wr\!\frak S_r\bigr)\Bigr].
\endalign
$$
This completes the proof. 
\qed
\enddemo

Recall that we have two kinds of orbifold Euler characteristics given
in (1-3). We apply formula (5-3) with $\varphi(M;G)
=\chi^{\text{orb}}_{\sssize\Gamma}(M;G)$, for various $\Gamma$ arising
as fundamental groups of real 2-dimensional surfaces which are
orientable or non-orientable, compact or non-compact. In the next
section, we apply formula (5-3) to $\chi_{\sssize\Gamma}(M;G)$
[Theorem 6-3]. We will see that these formulae are closely related but
very different.

\proclaim{Theorem 5-5} Let $G$ be a finite group, and let $M$ be a
$G$-manifold. For any group $\Gamma$, we have 
$$
\sum_{n\ge0}q^n\chi^{\text{orb}}_{\sssize\Gamma}(M^n;G\!\wr\!\frak S_n)
=\exp\Bigl[\sum_{r\ge1}\frac{q^r}r 
\Bigl\{\!\!\!\!\!\!\sum\Sb H \\ \ \ |\Gamma\!/\!H|=r \endSb 
\!\!\!\!\!\!
\chi_{\sssize H}^{\text{orb}}(M;G)\Bigr\}
\Bigr],
\tag5-4
$$
where the second summation on the right hand side runs over all index
$r$ subgroups of $\Gamma$. 
\endproclaim
\demo{Proof} We apply formula (5-3) with 
$\varphi=\chi^{\text{orb}}.$ By definition of
$\chi^{\text{orb}}(M;G)$ in (1-3), we have 
$$
\align
\sum_{n\ge0}q^n\chi^{\text{orb}}_{\Gamma}(M^n;G\!\wr\!\frak S_n)
&=\prod_{[H]}\prod_{[\rho]}\Big[\sum_{r\ge0}q^{|\Gamma\!/\!H|r}
\chi^{\text{orb}}\bigl((M^{\langle\rho\rangle})^r;
\text{\rm Aut}_{\Gamma\text{-}G}(P_{\rho})\!\wr\!\frak S_r\bigr)\Bigr] \\
&=\prod_{[H]}\prod_{[\rho]}
\bigg\{\sum_{r\ge0}q^{|\Gamma\!/\!H|r}\frac{\chi(M^{\langle\rho\rangle})^r}
{|\text{\rm Aut}_{\Gamma\text{-}G}(P_{\rho})|^rr!}\biggr\} \\
&=\prod_{[H]}\prod_{[\rho]}
\exp\biggl[q^{|\Gamma\!/\!H|}\frac{\chi(M^{\langle\rho\rangle})}
{|\text{\rm Aut}_{\Gamma\text{-}G}(P_{\rho})|}\biggr] \\
&=\exp\biggl[\sum_{r\ge1}q^r\!\!\!\!\!\!\!\!
\sum\Sb [H] \\ \ \ |\Gamma\!/\!H|=r \endSb 
\!\!\!\!\!\!\sum_{[\rho]}
\frac{\chi(M^{\langle\rho\rangle})}
{|\text{\rm Aut}_{\Gamma\text{-}G}(P_{\rho})|}\biggr].
\endalign
$$
In the last formula, $[H]$ runs over all conjugacy classes of index
$r$ subgroups of $\Gamma$, and for a given $[H]$, $[\rho]$ runs over
the set $\text{Hom}(H,G)/(N_{\Gamma}(H)\times G)$. 

Since the centralizer $C(\rho)$ is the isotropy subgroup
of the conjugation action by $G$ on $\text{Hom}(H,G)$ at $\rho$, we have
$\#(\rho)\cdot|C(\rho)|=|G|$, where $(\rho)$ is the $G$-conjugacy
class of $\rho$. Similarly, since $N^{\rho}_{\sssize\Gamma}(H)/H$ is
the isotropy subgroup of the $N_{\sssize\Gamma}(H)/H$-action on the
set of $G$-conjugacy classes $\text{Hom}(H,G)/G$ at $(\rho)$, the
length of the $N_{\sssize\Gamma}(H)/H$-orbit through $(\rho)$ is
$|N_{\sssize\Gamma}(H)/N_{\sssize\Gamma}^{\rho}(H)|$. Now we continue
our calculation. By Theorem 4-2, $|\text{\rm Aut}_{\Gamma\text{-}G}
(P_{\rho})|=|C(\rho)||N_{\sssize\Gamma}^{\rho}(H)/H|$. Thus, for a
fixed $[H]$, we have
$$
\multline
\!\!\!\!\!\!\!\!
\sum_{[\rho]}\frac{\chi(M^{\langle\rho\rangle})}
{|\text{\rm Aut}_{\Gamma\text{-}G}(P_{\rho})|}
=\sum_{[\rho]}\frac{\chi(M^{\langle\rho\rangle})}
{|C(\rho)||N_{\sssize\Gamma}^{\rho}(H)/H|} 
=\sum_{(\rho)}\frac{|N_{\sssize\Gamma}^{\rho}(H)/H|}
{|N_{\sssize\Gamma}(H)/H|}
\frac{\chi(M^{\langle\rho\rangle})}
{|C(\rho)||N_{\sssize\Gamma}^{\rho}(H)/H|}  \\
=\frac1{|N_{\sssize\Gamma}(H)/H|}
\sum_{(\rho)}\frac{\chi(M^{\langle\rho\rangle})}{|C(\rho)|}
=\frac1{|N_{\sssize\Gamma}(H)/H|}
\sum_{\rho}\frac{\chi(M^{\langle\rho\rangle})}
{\#(\rho)\cdot|C(\rho)|} \\
=\frac1{|N_{\sssize\Gamma}(H)/H|}\frac1{|G|}
\sum_{\rho}\chi(M^{\langle\rho\rangle})
=\frac1{|N_{\sssize\Gamma}(H)/H|}
\chi^{\text{orb}}_{\sssize H}(M;G).
\endmultline
$$
Here, $[\rho]$ runs over the set $\text{Hom}(H,G)/
(N_{\Gamma}(H)\times H)$, $(\rho)$ runs over the set of $G$-conjugacy
classes $\text{Hom}(H,G)/G$, and $\rho$ runs over the set of all
homomorphisms $\text{Hom}(H,G)$. The last equality is due to formula
(2-4). 

We convert the summation over the conjugacy classes $[H]$ of
index $r$ subgroups to a summation over index $r$ subgroups $H$.
Since there are $|\Gamma/N_{\Gamma}(H)|$ elements in the
$\Gamma$-conjugacy class $[H]$, 
$$
\align
\sum\Sb [H] \\ \ \ |\Gamma\!/\!H|=r\endSb 
\!\!\!\!\!\!\!\!
\frac1{|N_{\Gamma}(H)/H|}\chi^{\text{orb}}_{\sssize H}(M;G)
&=\!\!\!\!\!\!\!\!\sum\Sb H \\ \ \ |\Gamma\!/\!H|=r \endSb 
\!\!\!\!\!\!\!\!
\frac1{|\Gamma/N_{\Gamma}(H)||N_{\Gamma}(H)/H|}
\chi^{\text{orb}}_{\sssize H}(M;G) \\
&=\!\!\!\!\!\!\!\!\sum\Sb H \\ \ \ |\Gamma\!/\!H|=r \endSb 
\!\!\!\!\!\!\!\!
\frac1{|\Gamma\!/\!H|}
\chi^{\text{orb}}_{\sssize H}(M;G)
=\frac1r\!\!\!\!\!\!\sum\Sb H \\ \ \ |\Gamma\!/\!H|=r \endSb
\!\!\!\!\!\!
\chi^{\text{orb}}_{\sssize H}(M;G).
\endalign
$$
This completes the proof of formula (5-4). 
\qed
\enddemo

The above formula implies that to calculate the orbifold Euler
characteristic $\chi^{\text{orb}}_{\sssize\Gamma}(M^n;G\!\wr\!\frak S_n)$
of the $n$-fold symmetric product of an orbifold, we need to know
$\chi^{\text{orb}}_{\sssize H}(M;G)$ for every subgroup $H$ of $\Gamma$
whose index is at most $n$.

We specialize this formula. Letting $M=\text{pt}$, we obtain the
following combinatorial formula:
$$
\sum_{n\ge0}q^n\frac{|\text{Hom}(\Gamma,G\!\wr\!\frak S_n)|}
{|G|^nn!}
=\exp\biggl[\sum_{r\ge1}\frac{q^r}r
\Bigl\{\!\!\!\!\!\sum\Sb H \\ \ \ |\Gamma\!/\!H|=r \endSb 
\!\!\!\!\!\!
\frac{|\text{Hom}(H,G)|}{|G|}\Bigr\}\biggr].
\tag5-5
$$
This formula can be used to compute the number of homomorphisms into
wreath products. 

Next, we let $G$ be a trivial group in (5-4). Then we
get
$$
\sum_{n\ge0}q^n\chi^{\text{orb}}_{\sssize\Gamma}(M^n;\frak S_n)
=\Bigl(\exp\Bigl[\sum_{r\ge1}\frac{q^r}rj_r(\Gamma)
\Bigr]\Bigr)^{\chi(M)},
\tag5-6
$$
where $j_r(\Gamma)$ is the number of index $r$ subgroups of
$\Gamma$. This formula shows that when $M$ varies,
$\chi^{\text{orb}}_{\sssize\Gamma}(M^n;\frak S_n)$ depends only on the
Euler characteristic $\chi(M)$. Further specializing this formula to
the case $M=\text{pt}$, we get a well-known combinatorial formula
[\Cite{St}, p.76]:
$$
\sum_{n\ge0}q^n\frac{|\text{Hom}(\Gamma,\frak S_n)|}{n!}
=\exp\Bigl[\sum_{r\ge1}\frac{q^r}rj_r(\Gamma)\Bigr].
\tag5-7
$$
A related combinatorial formula can be found in physics
literature in the context of partition functions of permutation
orbifolds \cite{Ba}. 

We apply formula (5-4) to various groups $\Gamma$ to deduce
numerous consequences. 

\specialhead
(i)\qua Higher order ($p$-primary) orbifold Euler characteristic
\endspecialhead

These are orbifold Euler characteristics associated to $\Gamma=\Bbb
Z^d$ or $\Gamma=\Bbb Z_p^d$, where $\Bbb Z_p$ is the ring of $p$-adic
integers. The case $\Gamma=\Bbb Z^d$ corresponds to the case in which
the manifold $\Sigma$ is a $d$-dimensional torus $T^d$ in our
consideration of twisted sectors (2-11) and (2-12). 

We need to know the number $j_r(\Gamma)$ of index $r$ subgroups for
these groups $\Gamma$. This is well known, and we discussed these
numbers in our previous paper [\Cite{T},\,Lemma 4-4, Lemma 5-5].

\proclaim{Lemma 5-6} \rom{(1)}\qua For any $r\ge1$ and $d\ge1$, we have 
$$
j_r(\Bbb Z^d)=\!\!\!\!\!\!\!\!\!\!\!\!
\sum_{\ \ \ r_1r_2\cdots r_d=r}
\!\!\!\!\!\!\!\!\!\!\!\!r_2r_3^2\cdots r_d^{d-1},
\quad\text{and}\quad
j_r(\Bbb Z^d)=\sum_{m|r}m\cdot j_m(\Bbb Z^{d-1}). 
$$

\rom{(2)}\qua For any $r\ge0$ and $d\ge1$, we have 
$$
\botsmash{
j_{p^r}(\Bbb Z^d_p)=\!\!\!\!\!\!\!\!\!\!\!\!\!\!\!\!\!
\sum_{\ \ \ \ell_1+\ell_2+\cdots+\ell_d=r}
\!\!\!\!\!\!\!\!\!\!\!\!\!\!\!\!
p^{\ell_2}p^{2\ell_3}\cdots p^{(d-1)\ell_d},
\quad\text{and}\quad
j_{p^r}(\Bbb Z^d_p)=\!\!\!\sum_{0\le\ell\le r}\!\!\!
p^{\ell}\cdot j_{p^{\ell}}(\Bbb Z^{d-1}_p).
}
$$
\endproclaim

One of the main results in \cite{T} was the {\it inductive}
calculation of the generating function of
$\chi_{\sssize{\Gamma}}(M^n;G\!\wr\!\frak S_n)$ for $\Gamma=\Bbb Z^d,
\Bbb Z^d_p$. (See Theorem 6-5.) Now with our general formula (5-4),
the proof of the formula of generating functions for
$\{\chi_{\sssize{\Gamma}}^{\text{orb}}(M^n;G\!\wr\!\frak
S_n)\}_{n\ge0}$ with $\Gamma=\Bbb Z^d, \Bbb Z^d_p$ is a
straightforward corollary. Note that our proof of (5-4) was by a
direct proof, not inductive proof.

\proclaim{Theorem 5-7} Let $G$ be a finite group and let $M$ be a
$G$-manifold. For $d\ge1$,
$$
\aligned
\sum_{n\ge0}q^n\chi^{\text{orb}}_{\Bbb Z^d}(M^n;G\!\wr\!\frak S_n)
&=\Bigl[\prod_{r\ge1}(1-q^r)^{-j_r(\Bbb Z^{d-1})}\Bigr]
^{\chi^{\text{orb}}_{\Bbb Z^d}(M;G)} \\
\sum_{n\ge0}q^n\chi^{\text{orb}}_{\Bbb Z_p^d}(M^n;G\!\wr\!\frak S_n)
&=\Bigl[\prod_{r\ge0}(1-q^{p^r})^{-j_{p^r}(\Bbb Z^{d-1}_p)}\Bigr]
^{\chi^{\text{orb}}_{\Bbb Z_p^d}(M;G)}.
\endaligned
\tag5-8
$$
\endproclaim
\demo{Proof} First note that for any $r\ge1$, every  index $r$
subgroup of $\Bbb Z^d$ is isomorphic to $\Bbb Z^d$, and any subgroup
of $\Bbb Z^d_p$ has index a power of $p$ and any index $p^r$ subgroup
of $\Bbb Z^d_p$ is again isomorphic to $\Bbb Z^d_p$. Thus, by (5-4), 
$$
\aligned
\sum_{n\ge0}q^n\chi^{\text{orb}}_{\sssize\Bbb Z^d}(M^n;G\!\wr\!\frak S_n)
&=\exp\Bigl[\sum_{r\ge1}\frac{q^r}{r}j_r(\Bbb Z^d)
\chi^{\text{orb}}_{\sssize\Bbb Z^d}(M;G)\Bigr] \\
&=\Bigl(\exp\Bigl[\sum_{r\ge1}\frac{q^r}{r}j_r(\Bbb Z^d)\Bigr]\Bigr)
^{\chi^{\text{orb}}_{\sssize\Bbb Z^d}(M;G)}.
\endaligned
$$
For the remaining part, using the inductive formula in Lemma
5-6, we have
$$
\multline
\exp\Bigl[\sum_{r\ge1}\frac{q^r}rj_r(\Bbb Z^d)\Bigr]
=\exp\Bigl[\sum_{r\ge1}\frac{q^r}r
\sum_{m|r}m\cdot j_m(\Bbb Z^{d-1})\Bigr]  \\
=\exp\Bigl[\sum_{m\ge1}j_m(\Bbb Z^{d-1})
\Bigl(\sum_{\ell\ge1}\frac{q^{m\ell}}{\ell}\Bigr)\Bigr] 
=\prod_{m\ge1}(1-q^m)^{-j_m(\Bbb Z^{d-1})}.
\endmultline
$$
Similarly for the $p$-adic case. This completes the proof. 
\qed
\enddemo

The $d=0$ case in the above is elementary, and the formula was given in
(2-10).

\specialhead
(ii)\qua Higher genus orbifold Euler characteristic (orientable case)
\endspecialhead

Let $\Sigma_{g+1}$ be a real $2$-dimensional closed orientable genus
$g+1$ surface. Its fundamental group is the surface group
$\Gamma_{\!g+1}$ described in (1-7). By covering space theory,
conjugacy classes of index $r$ subgroups of $\Gamma_{\!g+1}$ are in
$1:1$ correspondence with the isomorphism classes of $r$-fold covering
spaces over $\Sigma_{g+1}$. By Hurwitz's Theorem, or by a simple
argument using Euler characteristic, we see that any such $r$-fold
covering space of $\Sigma_{g+1}$ is a closed orientable surface of
genus $rg+1$. Hence any index $r$ subgroup of $\Gamma_{\!g+1}$ is
always isomorphic to $\Gamma_{rg+1}$, although they may sit inside of
the group $\Gamma_{\!g+1}$ differently, and their conjugacy classes
can be different. A direct application of Theorem 5-5 gives the
next theorem.

\proclaim{Theorem 5-8} Let $g\ge0$. The generating function of the
genus $g+1$ orbifold Euler characteristic of the $n$-th
symmetric orbifold is given by
$$
\sum_{n\ge0}q^n\chi^{\text{orb}}_{\sssize\Gamma_{\!g+1}}
(M^n;G\!\wr\!\frak S_n)
=\exp\Bigl[\sum_{r\ge1}q^r\Bigl\{\frac{j_r(\Gamma_{\!g+1})}{r}
\chi^{\text{orb}}_{\sssize\Gamma_{rg+1}}(M;G)\Bigr\}\Bigr].
\tag5-9
$$
\endproclaim

The above formula requires that we know the numbers
$j_r(\Gamma_{\!g+1})$ for $g\ge0$. By letting $M=\text{pt}$, we get the
first combinatorial formula in (1-14). Further letting $G$ be trivial,
we get 
$$
\sum_{n\ge0}q^n\frac{\text{Hom}(\Gamma_{\!g+1},\frak S_n)|}{n!}
=\exp\Bigl[\sum_{r\ge1}\frac{j_r(\Gamma_{\!g+1})}{r}q^r\Bigr].
\tag5-10
$$
The left hand side was given in (1-15) in terms of the character
theory of $\frak S_n$. Thus, $j_r(\Gamma_{\!g+1})$ is, in principle,
calculable.

\specialhead
(iii)\qua Orbifold Euler characteristic associated to free groups
\endspecialhead

Let $\Gamma=F_{s+1}$ be a free group generated by $s+1$ elements,
$s\ge0$. Any index $r$ subgroup is isomorphic to $F_{rs+1}$, a free
group with $rs+1$ generators. This can be most easily seen
geometrically using covering spaces, as follows. Let $S$ be a torus
with $s$ disjointly embedded discs removed with $s\ge1$. Then
$\pi_1(S)\cong F_{s+1}$. Any conjugacy class of an index $r$ subgroup
corresponds to an isomorphism class of an $r$-fold covering space
$\tilde{S}$ which is a torus with $rs$ disjointly embedded discs
removed. Hence $\pi_1(\tilde{S})\cong F_{rs+1}$. Thus, any index $r$
subgroup of $F_{s+1}$ is again a free group isomorphic to
$F_{rs+1}$. Note that any two index $r$ subgroups of $F_{s+1}$ are not
necessarily conjugate to each other, although they are abstractly
isomorphic. Again applying Theorem 5-5 directly, we get

\proclaim{Theorem 5-9} Let $s\ge0$. Then 
$$
\sum_{n\ge0}q^n\chi^{\text{orb}}_{\sssize F_{s+1}}(M^n;G\!\wr\!\frak S_n)
=\exp\Bigl[\sum_{r\ge1}\Bigl\{\frac{j_r(F_{s+1})}{r}
\chi^{\text{orb}}_{\sssize F_{rs+1}}(M;G)\Bigr\}q^r\Bigr].
\tag5-11
$$
\endproclaim

Now we let $M=\text{pt}$ in the above formula. Since
$\text{Hom}(F_{r},G)\cong G^r$, we have a following combinatorial
formula
$$
\sum_{n\ge0}q^n(|G|^nn!)^s
=\exp\Bigl[\sum_{r\ge1}q^r\frac{j_r(F_{s+1})}{r}|G|^{rs}\Bigr].
\tag5-12
$$
How do we calculate $j_r(F_{s+1})$, the number of index $r$ subgroups
of a free group $F_{s+1}$? Letting $|G|=1$ in (5-12), we get
a well known formula [\Cite{St}, p.76]
$$
\sum_{n\ge0}q^n(n!)^s=\exp\Bigl[\sum_{r\ge1}q^r
\frac{j_r(F_{s+1})}{r}\Bigr].
\tag5-13
$$
This formula determines $j_r(F_{s+1})$ from known
quantities. Actually, the formulae (5-12) and (5-13) are easily seen
to be equivalent.

\specialhead
(iv)\qua  Higher genus orbifold Euler characteristic (non-orientable case) 
\endspecialhead

Let $N_{h+2}$ be a closed non-orientable (real 2-dimensional) surface
of genus $h+2$ for $h\ge0$. The fundamental group
$\Lambda_{h+2}=\pi_1(N_{h+2})$ is described by (1-9). We apply Theorem
5-5 with $\Gamma=\Lambda_{h+2}$. Since $N_1=\Bbb RP^2$ has fundamental
group $\Bbb Z/2$, this case directly follows from (5-4). Here we only
consider non-orientable surface groups of genus $2$ or greater.  First
we need to know more about subgroups of $\Lambda_{h+2}$. We discuss
this using covering spaces of $N_{h+2}$.

For $r\ge1$, by examining Euler characteristic, we can easily see that
any $r$-fold connected covering space of $N_{h+2}$ is either
non-orientable ($\cong N_{rh+2}$) or orientable ($\cong
\Sigma_{\frac{rh}2+1}$), and orientable covering spaces can occur only
when $r$ is even. We call subgroups of $\Lambda_{h+2}$ corresponding
to connected orientable covering spaces orientable subgroups. Let
$j_r(\Lambda_{h+2})^+$ denote the number of index $r$ orientable
subgroups of $\Lambda_{h+2}$. Similarly, let $j_r(\Lambda_{h+2})^-$
denote the number of index $r$ non-orientable subgroups. The surface
$N_{h+2}$ has a unique orientable double cover $\Sigma_{h+1}$, and any
orientable cover of $N_{h+2}$ is a cover of $\Sigma_{h+1}$. This means
that $j_{\text{odd}}(\Lambda_{h+2})^+=0$ and $j_{2r}(\Lambda_{h+2})^+
=j_r(\Gamma_{h+1})$ for all $r\ge1$ and $h\ge0$. Now we apply Theorem
5-5.

\proclaim{Theorem 5-10} Let $h\ge0$. Then
$$
\multline
\sum_{n\ge0}q^n\chi_{\sssize\Lambda_{h+2}}^{\text{orb}}
(M^n;G\!\wr\!\frak S_n)
=\topsmash{
\exp\Bigl[\sum_{r\ge1}\frac{q^r}{r}j_r(\Lambda_{h+2})
\chi^{\text{orb}}_{\sssize\Lambda_{h+2}}(M;G)
} \\
\smash{
+\sum_{r\ge1}\frac{q^{2r}}{2r}j_r(\Gamma_{h+1})
\bigl\{\chi^{\text{orb}}_{\sssize\Gamma_{rh+1}}(M;G)
-\chi^{\text{orb}}_{\sssize \Lambda_{2rh+2}}(M;G)\bigr\}\Bigr].
}
\endmultline
\tag5-14
$$
\endproclaim
\demo{Proof} Theorem 5-5 implies that 
$$
\multline
\botsmash{
\sum_{n\ge0}q^n\chi_{\sssize\Lambda_{h+2}}^{\text{orb}}
(M^n;G\!\wr\!\frak S_n)
}\\
=\exp\Bigl[\sum_{r\ge1}\frac{q^r}{r}
\bigl\{j_r(\Lambda_{h+2})^+
\chi^{\text{orb}}_{\sssize\Gamma_{\frac{rh}2+1}}(M;G) 
+j_r(\Lambda_{h+2})^-
\chi^{\text{orb}}_{\sssize \Lambda_{rh+2}}(M;G)\bigr\}\Bigr].
\endmultline
$$
We only have to rewrite this formula using
$j_r(\Lambda_{h+2})^-=j_r(\Lambda_{h+2})-j_r(\Lambda_{h+2})^+$, and a
fact that $j_r(\Lambda_{h+2})^+=j_{\frac r2}(\Gamma_{h+1})$ when $r$
is even and $j_r(\Lambda_{h+2})^+=0$ when $r$ is odd.
\qed
\enddemo

How do we calculate the number $j_r(\Lambda_{h+2})$? We let
$M=\text{pt}$ and let $G$ be a trivial group. Then the above formula
reduces to 
$$
\sum_{n\ge0}q^n\frac{|\text{Hom}(\Lambda_{h+2},\frak S_n)|}{n!}
=\exp\Bigl[\sum_{r\ge1}\frac{q^r}rj_r(\Lambda_{h+2})\Bigr].
\tag5-15
$$
By using the character theory of $\frak S_n$, one can show that for
all $h\ge-1$, 
$$
\frac{|\text{Hom}(\Lambda_{h+2},\frak S_n)|}{n!}
=\!\!\!\!\!\!\!\!\!\!\!
\sum_{[V]\in\text{Irred}(\frak S_n)\ }
\!\!\!\!\!\!\Bigl(\frac{n!}{\dim V}\Bigr)^{\!\!h}.
\tag5-16
$$
The above two formulae determine $j_r(\Lambda_{h+2})$. The formula
(5-16) is a special case of Kerber-Wagner Theorem \cite{KW}. Their
formula for a general finite group $G$, instead of $\frak S_n$, is 
$$
\frac{|\text{Hom}(\Lambda_{h+2},G)|}{|G|}
=\!\!\!\!\!\!\!\!
\sum_{\chi\in\text{Irred}(G)\ }
\!\!\!\!\!\!\Bigl(\frac{|G|}{\dim \chi}\Bigr)^{\!\!h}
\varepsilon_2(\chi)^{h+2},
\tag5-17
$$
where $\varepsilon_2(\chi)$ is the Schur indicator defined by 
$\varepsilon_2(\chi)=\sum_{g\in G}\chi(g^2)/|G|$. It is well known
that $\varepsilon_2(\chi)=1$ if $\chi$ is of real type, 
$\varepsilon_2(\chi)=0$ if $\chi$ is of complex type, and 
$\varepsilon_2(\chi)=-1$ if $\chi$ is of quaternionic type. For
details of this fact, see [\Cite{BtD},\,p.100]. It is well known that
any representation of $\frak S_n$ is of real type so that
$\varepsilon(\chi)=1$ for any irreducible character $\chi$ of $\frak
S_n$. 

\remark{Remark} The formulae (5-10) and (1-15), and formulae (5-15)
and (5-16) are very similar. By comparing these formulae it is
immediate that $j_r(\Gamma_{h+1})=j_r(\Lambda_{2h+2})$ for all $r\ge1$
and $h\ge0$. However, this does not necessarily mean that
$\chi^{\text{orb}}_{\sssize\Gamma_{\!g+1}}(M;G)
=\chi^{\text{orb}}_{\sssize\Lambda_{2g+2}}(M;G)$. In fact, they are in
general different, although in some cases they coincide. See (5-22)
for an example of such a case. 
\endremark

\specialhead
(v)\qua Klein bottle orbifold Euler characteristic
\endspecialhead

We specialize our formula (5-14) to $h=0$ case. The corresponding
non-orientable closed surface $N_2$ is the Klein bottle. By the remark
above, we have $j_r(\Gamma_1)=j_r(\Lambda_2)$ for all $r\ge1$. Since
by Theorem 5-6, $j_r(\Gamma_1)=j_r(\Bbb Z^2)=\sum_{\ell|r}\ell
=\sigma_1(r)$, the sum of positive divisors of $r$, the formula (5-14)
reduces to
$$
\multline
\sum_{n\ge0}q^n\chi_{\sssize\Lambda_{2}}^{\text{orb}}
(M^n;G\!\wr\!\frak S_n)\\
=\exp\Bigl[\botsmash{\sum_{r\ge1}}\frac{q^r}{r}\sigma_1(r)
\chi^{\text{orb}}_{\sssize\Lambda_{2}}(M;G) 
+\sum_{r\ge1}\topsmash{\frac{q^{2r}}{2r}}\sigma_1(r)
\bigl\{\chi^{\text{orb}}_{\sssize\Gamma_{1}}(M;G)
-\chi^{\text{orb}}_{\sssize \Lambda_{2}}(M;G)\bigr\}\Bigr].
\endmultline
$$
Rewriting exponentials in terms of infinite products as in the proof
of Theorem 5-7, we obtain

\proclaim{Theorem 5-11} Let $G$ be a finite group and let $M$ be a
$G$-manifold. Then 
$$
\sum_{n\ge0}q^n\chi^{\text{orb}}_{\sssize\Lambda_2}
(M^n;G\!\wr\!\frak S_n)
=\Bigl[\prod_{r\ge1}(1-q^{2r})\Bigr]^
{\frac{-1}2\chi^{\text{orb}}_{\Gamma_1}(M;G)}
\Bigl[\prod_{r\ge1}\Bigl(\frac{1+q^r}{1-q^r}\Bigr)\Bigr]^
{\frac12\chi^{\text{orb}}_{\Lambda_2}(M;G)}.
\tag5-18
$$
\endproclaim

Letting $M=\text{pt}$, we obtain a combinatorial formula for the
number of homomorphisms from the fundamental group of a Klein
bottle into a wreath product:
$$
\multline
\botsmash{
\sum_{n\ge0}q^n\frac{|\text{Hom}(\Lambda_2,G\!\wr\!\frak
S_n)|}{|G|^nn!}
}\\
=\Bigl[\prod_{\ell\ge1}(1-q^{2\ell})\Bigr]
^{-\frac12|\text{Hom}(\Gamma_1,G)|/|G|}
\biggl[\prod_{\ell\ge1}
\biggl(\frac{1+q^{\ell}}{1-q^{\ell}}\biggr)\biggr]
^{\frac12|\text{Hom}(\Lambda_2,G)|/|G|}. 
\endmultline
\tag5-19
$$
Letting $G$ be trivial in (5-18), we get 
$$
\sum_{n\ge0}q^n\chi^{\text{orb}}_{\sssize\Lambda_2}(M^n;\frak S_n)
=\Bigl[\prod_{\ell\ge1}(1-q^{\ell})\Bigr]^{-\chi(M)}.
\tag5-20
$$
On the other hand, when $d=2$ and $G$ is trivial, formula (5-8)
reduces to
$$
\sum_{n\ge0}q^n\chi_{\sssize\Bbb Z^2}^{\text{orb}}(M^n;\frak S_n)
=\bigl[\prod_{r\ge1}(1-q^r)\bigr]^{-\chi(M)}.
\tag5-21
$$
Comparing the above two formula, we see that
$\chi_{\sssize\Gamma_1}^{\text{orb}}(M^n;\frak S_n)
=\chi_{\sssize\Lambda_2}^{\text{orb}}(M^n;\frak S_n)$ for all $n\ge0$,
since $\Gamma_1=\Bbb Z^2$. This is no coincidence. In general, using
(5-6) and a fact that $j_r(\Gamma_{\!g+1})=j_r(\Lambda_{2g+2})$ for
$g\ge0$, $r\ge1$ (see the remark at the end of (iv) in this section),
we get
$$
\sum_{n\ge0}q^n\chi_{\sssize\Gamma_{\!g+1}}^{\text{orb}}
(M^n;\frak S_n)
=\sum_{n\ge0}q^n\chi_{\sssize\Lambda_{2g+2}}^{\text{orb}}
(M^n;\frak S_n), \qquad g\ge0.
\tag5-22
$$

\head
Orbifold invariants associated to $\Gamma$-sets and generating
functions:\strut\ infinite product formulae
\endhead

In section 5, we calculated the generating function
$\sum_{n\ge0}q^n\varphi_{\sssize\Gamma}(M^n;G\!\wr\!\frak S_n)$ for
$\varphi(M;G)=\chi^{\text{orb}}(M;G)$, where $G$ is a finite group,
$M$ is a $G$-manifold, and $\Gamma$ is any group. In this section, we
investigate the generating function with
$\varphi(M;G)=\chi(M;G)=\chi(M/G)$. The latter orbifold Euler
characteristic is more closely tied to geometry of $G$-action than the
former.

To describe the generating function
$\sum_{n\ge0}q^n\chi_{\sssize\Gamma}(M^n;G\!\wr\!\frak S_n)$, we need
to introduce a notion of generalized orbifold invariants associated to
transitive $\Gamma$-sets. Orbifold invariants associated to general
$\Gamma$-sets will be described later in this section in (6-13). Let
$H\subset\Gamma$ be a subgroup of finite index. Then the
$\varphi$-orbifold invariant associated to the $\Gamma$-isomorphism
class of a $\Gamma$-transitive set $\Gamma\!/\!H$ is defined by
$$
\varphi_{\sssize[\Gamma\!/\!H]}(M;G)\overset{\text{def}}\to=
\!\!\!\!\!\!\!\!\!\!\!\!\!\!\!\!\!\!\!\!\!\!\!\!\!\!\!\!\!\!\!\!\!\!\!\!
\sum_{\ \ \ \ \ \ \ \ \ \ [\rho]\in\text{Hom}(H,G)/(N_{\Gamma}(H)\times G)}
\!\!\!\!\!\!\!\!\!\!\!\!\!\!\!\!\!\!\!\!\!\!\!\!\!\!\!\!\!\!\!\!\!\!\!\!
\varphi\bigl(M^{\langle\rho\rangle};
\text{Aut}_{\Gamma\text{-}G}(P_{\rho})\bigr),
\tag6-1
$$
where $\pi_{\rho}:P_{\rho}=\Gamma\times_{\rho}G @>>> \Gamma\!/\!H$ is a
$\Gamma$-irreducible $G$-bundle, and $\text{Aut}_{\Gamma\text{-}G}
(P_{\rho})$ is the group of $\Gamma$-equivariant $G$-bundle
automorphisms of $P_{\rho}$ described in Theorem 4-2 and Theorem 4-4. 

When $H=\Gamma$, the conjugation on homomorphisms $\rho:\Gamma @>>> G$
by $N_{\Gamma}(\Gamma)=\Gamma$ is absorbed by the conjugation by $G$. 
Hence we have $\text{Hom}(\Gamma,G)/(N_{\Gamma}(\Gamma)\times G)=
\text{Hom}(\Gamma,G)/G$. Also in this case, 
$\text{Aut}_{\Gamma\text{-}G}(P_{\rho})$ reduces to $C_G(\rho)$: see
(4-2) for example. Hence $\varphi_{\sssize[\Gamma/\Gamma]}(M;G)
=\varphi_{\sssize\Gamma}(M;G)$, recovering the $\Gamma$-extended
orbifold invariant.

When $G$ is a trivial group $\{e\}$, $\rho$ is the trivial
homomorphism and $\text{Aut}_{\Gamma\text{-}G}(P_{\rho})$ reduces to
the group of $\Gamma$-automorphisms of the $\Gamma$-set $\Gamma\!/\!H$,
that is $H\backslash N_{\Gamma}(H)$. Thus,
$\varphi_{\sssize[\Gamma\!/\!H]}(M;\{e\}) =\varphi\bigl(M;H\backslash
N_{\Gamma}(H)\bigr)$, where $H\backslash N_{\Gamma}(H)$ acts trivially
on $M$ by Proposition 5-3.

The above definition (6-1) may seem rather unusual, but it is the
correct one. For example, the set over which $[\rho]$ runs has already
appeared in the classification theorem of $\Gamma$-irreducible
$G$-bundles in Theorem 3-6, and the above definition is closely
related to this theorem.

To explain the origin of the definition, we consider the map 
$$
\pi: \text{Hom}(\Gamma,G_n)/G_n @>>> 
\text{Hom}(\Gamma,\frak S_n)/\frak S_n,
$$
induced by the projection map $G_n @>>> \frak
S_n$. Since the set $\text{Hom}(\Gamma,\frak S_n)/\frak S_n$ can be
regarded as the set of isomorphism classes of $\Gamma$-sets of order
$n$, for any index $n$ subgroup $H$ of $\Gamma$, we may regard
$[\Gamma\!/\!H]\in \text{Hom}(\Gamma,\frak S_n)/\frak S_n$. The meaning of 
$\varphi_{\sssize[\Gamma\!/\!H]}(M;G)$ is clarified by the next
proposition. 

\proclaim{Proposition 6-1} For any subgroup $H$ of index $n$ in
$\Gamma$, we have 
$$
\varphi_{\sssize[\Gamma\!/\!H]}(M;G)
=\!\!\!\!\!\!\!\!\!\!\!\!\!\!\!\!
\sum_{\ \ \ \ [\theta]\in\pi^{-1}([\Gamma\!/\!H])}
\!\!\!\!\!\!\!\!\!\!\!\!\!\!\!\!
\varphi\bigl((M^n)^{\langle\theta\rangle};C_{G_n}(\theta)\bigr),
\tag6-2
$$
where $\pi:\text{\rm Hom}(\Gamma,G_n)/G_n @>>>
\text{\rm Hom}(\Gamma,\frak S_n)/\frak S_n$ is the obvious map. 
\endproclaim
\demo{Proof} We rewrite the quantity in the right hand side in terms of
$\Gamma$-equivariant $G$-bundles. By Theorem 3-1, any element in
$\pi^{-1}([\Gamma\!/\!H])$ is an isomorphism class of a
$\Gamma$-irreducible $G$-bundle over the $\Gamma$-set
$\Gamma\!/\!H$. By the Classification Theorem 3-6 this set is in 1:1
correspondence with $\text{Hom}(H,G)/(N_{\Gamma}(H)\times G)$, and the
isomorphism class of the $\Gamma$-irreducible $G$-bundle corresponding
to $[\theta]\in\pi^{-1}([\Gamma\!/\!H])$ is the isomorphism class of
$\pi_{\rho}:P_{\rho} @>>> \Gamma\!/\!H$, where $\rho:H @>>> G$ is
constructed from $\theta$ (see the proof of Theorem 3-1), essentially
by restriction to the isotropy subgroup $H\subset\Gamma$. Using Lemma
4-1 and Proposition 5-1, we have
$$
\sum_{\ \ \ [\theta]\in\pi^{-1}([\Gamma\!/\!H])}
\!\!\!\!\!\!\!\!\!\!\!\!\!\!\!\!
\varphi\bigl((M^n)^{\langle\theta\rangle};C_{G_n}(\theta)\bigr)
=\sum_{[\rho]}\varphi\bigl(M^{\langle\rho\rangle};
\text{Aut}_{\Gamma\text{-}G}(P_{\rho})\bigr),
$$
where $[\rho]$ runs over the set 
$\text{Hom}(H,G)/(N_{\Gamma}(H)\times G)$. Then formula (6-1) 
completes the proof. 
\qed
\enddemo

Note that $\varphi_{\sssize[\Gamma\!/\!H]}(M;G)$ is a partial sum of
the summation defining $\varphi_{\sssize\Gamma}(M^n;G_n)$ in
(5-1). The situation here will be clarified in (6-14) after we define
$\varphi_{\sssize[X]}(M;G)$ for general $\Gamma$-set $X$. 

We calculate an example of an orbifold invariant associated to a
$\Gamma$-set. We denote an orbifold invariant $\chi_{\sssize\Bbb
Z^d}(M;G)$ by $\chi^{\sssize(d)}(M;G)$. This is the $d$-th order
orbifold Euler characteristic of $(M;G)$ discussed in \cite{T}.

\proclaim{Lemma 6-2} Let $\Gamma$ be an abelian group, and let
$d\ge0$. Let $(M;G)$ be as before. Then for any subgroup $H$ of finite
index in $\Gamma$, we have
$$
\chi^{\sssize(d)}_{\sssize[\Gamma\!/\!H]}(M;G)
=|\Gamma\!/\!H|^d\cdot \chi^{\sssize(d)}_{\sssize H}(M;G).
\tag6-3
$$
In particular, $\chi_{\sssize[\Gamma\!/\!H]}(M;G)=\chi_{\sssize
H}(M;G)$, if $\Gamma$ is abelian. 
\endproclaim
\demo{Proof} Let $\rho:H @>>> G$. To calculate
$\chi^{(d)}_{\sssize[\Gamma\!/\!H]}(M;G)$ defined in (6-1), we need to
understand the group $\text{Aut}_{\Gamma\text{-}G}(P_{\rho})$.  The
group $T_{\rho}$ given in (4-4) is, in our case, isomorphic to
$\Gamma\times C_G(\rho)$ since $\Gamma$ is abelian. Then, by Theorem
4-4 we have $\text{Aut}_{\Gamma\text{-}G}(P_{\rho})\cong H\backslash
T_{\rho}\cong \Gamma\times_{\rho}C(\rho)$. From this, 
$|\text{Aut}_{\Gamma\text{-}G}(P_{\rho})|=|\Gamma\!/\!H||C(\rho)|$ and 
$$
\chi^{\sssize(d)}_{\sssize[\Gamma\!/\!H]}(M;G)
=\sum_{[\rho]}\chi^{(d)}\bigl(M^{\langle\rho\rangle};
\Gamma\times_{\rho}C(\rho)\bigr)
=\sum_{[\rho]}\frac1{|\Gamma\!/\!H||C(\rho)|}
\!\!\!\!\!\!\!\!\!\!\!\!\!\!\!\!\!\!\!\!\!\!\!\!\!\!\!\!\!
\sum_{\ \ \ \ \ \ \ \ \ \ \ \ \ \ 
\phi:\Bbb Z^{d+1} @>>> \Gamma\times_{\rho}C(\rho)}
\!\!\!\!\!\!\!\!\!\!\!\!\!\!\!\!\!\!\!\!\!\!\!\!\!\!\!\!\!\!
\chi\bigl((M^{\langle\rho\rangle})^{\langle\phi\rangle}\bigr),
$$
where the last equality is due to (2-4) and (2-5), and $[\rho]$ runs
over the set $\text{Hom}(H,G)/G$, where $H$ is abelian. By Proposition
5-3, the image of the natural map $\iota:\Gamma @>>>
\Gamma\times_{\rho}C(\rho)$ acts trivially on
$M^{\langle\rho\rangle}$. Let
$\{\gamma_i\}_{i\in\Gamma\!/\!H}\subset\Gamma$ be the representatives of
of the coset $\Gamma\!/\!H$. Then $\Gamma\times_{\rho}C(\rho)
=\coprod_{i\in\Gamma\!/\!H}\iota(\gamma_i)C(\rho)$. Note that two elements
$\iota(\gamma_{i_1})g_1$ and $\iota(\gamma_{i_2})g_2$ in
$\Gamma\times_{\rho}C(\rho)$, where $g_1, g_2\in C(\rho)$, commute if
and only if $g_1$ and $g_2$ commute in $C(\rho)$. Thus the above
summation becomes
$$
\align
\sum_{[\rho]}\frac1{|\Gamma\!/\!H||C(\rho)|} |\Gamma\!/\!H|^{d+1}
\!\!\!\!\!\!\!\!\!\!\!\!\!\!\!\!\!\!
\sum_{\ \ \ \ \ \ \phi':\Bbb Z^{d+1} @>>> C(\rho)}
\!\!\!\!\!\!\!\!\!\!\!\!\!\!\!\!\!\!
\chi\bigl((M^{\langle\rho\rangle})^{\langle\phi'\rangle}\bigr)
&=|\Gamma\!/\!H|^d\sum_{[\rho]}\chi^{\sssize(d)}
\bigl(M^{\langle\rho\rangle};C(\rho)\bigr)\\
&=|\Gamma\!/\!H|^d\chi^{\sssize(d)}_{\sssize H}(M;G). 
\endalign
$$
Letting $d=0$, we obtain the last statement. This completes the
proof.
\qed
\enddemo

Now we calculate the generating function
$\sum_{n\ge0}q^n\chi_{\sssize\Gamma}(M^n;G\!\wr\!\frak S_n)$. 

\proclaim{Theorem 6-3} Let $G$ be a finite group and let $M$ be a
$G$-manifold. For any group $\Gamma$, we have 
$$
\gathered
\sum_{n\ge0}q^n\chi_{\sssize\Gamma}(M^n;G\!\wr\!\frak S_n)
=\prod_{r\ge1}(1-q^r)^{\!\!-\!\!\!\!\sum\limits
_{\ \ \sssize{[H]_r}}\!\!\!\!
\chi_{\sssize[\Gamma\!/\!H]}(M;G)}, \\
\text{where}\quad 
\chi_{\sssize[\Gamma\!/\!H]}(M;G)
=\!\!\!\!\!\!\!\!\!\!\!\!\!\!\!\!\!\!\!\!\!\!\!\!\!\!\!\!\!\!\!\!
\sum_{\ \ \ \ \ \ [\rho]\in\text{\rm Hom}(H,G)\!/\!(N_{\Gamma}^{}(H)
\times G)}
\!\!\!\!\!\!\!\!\!\!\!\!\!\!\!\!\!\!\!\!\!\!\!\!\!\!\!\!\!
\chi\bigl(M^{\langle\rho\rangle}/
\text{\rm Aut}_{\Gamma\text{-}G}(P_{\rho})\bigr).
\endgathered
\tag6-4
$$
Here in the summation in the right hand side of the first identity,
$[H]_r$ means that $H$ runs over all the conjugacy classes of index
$r$ subgroups of $\Gamma$, and the action of $\text{\rm
Aut}_{\Gamma\text{-}G}(P_{\rho})$ on $M^{\langle\rho\rangle}$ in the
second formula is described in Proposition {\rm 5-3}.

When $\Gamma$ is an abelian group, 
$$
\sum_{n\ge0}q^n\chi_{\sssize\Gamma}(M^n;G\!\wr\!\frak S_n)
=\topsmash{
\prod_{r\ge1}(1-q^r)^{\!\!-\!\!\sum\limits_{\sssize{H_r}}
\!\chi_{\sssize H}(M;G)},
}
\tag6-5
$$
where in the summation in the right hand side, $H_r$ means that $H$
runs over all index $r$ subgroups in $\Gamma$.
\endproclaim
\demo{Proof} For (6-4), by Proposition 5-4 and Macdonald's
formula (2-9), we have
$$
\multline
\sum_{n\ge0}q^n\chi_{\sssize\Gamma}(M^n;G\!\wr\!\frak S_n)
=\prod_{[H]}\prod_{[\rho]}\Bigl[\sum_{r\ge0}q^{|\Gamma\!/\!H|r}
\chi\bigl((M^{\langle\rho\rangle})^r;
\text{\rm Aut}_{\Gamma\text{-}G}(P_{\rho})\!\wr\!\frak S_r\bigr)\Bigr] \\
=\prod_{[H]}\prod_{[\rho]}\Bigl[\sum_{r\ge0}q^{|\Gamma\!/\!H|r}
\chi\bigl(SP^r\bigl(M^{\langle\rho\rangle}/
\text{Aut}_{\Gamma\text{-}G}(P_{\rho})\bigr)\bigr)\\
=\prod_{[H]}\prod_{[\rho]}(1-q^{|\Gamma\!/\!H|})
^{-\chi(M^{\langle\rho\rangle}/
\text{Aut}_{\Gamma\text{-}G}(P_{\rho}))} 
=\topsmash{
\prod_{r\ge1}(1-q^r)
^{\!\!-\!\!\!\!\sum\limits_{\ [H]_r}
\!\!
\sum\limits_{\sssize[\rho]}
\chi(M^{\langle\rho\rangle};
\text{Aut}_{\Gamma\text{-}G}(P_{\rho}))},
}
\endmultline
$$
where in the last summation $[H]_r$ runs over all the conjugacy
classes of finite index in $\Gamma$, and $[\rho]$ runs over the set
$\text{Hom}(H,G)/(N_{\Gamma}(H)\times G)$. Notice that the summation
over $[\rho]$ in the exponent in the last line exactly gives
$\chi_{\sssize[\Gamma\!/\!H]}(M;G)$.

When $\Gamma$ is an abelian group, the formula (6-5) follows from
(6-4) in view of Lemma 6-2 with $d=0$. Since $\Gamma$ is abelian, each
subgroup $H$ is its own conjugacy class $[H]$. This completes the proof. 
\qed
\enddemo

Note that in the corresponding formula for
$\chi^{\text{orb}}_{\sssize\Gamma}$ given in (5-4), the summation
there was over subgroups $H$ rather than conjugacy classes $[H]$
used in (6-4).

\remark{Remark} From the above proof, one might think that the concept
of orbifold invariant associated to a finite $\Gamma$-set is just a
convenient symbol representing a complicated summation given in the
right hand side of (6-1). However, as we will see in \S7, this concept
has a definitely natural geometric origin analogous to twisted sectors
discussed in \S2. 
\endremark

\smallskip

Letting $M=\text{pt}$ or letting $G$ be trivial, we immediately have
the following corollary.

\proclaim{Corollary 6-4} \rom{(1)}\qua For any group $\Gamma$ and any
finite group $G$, we have
$$
\sum_{n\ge0}q^n|\text{\rm Hom}(\Gamma,G_n)/G_n|
=\prod_{r\ge1}(1-q^r)^{\!\!-\!\!\!\sum\limits_{\ [H]_r}
\!\!\!
|\text{\rm Hom}(H,G)/(N_{\Gamma}(H)\times G)|}.
\tag6-6
$$
Furthermore, when $\Gamma$ is abelian, we may omit the term
$N_{\Gamma}(H)$ from the above expression. 

\rom{(2)}\qua For any $\Gamma$ and any manifold $M$, we have 
$$
\sum_{n\ge0}q^n\chi_{\sssize\Gamma}(M^n;\frak S_n)
=\Bigl[\prod_{r\ge1}(1-q^r)^{-u_r(\Gamma)}\Bigr]^{\chi(M)},
\tag6-7
$$
where $u_r(\Gamma)$ is the number of conjugacy classes of index $r$
subgroups of $\Gamma$. 
\endproclaim

The proof for (1) is straightforward. For (2), all we need to note is
that when the group $G$ is trivial, (6-1) gives that
$\chi_{\sssize[\Gamma\!/\!H]}(M;\{e\})=\chi(M;H\backslash
N_{\Gamma}(H)) =\chi(M)$, since the group $H\backslash N_{\Gamma}(H)$
acts trivially on $M$ by Proposition 5-3.

If we further specialize (6-7) to the case $M=\text{pt}$ or if we let
$G=\{e\}$ in (6-6), then we
obtain the following well known formula [\Cite{St},\,p.76]:
$$
\sum_{n\ge0}q^n|\text{Hom}(\Gamma,\frak S_n)/\frak S_n|
=\prod_{r\ge1}(1-q^r)^{-u_r(\Gamma)}.
\tag6-{$6'$}
$$
Recall that $|\text{Hom}(\Gamma,\frak S_n)/\frak S_n|
=|\text{Hom}(\Gamma\times\Bbb Z,\frak S_n)|/n!$. The above formula can
be used to calculate the number of conjugacy classes of a given index
in the orientable surface group $\Gamma_{\!g+1}$of genus $g+1$, in the
free group $F_{s+1}$ on $s+1$ generators, and in the non-orientable
surface group $\Lambda_{h+2}$ of genus $h+2$. We discuss details in
\S8, where a general formula to calculate $u_r(\Gamma)$ is proved. 

We apply Theorem 6-3 to various groups $\Gamma$ to obtain various
results.

\specialhead
(i)\qua Higher order ($p$-primary) orbifold Euler characteristic of
symmetric orbifolds
\endspecialhead

Suppose $\Gamma=\Bbb Z^{d}$ or $\Bbb Z^d_p$ for $d\ge0$. Then the
corresponding $\Gamma$-extended orbifold Euler characteristic is what
we call higher order ($p$-primary) orbifold Euler characteristic
denoted by $\chi^{\sssize(d)}(M;G)$ or $\chi^{\sssize(d)}_{\sssize
p}(M;G)$, respectively, in \cite{T}. Since any index $r$ subgroup $H$
of $\Bbb Z^d$ is isomorphic to $\Bbb Z^d$ for any $r\ge1$, we get the
following theorem as a direct consequence of Theorem 6-3 for abelian
$\Gamma$. Similarly for the case $\Gamma=\Bbb Z^d_p$. In \cite{T}, the
same theorem was proved by induction.

\proclaim{Theorem 6-5} For any $d\ge0$, we have 
$$
\aligned
\sum_{n\ge0}q^n\chi^{\sssize(d)}(M^n;G\!\wr\!\frak S_n)
&=\Bigl[\prod_{r\ge1}(1-q^r)^{-j_r(\Bbb Z^d)}\Bigr]
^{\chi^{\sssize(d)}(M;G)} \\
\sum_{n\ge0}q^n\chi^{\sssize(d)}_{\sssize p}(M^n;G\!\wr\!\frak S_n) 
&=\Bigl[\prod_{r\ge0}(1-q^{p^r})^{-j_{p^r}(\Bbb Z^d_p)}\Bigr]
^{\chi^{\sssize(d)}_p(M;G)}.
\endaligned
\tag6-8
$$
\endproclaim
\demo{Proof} All we need to note is that since $\Gamma=\Bbb Z^d$ is
abelian, we can apply formula (6-5), and that the number $u_r(\Bbb
Z^d)$ of conjugacy classes of index $r$ subgroups of $\Bbb Z^d$ is the
same as the number $j_r(\Bbb Z^d)$ of index $r$ subgroups. The
corresponding statement for $\Bbb Z^d_p$ is also valid.
\qed
\enddemo

\specialhead
(ii)\qua Higher genus orbifold Euler characteristic of symmetric orbifolds
\endspecialhead

A direct application of Theorem 6-3 to the case $\Gamma=\Gamma_{\!g+1}$
for $g\ge0$ gives the following.

\proclaim{Theorem 6-6} Let $g\ge0$. Then, with $[H]_r$ as before, 
$$
\botsmash{
\sum_{n\ge0}q^n\chi_{\sssize\Gamma_{\!g+1}}(M^n;G\!\wr\!\frak S_n)
=\prod_{r\ge1}(1-q^r)^{\!\!-\!\!\!\!\sum\limits_{\ \ [H]_r}
\!\!\!\!
\chi_{\sssize[\Gamma_{\!g+1}\!/\!H]}^{}(M;G)}.
}
\tag6-9
$$
\endproclaim

Although all index $r$ subgroups of $\Gamma_{\!g+1}$ are isomorphic to
each other ($\cong\Gamma_{rg+1}$), they may not be conjugate to each
other. This is why we cannot be any more concrete than the above
expression. In the formula (6-6) with $\Gamma=\Gamma_{\!g+1}$,
normalizers of index $r$ subgroups appear. The size of the normalizer
$N_{\Gamma_{\!g+1}}(H)$ can be different for different $[H]$ of the
same index. Geometrically, the group $N_{\Gamma}^{}(H)/H$ is the group
of deck transformations of the covering space $\Sigma_{rg+1} @>>>
\Sigma_{g+1}$ corresponding to the conjugacy class of the index $r$
subgroup $H$. Thus $\chi_{\sssize[\Gamma_{\!g+1}\!/\!H]}(M;G)$ really
depends on the conjugacy class $[H]$.

The summation in the exponent of the right hand side of (6-9) can be
thought of as a summation over all isomorphism classes of index $r$
covering spaces over $\Sigma_{g+1}$. See \S7 for details.

A similar formula is valid when $\Gamma$ is the fundamental group
$\Lambda_{h+2}$ of non-orientable genus $h+2$ closed surface, or a
free group $F_{s+1}$ on $s+1$ generators. We omit their explicit
expressions here. Again, any further analysis requires information on
normalizers of subgroups.

\specialhead
(iii)\qua Higher order higher genus orbifold Euler characteristic of
symmetric orbifolds
\endspecialhead

We combine the previous two cases and consider the generating function
for the invariant $\chi^{\sssize(d)}_{\sssize\Gamma}(M;G)
=\chi_{\sssize\Bbb Z^d\times\Gamma}(M;G)$. This equality is due to
(2-2). We also consider the case for 
$\chi^{\sssize(d)}_{\sssize p,\Gamma}$. 

\proclaim{Theorem 6-7} For any $d\ge0$ and any $\Gamma$, with
$[H]_{\ell}$ as before, we have 
$$
\aligned
\sum_{n\ge0}q^n\chi^{\sssize(d)}_{\sssize\Gamma}(M^n;G\!\wr\!\frak S_n)
&=\prod_{r,\ell\ge1}\Bigl[(1-q^{r\ell})^{-j_r(\Bbb Z^d)}\Bigr]
^{\!\!\!\sum\limits_{\ [H]_{\ell}}\!\!\!
\!\!\!\!
\chi_{\sssize[\Gamma\!/\!H]}^{\sssize(d)}(M;G)}\\
\sum_{n\ge0}q^n\chi^{\sssize(d)}_{\sssize p,\Gamma}(M^n;G\!\wr\!\frak S_n)
&=\topsmash{
\prod\Sb r\ge0 \\ \ell\ge1 \endSb
\Bigl[(1-q^{p^r\ell})^{-j_{p^r}(\Bbb Z^d_p)}\Bigr]
^{\!\!\!\sum\limits_{\ \ [H]_{\ell}}\!\!\!\!
\chi_{\sssize p,[\Gamma\!/\!H]}^{\sssize(d)}(M;G)}.
}
\endaligned
\tag6-10
$$
\endproclaim
\demo{Proof} By Proposition 5-4 and Theorem 6-5, we have 
$$
\align
\sum_{n\ge0}q^n\chi_{\sssize\Gamma}^{\sssize(d)}(M^n;G\!\wr\!\frak S_n)
&=\prod_{[H]}\prod_{[\rho]}\Big[\sum_{r\ge0}q^{|\Gamma\!/\!H|r}
\chi^{\sssize(d)}\bigl((M^{\langle\rho\rangle})^r;
\text{\rm Aut}_{\Gamma\text{-}G}(P_{\rho})\!\wr\!\frak S_r\bigr)\Bigr] \\
&=\prod_{[H]}\prod_{[\rho]}
\Bigl[\prod_{r\ge1}(1-q^{r|\Gamma\!/\!H|})^{-j_r(\Bbb Z^d)}\Bigr]
^{\chi^{\sssize(d)}(M^{\langle\rho\rangle};
\text{Aut}_{\Gamma\text{-}G}(P_{\rho}))} \\
&=\prod_{[H]}
\Bigl[\prod_{r\ge1}(1-q^{r|\Gamma\!/\!H|})^{-j_r(\Bbb Z^d)}\Bigr]
^{\chi^{\sssize(d)}_{\sssize[\Gamma\!/\!H]}(M;G)} \\
&=\prod_{\ell\ge1}
\topsmash{
\Bigl[\prod_{r\ge1}(1-q^{r\ell})^{-j_r(\Bbb Z^d)}\Bigr]
^{\!\sum\limits_{\ [H]_{\ell}}\!
\!\!\chi_{\sssize[\Gamma\!/\!H]}^{\sssize(d)}(M;G)},
}
\endalign
$$
here, for the third equality, we summed over
$[\rho]\in\text{Hom}(H,G)/(N_{\Gamma}(H)\times G)$, and for the fourth
equality, we first summed over conjugacy classes of index $\ell$
subgroups of $\Gamma$. The case for $\chi^{\sssize(d)}_{\sssize
p,\Gamma}$ is similar.  This completes the proof.
\qed
\enddemo

The generating functions of two types of orbifold Euler
characteristics of symmetric orbifolds given in Theorem 5-5 and
Theorem 6-3 seem very different.  But by (1-5), these two orbifold
Euler characteristics are related by a simple formula
$\chi_{\sssize\Gamma\times\Bbb Z}^{\text{orb}}(M;G)=
\chi_{\sssize\Gamma}(M;G)$. Thus we expect to get
a nontrivial identity by comparing these two formulae of generating
functions. We compare these two formulae for $M=\text{pt}$. The
corresponding formulae are (5-5) and (6-6).

Now, applying (2-6) to the wreath product $G_n$, we obtain
$$
\frac{|\text{Hom}(\Gamma\times\Bbb Z, G_n)|}{|G_n|}
=|\text{Hom}(\Gamma,G_n)/G_n|.
$$
This is the link between (5-5) and (6-6). By directly comparing the
right hand sides of corresponding formulae, we obtain the identity:
$$
\exp\biggl[\sum_{r\ge1}\frac{q^r}r
\Bigl\{\!\!\!\!\!\!\!\!\!\!\!
\sum\Sb K \\ \ \ \ [\Gamma\times\Bbb Z:K]=r \endSb 
\!\!\!\!\!\!\!\!\!
\frac{|\text{Hom}(K,G)|}{|G|}\Bigr\}\biggr]
=\topsmash{
\prod_{r\ge1}(1-q^r)^{\!\!-\!\!\sum\limits_{\ [H]_{r}}
\!\!\!
|\text{\rm Hom}(H,G)/(N_{\Gamma}(H)\times G)|}.
}
\tag6-11
$$
Expanding this identity, we see that the following proposition must
hold. Although the identity (6-11) is its proof, we also give a group
theoretic proof to show what kind of group theoretic ingredients are 
involved. 

\proclaim{Proposition 6-8} Let $G$ be a finite group and let
$\Gamma$ be an arbitrary group. For any $r\ge1$, we have
$$
\sum_{[\Gamma\times\Bbb Z;K]=r}\!\!\!\!\!\!\!\!
\frac{|\text{\rm Hom}(K,G)|}{|G|}
=\!\sum_{\ell|r}\ell\Bigl\{
\!\!\!\!\!\!\!
\botsmash{
\sum\Sb [H] \\ \ \ |\Gamma\!/\!H|=\ell \endSb
}
\!\!\!\!\!\!
|\text{\rm Hom}(H,G)/(N_{\Gamma}(H)\times G)|\Bigr\}.
\tag6-12
$$
When $G$ is trivial, this formula reduces to a well-known formula 
$j_r(\Gamma\times\Bbb Z)
=\botsmash{\sum\limits_{\ell|r}}\ell\cdot u_{\ell}(\Gamma)$. 
\endproclaim
\demo{Proof} Let $\ell|r$. For an index $\ell$ subgroup $H\subset
\Gamma$ and $\overline{z}\in N_{\Gamma}^{}(H)/H$, let
$$
K_{(H,\overline{z})}=\{(w,ra/\ell)\in N_{\Gamma}^{}(H)\times\Bbb Z \mid
\overline{w}=\overline{z}^a \text{ in }
N_{\Gamma}^{}(H)/H\}\subset\Gamma\times\Bbb Z.
$$
Then $K_{(H,\overline{z})}$ is a subgroup of index $r$ in
$\Gamma\times\Bbb Z$, and any index $r$ subgroup of $\Gamma\times\Bbb
Z$ is of this form [\Cite{Su},\,p.140]. Thus, index $r$ subgroups of
$\Gamma\times\Bbb Z$ are parametrized by the pair $(H,\overline{z})$,
where the index of $H$ divides $r$ and $\overline{z}$ is an arbitrary
element in $N_{\Gamma}^{}(H)/H$. Let $z\in N_{\Gamma}^{}(H)$ be any
element whose reduction mod $H$ is $\overline{z}$. Then any element in
$K_{(H,\overline{z})}$ is of the form $(z,r/\ell)^a\cdot(h,0)
\in N_{\Gamma}^{}(H)\times\Bbb Z$ for some unique $a\in\Bbb Z$ and
$h\in H$. Now let $\phi:K_{(H,\overline{z})} @>>> G$ be a
homomorphism. Let $\rho=\phi|_H:H @>>> G$ be its restriction. Let
$\phi\bigl((z,r/\ell)\bigr)=g\in G$. Since we have
$(z,r/\ell)(h,0)(z,r/\ell)^{-1}=(zhz^{-1},0)$ in
$K_{(H,\overline{z})}$, we have $g\rho(h)g^{-1}=\rho(zhz^{-1})$ in
$G$. This means that $\rho^z$, defined by $\rho^z(h)=\rho(zhz^{-1})$
for $h\in H$, and $\rho$ are $G$-conjugate, which we denote by
$\rho^z\sim\rho$. Conversely, given any $\rho:H @>>> G$ and any $g\in
G$ satisfying a relation $g\rho(h)g^{-1}=\rho^z(h)$ for any $h\in H$,
we can define a homomorphism $\phi:K_{(H,\overline{z})} @>>> G$ by
$\phi\bigl((z^ah,ra/\ell)\bigr)=g^a\rho(h)$. Since for a given $\rho:H
@>>> G$ such that $\rho^z\sim\rho$, there are $|C(\rho)|$ choices of
$g\in G$ satisfying $g\rho g^{-1}=\rho^z$, we see that
$|\text{Hom}(K_{(H,\overline{z})},G)|=\sum_{\rho}|C(\rho)|$, where the
summation runs over all $\rho: H @>>> G$ such that
$\rho^z\sim\rho$. Hence,
$$
\multline
\sum_{\overline{z}\in N_{\Gamma}^{}(H)/H}
\!\!\!\!\!\!\!\!\!\frac{|\text{Hom}(K_{(H\overline{z})}^{},G)|}{|G|}
=\!\!\!\!\!\!\!\!\!\!
\sum_{\overline{z}\in N_{\Gamma}^{}(H)/H\ \ \ \ }
\!\!\!\!\!\!\!\!
\sum\Sb \rho:H @>>> G \\\rho^z\sim \rho \endSb
\!\!\!\!\!\frac{|C(\rho)|}{|G|}
=\!\!\!\!\!\!\!\!
\sum_{\rho:H @>>> G\ \ }
\!\!\!\!\!\!\!\!\frac{|C(\rho)|}{|G|}\!\!\!\!\!\!\!\!\!\!
\sum\Sb \ \ \ \overline{z}\in N_{\Gamma}^{}(H)/H 
\\ \rho^z\sim\rho \endSb 
\!\!\!\!\!\!\!\!\!\!\!\!1 \\
=\!\!\!\!\!\!\!\sum_{\rho:H @>>> G\ \ }
\!\!\!\!\!\!\!\frac{|C(\rho)|}{|G|} |N_{\Gamma}^{\rho}(H)/H|
=\!\!\!\!\!\!\!\!\!\sum_{\rho:H @>>> G\ \ }
\!\!\!\!\frac1{\#(\rho)} |N_{\Gamma}^{\rho}(H)/H|
=\!\!\!\!\!\!\!\!\!\!\!\!\!\!\!\!\!\!\!\!\!\!
\sum_{\ \ \ \ \ \ (\rho)\in\text{Hom}(H,G)/G}
\!\!\!\!\!\!\!\!\!\!\!\!\!\!\!\!\!\!\!\!\!
|N_{\Gamma}^{\rho}(H)/H|.
\endmultline
$$
Here $(\rho)$ denotes the $G$-conjugacy class of $\rho$. The third
equality above is by definition of $N^{\rho}_{\sssize\Gamma}(H)$ in
(4-1). On the other hand, $|\text{Hom}(H,G)/(N_{\Gamma}(H)\times G)|$
is equal to the number of orbits under the $N_{\Gamma}(H)$-action on
$\text{Hom}(H,G)/G$. At $(\rho)\in \text{Hom}(H,G)/G$, the
isotropy subgroup of $N_{\Gamma}^{}(H)$-action is
$N_{\Gamma}^{\rho}(H)$. Thus,
$$
\multline
|\text{Hom}(H,G)/(N_{\Gamma}(H)\times G)|
=\!\!\!\!\!\!\!\!\!\!\!\!\!\!\!\!\!\!
\sum_{\ \ \ (\rho)\in\text{Hom}(H,G)/G}
\!\!\!\!\!\!\!\!\!\!\!\!\!\!\!\!\!
1/|(\rho)\text{-orbit}|
=\!\!\!\!\!\!\!\!\!\!\!\!\!\!\!\!\!\!\!\!\!\!\!\!\!\!
\sum_{\ \ \ \ \ \ \ \ \ \ (\rho)\in\text{Hom}(H,G)/G}
\!\!\!\!\!\!\!\!\!\!\!\!\!\!\!\!\!\!\!\!\!\!\!\!\!
|N_{\Gamma}^{\rho}(H)/N_{\Gamma}^{}(H)| \\
=(1/|N_{\Gamma}^{}(H)/H|)
\!\!\!\!\!\!\!\!\!\!\!\!\!\!\!\!\!\!\!\!
\sum_{\ \ \ \ \ \ (\rho)\in\text{Hom}(H,G)/G}
\!\!\!\!\!\!\!\!\!\!\!\!\!\!\!\!\!\!\!\!
|N_{\Gamma}^{\rho}(H)/H|
=\frac{\#[H]}{\ell}
\!\!\!\!\!\!\!\!\!\!\!
\sum_{\ \ \overline{z}\in N_{\Gamma}^{}(H)/H}
\!\!\!\!\!\!\!\!\!\!\!\!
\frac{|\text{Hom}(K_{(H,\overline{z})},G)|}{|G|},
\endmultline
$$
where $|(\rho)\text{-orbit}|$ denotes the length of the
$N_{\Gamma}^{}(H)$-orbit in $\text{Hom}(H,G)/G$\break through $(\rho)$. The
last equality is by the above calculation. Since index $r$ subgroups
of $\Gamma\times\Bbb Z$ are parametrized by pairs $(H,\overline{z})$,
where $H$ is a subgroup of $\Gamma$ whose index divides $r$ and
$\overline{z}\in N_{\Gamma}^{}(H)/H$, we have
$$
\multline
\sum_{[\Gamma\times\Bbb Z:K]=r}
\!\!\!\!\!\!\!\!\frac{|\text{Hom}(K,G)|}{|G|}
=\sum_{\ell|r}\!\sum_{H_{\ell}}
\!\!\!\!\!\!\!\!\!\!\!\!\!\!\!
\sum_{\ \ \ \ \ \ \ \ \overline{z}\in N_{\Gamma}^{}(H)/H}
\!\!\!\!\!\!\!\!\!\!\!\!\!\!\!\!
\frac{|\text{Hom}(K_{(H,\overline{z})},G)|}{|G|} \\
=\sum_{\ell|r}\sum_{H_{\ell}}
\frac{\ell}{\#[H]}
|\text{Hom}(H,G)/(N_{\Gamma}^{}(H)\times G)|\\
=\sum_{\ell|r}\!\sum_{\ [H]_{\ell}}\!
\ell\, |\text{Hom}(H,G)/(N_{\Gamma}^{}(H)\times G)|.
\endmultline
$$
This completes the proof. 
\qed
\enddemo
From the proof, it is clear that $K_{(H,\overline{z})}$ is a
semi-direct product $H\rtimes\Bbb Z$, where $\Bbb Z\subset
K_{(H,\overline{z})}$ is a subgroup generated by the element
$(z,r/\ell)\in K_{(H,\overline{z})}$.

We next describe orbifold invariants associated to general
$\Gamma$-sets. Let $X$ be any $\Gamma$-set of order $n$. Then for any
given multiplicative invariant $\varphi(M;G)$, we define the orbifold
invariant associated to the $\Gamma$-isomorphism class $[X]$ by
$$
\varphi_{\sssize[X]}(M;G)\overset{\text{def}}\to=
\!\!\!\!\!\!
\sum_{\ \ [P @>>> X]}\!\!\!\!\!
\varphi\bigl(\Cal S[P\times_GM]^{\Gamma};
\text{Aut}_{\Gamma\text{-}G}(P)\bigr)
\cong\!\!\!\!\!\!\!\!\!\!\!\sum_{\ \ \ [\theta]\in\pi^{-1}([X])}
\!\!\!\!\!\!\!\!\!\!\!
\varphi\bigl((M^n)^{\langle\theta\rangle};C_{G_n}(\theta)\bigr),
\tag6-13
$$
where the first summation is over all isomorphism classes of
$\Gamma$-equivariant $G$-principal bundles $P$ over $X$, and $\pi:
\text{Hom}(\Gamma,G_n)/G_n @>>>
\text{Hom}(\Gamma,\frak S_n)/\frak S_n$ is as before. When $X$ is an
empty set $\emptyset$, we let $\varphi_{\sssize\emptyset}(M;G)=1$. The
second isomorphism is due to Theorem 3-1, Lemma 4-1 and Proposition
5-1. The details of the proof are similar to the proof of Proposition
6-1. We could define $\varphi_{\sssize[\Gamma\!/\!H]}(M;G)$ by the
above conceptual formula instead of (6-1), but (6-1) is more
explicit. The geometric definition (6-13) in terms of
$\Gamma$-equivariant $G$-principal bundles has an advantage in its
simplicity. 

In view of Theorem 3-1 and (6-13), generating functions for
$\varphi_{\sssize\Gamma}$ and $\varphi_{\sssize[X]}$ are related by 
$$
\sum_{n\ge0}q^n\varphi_{\sssize\Gamma}(M^n;G\!\wr\!\frak S_n)
=\sum_{[X]}q^{\sssize|X|}\varphi_{\sssize[X]}(M;G),
\tag6-14
$$
where on the right hand side, $[X]$ runs over the set of all the
isomorphism classes of finite $\Gamma$-sets. 

Next we discuss multiplicativity of $\varphi_{\sssize[X]}(M;G)$ with
respect to $[X]$. It is not true that we have 
$\varphi_{\sssize[X]}(M;G)=\varphi_{\sssize[X_1]}(M;G)\cdot
\varphi_{\sssize[X_2]}(M;G)$ whenever $X=X_1\coprod X_2$. However, when
$X_1$ and $X_2$ are ``prime'' to each other, namely when $X_1$ and
$X_2$ do not contain the same $\Gamma$-transitive sets in common in
their decomposition into $\Gamma$-orbits, the above multiplicativity
is valid.

\proclaim{Proposition 6-9} Let $\varphi(M;G)$ be a multiplicative
orbifold invariant for\break $(M;G)$. Let $\Gamma$ be a group. For any
finite $\Gamma$-set $X$, let $X=\botsmash{\coprod_{[H]}}
r(H)\cdot\Gamma\!/\!H$ be its decomposition into $\Gamma$-orbits, where
$r(H)$ is the number of $\Gamma$-orbits which are isomorphic to
$\Gamma\!/\!H$. Then
$$
\varphi_{\sssize[X]}(M;G)=\prod_{[H]}
\varphi_{\sssize r(H)[\Gamma\!/\!H]}(M;G).
\tag6-15
$$
Consequently, the generating function \rom{(6-14)} decomposes into a
product
$$
\sum_{n\ge0}q^n\varphi_{\sssize\Gamma}(M^n;G\!\wr\!\frak S_n)
=\prod_{[H]}\Bigl\{\sum_{r\ge0}q^{|\Gamma\!/\!H|r}
\varphi_{\sssize r[\Gamma\!/\!H]}(M;G)\Bigr\}.
\tag6-16
$$
Furthermore, the generating function of 
$\varphi_{\sssize r[\Gamma\!/\!H]}(M;G)$ for $r\ge0$ is given by 
$$
\sum_{r\ge0}q^r\varphi_{\sssize r[\Gamma\!/\!H]}(M;G)
=\prod_{[\rho]}\Bigl\{
\sum_{r\ge0}q^r\varphi\bigl((M^{\langle\rho\rangle})^r;
\text{\rm Aut}_{\Gamma\text{-}G}(P_{\rho})
\!\wr\!\frak S_r\bigr)\Bigr\},
\tag6-17
$$
where $[\rho]$ runs over the set $\text{\rm
Hom}(H,G)/(N_{\Gamma}(H)\times G)$. 
\endproclaim
\demo{Proof} For any $\Gamma$-equivariant $G$-bundle $P$ over
$X$, let $P_H$ be the part of $P$ above $r(H)[\Gamma\!/\!H]$. Then, 
using the multiplicativity of $\varphi(\ \cdot\ ;\ \cdot\ )$, we have 
$$
\align
\varphi_{\sssize[X]}(M;G)
&=\!\!\!\!\!\!\!\sum_{\ \ [P @>>> X]}\!\!\!\!\!\!\!
\varphi\bigl(\Cal S[P\times_GM]^{\Gamma};
\text{Aut}_{\Gamma\text{-}G}(P)\bigr)\\
&=\!\!\!\!\!\!\!\sum_{[P @>>> X]\ \ }\!\!\!\!\!\prod_{[H]}\varphi
\bigl(\Cal S[P_H\times_GM]^{\Gamma};
\text{Aut}_{\Gamma\text{-}G}(P_H)\bigr) \\
&=\prod_{[H]}\!\!\!\!\!\!\!\!\!\!\!\!\!\!\!\!\!\!\!\!\!\!\!\!
\sum_{\ \ \ \ \ \ \ \ \ \ \ \ [P_H @>>> r(H)\cdot\Gamma\!/\!H]}
\!\!\!\!\!\!\!\!\!\!\!\!\!\!\!\!\!\!\!\!\!\!\!\!\!\!\!
\varphi\bigl(\Cal S[P_H\times_GM]^{\Gamma};
\text{Aut}_{\Gamma\text{-}G}(P_H)\bigr)
=\prod_{[H]}\varphi_{\sssize r(H)[\Gamma\!/\!H]}(M;G),
\endalign
$$
where in the fourth expression, $[P_H]$ rums over the set of all
isomorphism classes of $\Gamma$-$G$ bundles over the $\Gamma$-set
$r(H)[\Gamma\!/\!H]$. The formula (6-16) follows from this.

The proof of (6-17) is analogous to the argument in the proof of
Proposition 5-4, and goes as follows. For any $\Gamma$-equivariant
$G$-bundle over $\amalg^r\Gamma\!/\!H$, let $r_{\rho}$ be the number of
$\Gamma$-irreducible $G$-bundle $P_{\rho} @>>> \Gamma\!/\!H$ appearing in
the irreducible decomposition of $P$, where
$[\rho]\in\text{Hom}(H,G)/(N_{\Gamma}^{}(H)\times G)$. Then,
decomposing $\Gamma$-equivariant $G$-bundles $P @>>> r\cdot\Gamma\!/\!H$
according to $[\rho]$, we have
$$
\align
\sum_{r\ge0}q^r\varphi_{\sssize r[\Gamma\!/\!H]}(M;G)
&=\sum_{r\ge0}q^r\!\!\!\!\!
\sum_{\sum\limits_{[\rho]}r_{\!\rho}=r\ \ }
\!\!\!\!
\prod_{[\rho]}\varphi\bigl((M^{\langle\rho\rangle})^{r_{\!\rho}};
\text{Aut}_{\Gamma\text{-}G}(P_{\rho})\!\wr\!\frak S_{r_{\!\rho}}\bigr) \\
&=\!\sum_{r_{\!\rho}\ge0}\prod_{[\rho]}q^{r_{\!\rho}}
\varphi\bigl((M^{\langle\rho\rangle})^{r_{\!\rho}};
\text{Aut}_{\Gamma\text{-}G}(P_{\rho})\!\wr\!\frak S_{r_{\!\rho}}\bigr)\\
&=\prod_{[\rho]}\Bigl\{
\sum_{r\ge0}q^r\varphi\bigl((M^{\langle\rho\rangle})^r;
\text{\rm Aut}_{\Gamma\text{-}G}(P_{\rho})
\!\wr\!\frak S_r\bigr)\Bigr\}.
\endalign
$$
This completes the proof.  
\qed
\enddemo
 
Note that formulae (6-14), (6-16), and (6-17) essentially give a
proof of Proposition 5-4, from the point of view of $\Gamma$-sets. 

Next we examine the case in which $G$ is the trivial group. We have
already discussed this case for a transitive $\Gamma$-set $X$. The
proof of the following lemma is straightforward.

\proclaim{Lemma 6-10} With the same notations as above, 
suppose $G=\{e\}$, the trivial group. Then 
$$
\varphi_{\sssize[X]}(M;\{e\})=\prod_{[H]}\varphi\bigl(M^{r(H)};
\bigl(H\backslash N_{\Gamma}^{}(H)\bigr)\!\wr\!\frak S_{r(H)}\bigr),
\tag6-18
$$
where $H\backslash N_{\Gamma}^{}(H)$ acts trivially on $M$. For the
case of the orbifold Euler characteristic, we have
$$
\chi_{\sssize[X]}(M;\{e\})=\prod_{[H]}\chi\bigl(SP^{r(H)}(M)\bigr).
\tag6-19
$$
In particular, if $X$ is a transitive $\Gamma$-set, then
$\chi_{\sssize[X]}(M;\{e\})=\chi(M)$. 
\endproclaim

\head
Orbifold invariants associated to covering spaces
\endhead

Let $\Sigma$ be a connected manifold and let $M$ be a $G$-manifold
where $G$ is a finite group. In \S2, we considered a twisted space
$\Bbb L_{\Sigma}^{}(M;G)$ defined by
$$
\Bbb L_{\Sigma}(M;G)=\coprod_{[P]}
\bigl[\text{Map}_G(P,M)/\text{Aut}_G^{}(P)\bigr],
\tag7-1
$$
where $[P]$ runs over the set of isomorphism classes of $G$-principal
bundles over $\Sigma$, and $\text{Aut}_G^{}(P)$ is the group of
$G$-bundle automorphisms inducing the identity map on the base
$\Sigma$. Let $\Gamma=\pi_1(\Sigma)$. We defined the
orbifold invariant associated to a group $\Gamma$ in (2-1) as the
orbifold invariant of the subset of $\Bbb L_{\Sigma}(M;G)$ consisting
of locally constant $G$-equivariant maps described in (2-14):
$$
\varphi_{\sssize\Gamma}(M;G)
=\varphi_{\text{orb}}^{}\bigl(\Bbb L_{\Sigma}(M;G)_{\text{l.c.}}\bigr)
=\!\!\!\!\!\!\!\!\!\!\!\!\!\!\!\!\!\!\!\!\!\!
\sum_{\ \ \ \ \ \ \ [\phi]\in\text{Hom}(\Gamma,G)/G}
\!\!\!\!\!\!\!\!\!\!\!\!\!\!\!\!\!\!\!\!\!
\varphi\bigl(M^{\langle\phi\rangle};C(\phi)\bigr),
\tag7-2
$$
where l.c. stands for ``locally constant.''

We generalize the above construction of the twisted space $\Bbb
L_{\Sigma}(M;G)$ associated to $\Sigma$ to twisted spaces associated
to finite covering spaces of $\Sigma$. Let
$\pi_{\sssize\Sigma'\!\!/\Sigma}^{}:\Sigma' @>>> \Sigma$ be a (not
necessarily connected) finite covering space of $\Sigma$. The idea of
constructing twisted spaces in this new context is the same as before:
we consider liftings of maps $\text{Map}(\Sigma',M/G)$ to maps from
$G$-principal bundles $P$ over $\Sigma'$ to $M$. So let $\pi:P @>>>
\Sigma'$ be a $G$-principal bundle over $\Sigma'$ and consider the set
$\text{Map}_G(P,M)$ of $G$-equivariant maps from $P$ to $M$. Now the
equivalence relation we introduce among these maps is different from
the one in \S2. Given two $G$-maps $\gamma_1:P_1 @>>> M$ and
$\gamma_2:P_2 @>>> M$, $\gamma_1$ and $\gamma_2$ are said to be
equivalent in the present context if there exists a $G$-bundle
isomorphism $\alpha:P_1 @>>> P_2$ inducing a {\it deck transformation}
of $\Sigma'$ over $\Sigma$ such that
$\gamma_2\circ\alpha=\gamma_1$. We call such $\alpha$ a $G$-bundle
isomorphism {\it over a covering space} $\Sigma'/\Sigma$.

Let $\text{Aut}_G^{}(P)_{\Sigma'\!\!/\Sigma}$ be the group of
$G$-automorphisms of $P$ over the covering space
$\Sigma'\!\!/\Sigma$. The set of equivalence classes denoted by $\Bbb
L_{\Sigma'\!\!/\Sigma}^{}(M;G)$ is our new twisted space associated to
the covering space $\Sigma'\!\!/\Sigma$. Thus,
$$
\Bbb L_{\Sigma'\!\!/\Sigma}^{}(M;G)
=\!\!\!\!\!\!\!\!\!\coprod_{[P @>>> \Sigma'\!\!/\Sigma]}\!\!\!\!\!\!\!\!
\bigl[\text{Map}_G^{}(P,M)/\text{Aut}_G^{}(P)_{\Sigma'\!\!/\Sigma}\bigr].
\tag7-3
$$
Here the disjoint union runs over the set of all $G$-bundle
isomorphism classes over the covering space $\Sigma'\!\!/\Sigma$ in
the above sense. We consider the subset of locally constant
$G$-equivariant maps:
$$
\Bbb L_{\Sigma'\!\!/\Sigma}^{}(M;G)_{\text{l.c.}}^{}\!\!
=\!\!\!\!\!\!\!\!\!
\coprod_{[P @>>> \Sigma'\!\!/\Sigma]}
\!\!\!\!\!\!\!\!
\bigl[\text{Map}_G^{}(P,M)_{\text{l.c.}}^{}
/\text{Aut}_G^{}(P)_{\Sigma'\!\!/\Sigma}\bigr].
\tag7-4
$$
Here note that the group $\text{Aut}_G^{}(P)_{\Sigma'\!\!/\Sigma}$
still acts on the space $\text{Map}_G^{}(P,M)_{\text{l.c.}}^{}$ of
locally constant $G$-equivariant maps. By definition, the orbifold
invariant associated to a finite covering $\Sigma' @>>> \Sigma$ is the
$\varphi$-invariant of $\Bbb
L_{\Sigma'\!\!/\Sigma}^{}(M;G)_{\text{l.c.}}^{}$, viewed
equivariantly rather than as a quotient:
$$
\varphi_{\sssize[\Sigma'\!\!/\Sigma]}^{}(M;G)=
\varphi_{\text{orb}}
\bigl(\Bbb L_{\Sigma'\!\!/\Sigma}^{}(M;G)_{\text{l.c.}}^{}\bigr)
\!\!
\overset{\text{def}}\to=
\!\!\!\!\!\!\!\!
\sum_{[P @>>> \Sigma'\!\!/\Sigma]}
\!\!\!\!\!\!\!\!
\varphi\bigl(\text{Map}_G^{}(P,M)_{\text{l.c.}}^{};
\text{Aut}_G^{}(P)_{\Sigma'\!\!/\Sigma}\bigr). 
\tag7-5
$$
We analyze quantities appearing in the above formula to prove the
following theorem, which is the main result of this section. 

\proclaim{Theorem 7-1} Let $\pi_{\sssize\Sigma'\!\!/\Sigma}^{}:
\Sigma' @>>> \Sigma$ be a finite \rom{(}not necessarily
connected\rom{)} covering space over a connected manifold $\Sigma$
with a base point $x_0$. Let $\Gamma=\pi_1(\Sigma,x_0)$. Let
$X=\Sigma'_{x_0}$ be the fibre over $x_0\in\Sigma$, and let $[X]$ be
its isomorphism class as a $\Gamma$-set. Then, the orbifold invariant
associated to the covering $\Sigma'\!\!/\Sigma$ defined in \rom{(7-5)}
is equal to the orbifold invariant associated to the isomorphism class
of the $\Gamma$-set $[X]$ defined in \rom{(6-13)}\rom{:}
$$
\varphi_{\sssize[\Sigma'\!\!/\Sigma]}^{}(M;G)=
\varphi_{\sssize[X]}^{}(M;G).
\tag7-6
$$

The invariant $\varphi_{\sssize[\Sigma'\!\!/\Sigma]}^{}(M;G)$ is
multiplicative in the following sense. Let $\pi_1:\Sigma_1' @>>>
\Sigma$ and $\pi_2: \Sigma_2' @>>> \Sigma$ be two coverings over
$\Sigma$. Suppose $\Sigma_1'$ and $\Sigma_2'$ do not contain any
isomorphic connected coverings in common. Then
$$
\varphi_{\sssize[(\Sigma_1'\cup\Sigma_2')/\Sigma]}^{}(M;G)
=\varphi_{\sssize[\Sigma_1'\!/\Sigma]}^{}(M;G)
\varphi_{\sssize[\Sigma'_2\!/\Sigma]}^{}(M;G).
\tag7-7
$$
\endproclaim

Since any finite $\Gamma$-set $X$ can be a fibre of some finite
covering space $\Sigma'=\widetilde{\Sigma}\times_{\Gamma}X @>>>
\Sigma$, where $\widetilde{\Sigma}$ is the universal covering space of
$\Sigma$, the identity (7-6) gives a geometric meaning of orbifold
invariants associated to $\Gamma$-sets in terms of covering
spaces. The formula (7-7) corresponds to the formula (6-15) for
orbifold invariants associated to $\Gamma$-sets.

To prove Theorem 7-1, we first recall some basic facts on (universal)
covering spaces. Let $\widetilde{\pi}:\widetilde{\Sigma} @>>> \Sigma$
be the universal covering space of $\Sigma$, and let $x_0\in\Sigma$ be
a base point. Let $\Gamma=\pi_1(\Sigma,x_0)$ be the fundamental group
of $\Sigma$. One convenient description of $\tilde{\Sigma}$ is the
following one:
$$
\widetilde{\Sigma}=\{(x,[\gamma]) \mid x\in\Sigma \text{ and $[\gamma]$ 
is the homotopy class of a path $\gamma$ from $x$ to $x_0$} \}.
\tag7-8
$$
We can introduce a suitable topology on $\widetilde{\Sigma}$ so that
the projection map $\widetilde{\Sigma} @>>> \Sigma$ is a covering
map. In this description, the point
$(x,[\gamma])\in\widetilde{\Sigma}$ is the end point of the lift in
$\widetilde{\Sigma}$ of $\gamma^{-1}$ starting at the natural base
point $y_0=(x_0,[c_{x_0}])$ of $\widetilde{\Sigma}$, where $c_{x_0}$
is the constant path at $x_0\in\Sigma$. The fundamental group $\Gamma$
acts from the {\it right} on $\widetilde{\Sigma}$ as deck
transformations by
$$
\aligned
\widetilde{\Sigma}\times\Gamma & @>{\hphantom{cons}}>> 
\widetilde{\Sigma} \\
\bigl((x,[\gamma]),[\eta&]\bigr)  \mapsto (x,[\gamma][\eta]),
\endaligned
\tag7-9
$$
where $[\eta]\in\Gamma$. This action is well defined because $\gamma$
is a path from $x$ to the base point $x_0$ and $\eta$ is a loop at
$x_0$. Here, we regard $\widetilde{\Sigma}$ as a right
$\Gamma$-principal bundle. This is the reason why we let $\Gamma$ act
on $\widetilde{\Sigma}$ on the right.  On the other hand, the group
$\Gamma$ also acts on the fibre
$\widetilde{F}_{x_0}=\widetilde{\pi}^{-1}(x_0)$ from the {\it left} by
$$
\aligned
\Gamma\times\widetilde{F}_{x_0} & @>{\hphantom{cons}}>> 
\widetilde{F}_{x_0} \\
\bigl([\eta], (x_0,[\gamma])\bigr)&  \mapsto (x_0,[\eta][\gamma]).
\endaligned
\tag7-10
$$
Again, this is well defined because the path $\gamma$ is now a loop
based at $x_0$. Note that for a given $y=(x_0,[\gamma])$, the result
of the action of $[\eta]$ on $y$, $[\eta]\cdot y$, is the end point of
the lift in $\widetilde{\Sigma}$ of $\eta^{-1}$ starting at
$y$. Obviously, the left action of $\Gamma$ on $\widetilde{F}_{x_0}$,
and the right action of $\Gamma$ as deck transformations
restricted to $\widetilde{F}_{x_0}$, commute. In some literature, the
action of $\Gamma$ on the fibre is defined from the right. We decided
to use the above left action in this paper to make it explicit that
the action of $\Gamma$ on the fibre commutes with its action on
$\widetilde{\Sigma}$ as deck transformations. 

The story is similar for any connected covering space
$\pi_{\sssize\Sigma'\!\!/\Sigma}: \Sigma' @>>> \Sigma$. Any such
$\Sigma'$ is isomorphic to a covering space of the form
$\widetilde{\Sigma}/H$ for some subgroup $H\subset\Gamma$, determined
up to conjugacy, and can be described as $\Sigma'=\{(x,[\gamma]H)\}$
where $x\in\Sigma$ and $[\gamma]$ is the homotopy class of paths from
$x$ to the base point $x_0$. The point $y_0=(x_0,H)\in\Sigma'$ can be
used as the natural base point of
$\Sigma'\cong\widetilde{\Sigma}/H$. The group of deck transformations
for the covering space $\Sigma'@>>> \Sigma$ given by
$G(\Sigma'/\Sigma)\cong N_{\Gamma}^{}(H)/H$ acts on $\Sigma'$ from the
right, and $\Gamma$ acts on the fibre
$F_{x_0}=\pi_{\sssize\Sigma'\!\!/\Sigma}^{-1}(x_0)$ over $x_0$ from
the left in the same way as for $\widetilde{\Sigma}$:
$$
\aligned
\Gamma\times {F}_{x_0} & @>{\hphantom{cons}}>> 
{F}_{x_0} \\
\bigl([\eta], (x_0,[\gamma]H)\bigr)& \mapsto (x_0,[\eta][\gamma]H),
\endaligned
\qquad
\aligned
\Sigma'\times (N_{\Gamma}^{}(H)/H) & @>{\hphantom{cons}}>> \Sigma' \\
\bigl((x,[\gamma]H),[\eta]H\bigr)&  \mapsto (x,[\gamma][\eta]H).
\endaligned
\tag7-11
$$
For any $y\in F_{x_0}$, the point $[\eta]\cdot y\in F_{x_0}$ is
obtained as the end point of the lift in $\Sigma'$ of $\eta^{-1}$
starting at $y\in F_{x_0}$. Similarly, given
$y=(x,[\gamma]H)\in\Sigma'$, the result of the action of the deck
transformation by $[\eta]H$ on $y$ is the point $y\cdot
[\eta]H$ obtained as the end point of the lift in $\Sigma'$ of
$\gamma^{-1}$ starting at $[\eta]\cdot y_0=(x_0,[\eta]H)$.

If $\Sigma'$ is not connected, the above facts apply to each connected
component. Note that the intersection of a connected component of
$\Sigma'$ with the fibre $F_{x_0}$ is a single $\Gamma$-orbit. 

Next we study basic properties of $G$-principal bundles over covering
spaces. Let $\pi': P @>>> \Sigma'$ be a $G$-principal bundle over a
not necessarily connected covering space $\Sigma'$ over $\Sigma$. Let
$\pi=\pi_{\sssize\Sigma'\!\!/\Sigma}\circ\pi': P @>>> \Sigma$. This is
a covering space over $\Sigma$ on which $G$ acts as a group of deck
transformations. Let $P_{x_0}=\pi^{-1}(x_0)$ and
$\Sigma'_{x_0}=\pi_{\sssize\Sigma'\!\!/\Sigma}^{-1}(x_0)$ be fibres
over $x_0$. Then $P_{x_0}$ has the structure of a $\Gamma$-equivariant
$G$-principal bundle over $\Sigma'_{x_0}$. It is easy to see that the
bundle $P$ is completely determined by its fibre $P_{x_0}$, in the
following way. 

\proclaim{Proposition 7-2} Let $\pi': P @>>> \Sigma'$ be a
$G$-principal bundle over a covering space
$\pi_{\sssize\Sigma'\!\!/\Sigma}^{}:
\Sigma' @>>> \Sigma$, where $\Sigma$ is a connected manifold with a
base point $x_0$. The the fibre $P_{x_0}$ has the structure of a
$\Gamma$-equivariant $G$-principal bundle over $\Sigma'_{x_0}$, and
$\pi': P @>>> \Sigma$ is determined by the $\Gamma$-equivariant
$G$-bundle structure of the fibre $P_{x_0} @>>> \Sigma'_{x_0}$. That
is, we have the following commutative diagram\rom{:}
$$
\CD
\widetilde{\Sigma}\times_{\Gamma}P_{x_0} @>{\cong}>{G}> P  \\
@VVV @VVV \\
\widetilde{\Sigma}\times_{\Gamma}\Sigma'_{x_0}
@>{\cong}>{N_{\Gamma}^{}(H)/H}> \Sigma' \\
@VVV @VVV \\
\widetilde{\Sigma}\times_{\Gamma}\{x_0\} @>{\cong}>> \Sigma
\endCD
\tag7-12
$$
Here the top horizontal map is a $G$-equivariant map sending
$\bigl[(x,[\gamma]),p\bigr]$ to $p'\in P$ over $x\in\Sigma$ obtained
as the end point of the lift in $P$ of $\gamma^{-1}$ starting at $p\in
P$. The middle horizontal map is $N_{\Gamma}^{}(H)/H$-equivariant map
sending $\bigl[(x,[\gamma]), (x_0,[\eta]H)\bigr]$ to
$\bigl(x,[\gamma][\eta]H\bigr)$.
\endproclaim
\demo{Proof} The proof is routine. So we only give a brief proof for
the top horizontal map. First we check that this correspondence is
well defined. For any $[\eta]\in\Gamma$, we consider two pairs
$\bigl((x,[\gamma][\eta]),p\bigr)$ and
$\bigl((x,[\gamma]),[\eta]p\bigr)$ in $\widetilde{\Sigma}\times
P_{x_0}$, and compare the corresponding points in $P$ given by the
procedure described above. The point $p'$ corresponding to the first
pair is the end point of the lift in $P$ of
$(\gamma\eta)^{-1}=\eta^{-1}\gamma^{-1}$ starting at $p$. Since
$[\eta]\cdot p$ is the end point of the lift in $P$ of $\eta^{-1}$
starting at $p$, the point $p'$ is the same as the end point $p''$ of
the lift in $P$ of $\gamma^{-1}$ starting at $[\eta]\cdot p$. But
$p''$ is also the point corresponding to the second pair above. This
proves that the top horizontal correspondence is well defined.

For $G$-equivariance, the point corresponding to $[(x,[\gamma]),p]$ is
the end point $p_1$ of the lift $\widetilde{\gamma}^{-1}_1$ in $P$ of
$\gamma^{-1}$ starting at $p$. Now the point corresponding to
$[(x,[\gamma]),pg]$ for $g\in G$ is the end point $p_2$ of the lift
$\widetilde{\gamma}^{-1}_2$ of $\gamma^{-1}$ starting at $pg$.  Since
both paths $\widetilde{\gamma}^{-1}_1\cdot g$ and
$\widetilde{\gamma}^{-1}_2$ are lifts of $\gamma^{-1}$ and have the
same starting point, they must coincide. Hence their end points are
the same points and we have $p_2=p_1g$. Hence this correspondence
is right $G$-equivariant.

To see surjectivity, for any point $p'\in P$ over $x\in\Sigma$, we
choose any path $\widetilde{\gamma}$ in $P$ from $p'$ to a point $p$
in the fibre $P_{x_0}$ over $x_0$. Let $\gamma$ be the path in $\Sigma$
obtained by projecting $\widetilde{\gamma}$ into $\Sigma$. So
$\gamma$ is a path from $x$ to $x_0$. Then $p'\in P$ is the point
corresponding to
$[(x,[\gamma]),p]\in\widetilde{\Sigma}\times_{\Gamma}P_{x_0}$.

For injectivity, suppose two points $[(x_1,[\gamma_1]),p_1]$ and
$[(x_2,[\gamma_2]),p_2]$ in $\widetilde{\Sigma}\times_{\Gamma}P_{x_0}$
correspond to the same point $p'$ in $P$. We first note that
$x_1=\pi(p')=x_2$. Next, let $\widetilde{\gamma_1}$ and
$\widetilde{\gamma_2}$ be lifts in $P$ of $\gamma_1,\gamma_2$ starting
at the same point $p'$ and ending at $p_1$ and $p_2$ in $P_{x_0}$,
respectively. Then a path
$\widetilde{\gamma}_2^{-1}\cdot\widetilde{\gamma}_1$ is the lift of a
loop $\eta=\gamma_2^{-1}\gamma_1$ based at $x_0$. So
$[\gamma_1]=[\gamma_2][\eta]$ and $p_2=[\eta]\cdot p_1$. This implies
that $[(x_1,[\gamma_1]),p_1]=[(x_1,[\gamma_2][\eta]),p_1]
=[(x_2,[\gamma_2]),[\eta]\cdot p_1]=[(x_2,[\gamma_2]),p_2]$. This
proves injectivity.

The proof for the middle horizontal map is similar.
\qed
\enddemo

This proposition reduces the comparison of $G$-bundles over covering
spaces to the comparison of fibres as $\Gamma$-equivariant
$G$-bundles. 

\proclaim{Proposition 7-3}  Let $\Sigma' @>>> \Sigma$ be a covering
space over a connected manifold $\Sigma$ with a base point $x_0$. Two
$G$-principal bundles $\pi_1':P_1 @>>>
\Sigma'$ and $\pi_2':P_2 @>>> \Sigma'$ are isomorphic as $G$-bundles
over the covering space $\Sigma'\!\!/\Sigma$ by an isomorphism
$\alpha:P_1 @>{\cong}>> P_2$ inducing a deck transformation of
$\Sigma'\!\!/\Sigma$ if and only if its restriction to fibres
$\alpha_{x_0}:P_{1,x_0} @>{\cong}>> P_{2,x_0}$ is an isomorphism as
$\Gamma$-equivariant $G$-bundles over the $\Gamma$-set
$\Sigma'_{x_0}$. Thus, we have the following bijective
correspondence\rom{:}
$$
\multline
\biggl\{\foldedtext\foldedwidth{2.6in}{Isomorphism classes of
$G$-principal bundles $\pi:P @>>> \Sigma'$ over a covering space
$\Sigma'\!\!/\Sigma$}\biggr\} \\
\underset{\text{\rm onto}}\to{\overset{1:1}\to\longleftrightarrow}
\biggl\{\foldedtext\foldedwidth{2.3in}{Isomorphism classes of
$\Gamma$-equivariant $G$-principal bundles $P_{x_0} @>>> 
\Sigma'_{x_0}$}\biggr\}.
\endmultline
\tag7-13 
$$
In terms of isomorphisms, we have the following correspondences\rom{:}
$$
\gather
\text{\rm Iso}_G(P_1,P_2)_{\Sigma'\!\!/\Sigma}^{} 
\underset{\text{\rm onto}}\to{\overset{1:1}\to\longleftrightarrow}
\text{\rm Iso}_{\Gamma\text{-}G}^{}(P_{1,x_0},P_{2,x_0}),
\tag7-14  \\
\text{\rm Aut}_G(P)_{\Sigma'\!\!/\Sigma}^{} 
\underset{\text{\rm onto}}\to{\overset{1:1}\to\longleftrightarrow}
\text{\rm Aut}_{\Gamma\text{-}G}^{}(P_{x_0}).\ \ \ \ \ 
\tag7-15
\endgather
$$
\endproclaim
\demo{Proof} The correspondence (7-15) is a special case of (7-14)
with $P_1=P_2$. The correspondence (7-13) follows immediately from
(7-14) by checking whether the sets are empty or not. Thus it suffices
to prove bijectivity of the correspondence in (7-14). 

Let $\alpha:P_1 @>>> P_2$ be a $G$-bundle isomorphism over the
covering space $\Sigma'\!\!/\Sigma$. Since $\alpha$ is an isomorphism
of covering spaces over $\Sigma$, it induces a $\Gamma=\pi_1(\Sigma)$-isomorphism of
fibres $P_{1,x_0}$ and $P_{2,x_0}$. Since $\alpha$ is also a fibre
preserving $G$-equivariant map, the restriction $\alpha_{x_0}:
P_{1,x_0} @>>> P_{2,x_0}$ is a $\Gamma$-equivariant $G$-bundle
isomorphism. This defines the correspondence in (7-14) in one
direction.

Conversely, suppose we are given two $G$-principal bundles $\pi_i: P_i
@>>> \Sigma'$ for $i=1,2$ together with a $\Gamma$-equivariant
$G$-bundle isomorphism $\phi: P_{1,x_0} @>{\cong}>> P_{2,x_0}$ between
the fibres over $x_0$. By Proposition 7-2, we have the following
$G$-bundle isomorphisms over the covering space $\Sigma'\!\!/\Sigma$
between $P_1$ and $P_2$:
$$
P_1 @<{\cong}<< \widetilde{\Sigma}\times_{\Gamma}P_{1,x_0}
@>{\cong}>{1\times_{\Gamma}\phi}> 
\widetilde{\Sigma}\times_{\Gamma}P_{2,x_0} 
@>{\cong}>> P_2. 
$$
This defines the correspondence in the other direction in (7-14). 
It is straightforward to check that the above two correspondences are
inverse to each other. This completes the proof. 
\qed
\enddemo

Combining with Classification Theorem 3-6, Theorem 4-2, and Theorem
4-4, we obtain the following corollary for the case of a connected
covering space $\Sigma' @>>> \Sigma$. 

\proclaim{Corollary 7-4}  Let $\pi_{\sssize\Sigma'\!\!/\Sigma}^{}:
\Sigma' @>>> \Sigma$ be a connected finite covering space 
corresponding to the conjugacy class $[H]$ in
$\Gamma=\pi_1(\Sigma,x_0)$, where $x_0\in\Sigma$ is a base
point. 

\rom{(1)}\qua We have the following bijective correspondence\rom{:}
$$
\biggl\{\foldedtext\foldedwidth{2,7in}{Isomorphism classes of 
$G$-principal bundles $\pi':P @>>> \Sigma'$ over the covering space 
$\Sigma'\!\!/\Sigma$}\biggr\}
\underset{\text{\rm onto}}\to{\overset{1:1}\to\longleftrightarrow}
\text{\rm Hom}(H,G)/(N_{\Gamma}^{}(H)\times G).
\tag7-16
$$

\rom{(2)}\qua Let $\pi':P @>>> \Sigma'$ be a $G$-principal bundle over a
connected covering space $\Sigma'\!\!/\Sigma$ whose holonomy
homomorphism is given by a homomorphism $\rho: H @>>> G$, unique up to
conjugation by $G$. Then the group of $G$-automorphisms of $P$ over
the covering space $\Sigma'\!\!/\Sigma$ \rom{(}inducing deck
transformations on $\Sigma'$\rom{)} is given by
$$
\text{\rm Aut}_G(P)_{\Sigma'\!\!/\Sigma}^{}
\cong\text{\rm Aut}_{\Gamma\text{-}G}^{}(P_{\rho})=
\text{\rm Aut}_{\Gamma\text{-}G}^{}(\Gamma\times_{\rho}G),
\tag7-17
$$
where the structure of the group in the right hand side is described
in \rom{(4-2)} and \rom{(4-5)}. 
\endproclaim

Now we go back to (7-5) and examine the space
$\text{Map}_G^{}(P,M)_{\text{l.c.}}$ of locally constant
$G$-equivariant maps from a $G$-bundle $P$ over $\Sigma'$ into $M$. 

\proclaim{Lemma 7-5} Let $\pi':P @>>> \Sigma'\!\!/\Sigma$ be a
$G$-bundle over a covering space. Then for any $G$-manifold $M$,
restriction of any locally constant $G$-map $\alpha:P @>>> M$ to the
fibre $\alpha_{x_0}: P_{x_0} @>>> M$ over $x_0\in\Sigma$ is a
$\Gamma$-invariant $G$-map. This correspondence gives rise to the
following bijective correspondences\rom{:}
$$
\text{\rm Map}_G^{}(P,M)_{\text{\rm l.c.}} @>{1:1}>{\text{\rm onto}}> 
\text{\rm Map}_G^{}(P_{x_0},M)^{\Gamma}
\underset\text{\rm onto}\to{\overset{1:1}\to\longleftrightarrow}
\Cal S[P_{x_0}\times_G M]^{\Gamma},
\tag7-18
$$
where $\Gamma$ acts on $M$ trivially and, $\Cal S[P_{x_0}\times_G M]$ denotes the set of sections of the
fibre bundle $P_{x_0}\times_G M @>>> \Sigma'_{x_0}$ with fibre $M$. 
\endproclaim
\demo{Proof} First, we show that the correspondence is well
defined. Let $\alpha:P @>>> M$ be a locally constant $G$-equivariant
map. Thus, $\alpha$ is constant on each connected component of
$P$. Since the connected component of $P$ is in bijective
correspondence with $\Gamma$-orbits in the fibre $P_{x_0}$ over
$x_0\in\Sigma$ (the $\Gamma$-orbit corresponding to a connected
component of $P$ is obtained by intersecting the component with the
fibre $P_{x_0}$), the restriction $\alpha_{x_0}: P_{x_0} @>>> M$ is
constant on each $\Gamma$-orbit. Hence $\alpha_{x_0}$ is
$\Gamma$-invariant. This shows that the correspondence in (7-18) is
well defined.

The reverse correspondence is given as follows. For any
$\Gamma$-invariant $G$-map $\phi:P_{x_0} @>>> M$, first we note that
$\phi$ is constant on each $\Gamma$-orbit in $P_{x_0}$. Then, we can
define a locally constant $G$-map $\alpha_{\phi}:P @>>> M$ by letting
the value of $\alpha_{\phi}$ on a connected component $C$ of $P$ to be 
$\phi(C\cap P_{x_0})$. We can easily check that the above two
correspondences are inverse to each other. This completes the proof. 
\qed
\enddemo

In the above correspondence, note that the topological condition of
local constancy of a $G$-map $\alpha:P @>>> M$ translates to an
algebraic condition of $\Gamma$-invariance of its restriction
$\alpha_{x_0}:P_{x_0} @>>> M$ to a fibre. Note that
$\pi^{-1}(x_0)=\Sigma_{x_0}$, where $\pi:\Sigma' @>>> \Sigma$, is a
$\Gamma$-set. 

\demo{\rom{(}Proof of Theorem \rom{7-1}\rom{)}} By (7-5), 
Proposition 7-3, and Lemma 7-5, we can rewrite quantities appearing in
the definition of $\varphi_{\sssize[\Sigma'\!\!/\Sigma]}^{}(M;G)$ in
terms of their restrictions to the fibre over $x_0\in\Sigma$:
$$
\align
\varphi_{\sssize[\Sigma'\!\!/\Sigma]}^{}(M;G)
&=\!\!\!\!\!\!\!\!
\sum_{[P @>>> \Sigma'\!\!/\Sigma]}
\!\!\!\!\!\!\!\!
\varphi\bigl(\text{Map}_G^{}(P,M)_{\text{l.c.}}^{};
\text{Aut}_G^{}(P)_{\Sigma'\!\!/\Sigma}\bigr)\\
&=\!\!\!\!\!\!\!\!\!\!\!\!\!
\sum_{\ \ \ [P_{x_0} \!\!@>>> X]_{\sssize\Gamma\text{-}G}}
\!\!\!\!\!\!\!\!\!\!\!\!\!
\varphi\bigl(\text{Map}_G(P_{x_0},M)^{\Gamma};
\text{Aut}_{\Gamma\text{-}G}^{}(P_{x_0})\bigr),
\endalign
$$
where $[P_{x_0} @>>> X]_{\sssize\Gamma\text{-}G}$ runs over all
isomorphism classes of $\Gamma$-equivariant $G$-principal bundles over
a $\Gamma$-set $X$. But this is precisely our definition of\break
$\varphi_{\sssize[X]}(M;G)$ given in (6-13).

Next, we prove the multiplicativity of
$\varphi_{\sssize[\Sigma'\!\!/\Sigma]}^{}(M;G)$. Let $P @>>>
\Sigma_1'\coprod\Sigma_2'$ be a $G$-bundle over a disjoint union of
two covering spaces of $\Sigma$. Let $P_1$ and $P_2$ be the
restrictions of $P$ to $\Sigma_1'$ and $\Sigma_2'$, respectively. Then
the mapping space splits into a product
$\text{Map}_G(P,M)_{\text{l.c.}}
=\text{Map}_G(P_1,M)_{\text{l.c.}}\times
\text{Map}_G(P_2,M)_{\text{l.c.}}$. When
$\Sigma_1'$ and $\Sigma_2'$ do not have common isomorphic component as
a covering space over $\Sigma$, the automorphism group also splits:
$\text{Aut}_G(P)_{(\Sigma_1'\cup\Sigma_2')/\Sigma}^{}
=\text{Aut}_G(P_1)_{\Sigma_1'/\Sigma}^{}\times
\text{Aut}_G(P_2)_{\Sigma_2'/\Sigma}^{}$, because there are no
nontrivial $G$-bundle maps between a connected component of $P_1$ and
a connected component of $P_2$. Thus,
{\small$$
\multline
\varphi_{\sssize[(\Sigma_1'\cup\Sigma_2')\!/\Sigma]}^{}(M;G)
=\!\!\!\!\!\!\!\!\!\!\!\!\!\!\!\!\!\!\!
\sum_{[P_1\amalg P_2 @>>> \Sigma_1'\amalg\Sigma_2']}
\!\!\!\!\!\!\!\!\!\!\!\!\!\!\!\!\!\!
\varphi\bigl(\text{Map}_G(P_1\amalg P_2, M)_{\text{l.c.}}; 
\text{Aut}_G(P_1\amalg P_2)_{(\Sigma_1'\cup\Sigma_2')
\!/\Sigma}\bigr) \\
=\!\Bigl[\!\!\!\!\!\!\!
\sum_{[P_1 @>>> \Sigma_1']}
\!\!\!\!\!\!
\varphi\bigl(\text{Map}_G(P_1, M)_{\text{l.c.}}; 
\text{Aut}_G(P_1)_{\Sigma_1'/\Sigma}\bigr)\!\Bigr]\!
\Bigl[\!\!\!\!\!\!\!\!\!
\sum_{\ \ \ [P_2 @>>> \Sigma_2']}
\!\!\!\!\!\!\!\!\!
\varphi\bigl(\text{Map}_G(P_2, M)_{\text{l.c.}}; 
\text{Aut}_G(P_2)_{\Sigma_2'\!/\Sigma}\bigr)\!\Bigr]  \\
=\varphi_{\sssize[\Sigma_1'\!/\Sigma]}^{}(M;G)\cdot 
\varphi_{\sssize[\Sigma_2'\!/\Sigma]}^{}(M;G).
\hphantom{\qquad\qquad\qquad\qquad\qquad}
\endmultline
$$}
This completes the proof. 
\qed
\enddemo

Theorem 7-1 allows us to rewrite identities (6-14) and (6-4) in
terms of covering spaces over $\Sigma$.

\proclaim{Corollary 7-6} Let $\Sigma$ be a connected manifold and let
$M$ be a $G$-manifold for a finite group $G$. Then,
$\Sigma$-associated orbifold invariant of the $n$-fold symmetric
orbifold $\varphi_{\sssize\Gamma}^{}(M^n;G\!\wr\!\frak S_n)$, where
$\Gamma=\pi_1(\Sigma)$, can be written in terms of orbifold invariants
of $(M;G)$ associated to up to $n$-fold \rom{(}not necessarily
connected\rom{)} covering spaces $\Sigma' @>>>
\Sigma$ of $\Sigma$. In fact, we have the following identity between
two generating functions\rom{:}
$$
\sum_{n\ge0}q^n\varphi_{\sssize\Gamma}^{}(M^n;G\!\wr\!\frak S_n)
=\!\!\!\!\!\!\sum_{[\Sigma' @>>> \Sigma]}
\!\!\!\!\!q^{\sssize|\Sigma'\!\!/\Sigma|}
\varphi_{\sssize[\Sigma'\!\!/\Sigma]}^{}(M;G).
\tag7-19
$$
Here $\varphi_{\sssize\Gamma}(M^n;G\!\wr\!\frak S_n)$ corresponds to
trivial covering $\varphi_{\sssize[\Sigma\!/\Sigma]}^{}
(M^n;G\!\wr\!\frak S_n)$, and the second summation is over all
isomorphism classes of \rom{(}not necessarily connected\rom{)}
covering spaces $\Sigma'$ over $\Sigma$.

For the case of orbifold Euler characteristics, we even have an
infinite product expansion\rom{:}
$$
\sum_{n\ge0}q^n\chi_{\sssize\Gamma}^{}(M^n;G\!\wr\!\frak S_n)
=\!\!\!\!\!\!\sum_{[\Sigma' @>>> \Sigma]}
\!\!\!\!\!q^{\sssize|\Sigma'\!\!/\Sigma|}
\chi_{\sssize[\Sigma'\!\!/\Sigma]}^{}(M;G)
=\!\!\!\!\!\!\!\!\!\!\!\!\!\!
\prod_{\ \ \ [\Sigma'' @>>> \Sigma]_{\text{conn.}}}
\!\!\!\!\!\!\!\!\!\!\!\!\!
(1-q^{\sssize|\Sigma''\!\!/\Sigma|})
^{-\chi_{\sssize[\Sigma''\!\!/\Sigma]}(M;G)},
\tag7-20
$$
where the product is over all isomorphism classes of finite
connected covering spaces $\Sigma''$ over $\Sigma$.
\endproclaim

The above identity reveals a close and precise connection between
symmetric products and covering spaces in our particular context.

\head
Number of conjugacy classes of subgroups of a given index\strut
\endhead

As before, let $u_r(\Gamma)$ be the number of conjugacy classes of
subgroups of index $r$ in $\Gamma$ and let $j_r(\Gamma)$ be the number
of index $r$ subgroups of $\Gamma$. In this section, we use our
combinatorial formulae obtained as corollaries of our topological
formulae to compute $u_r(\Gamma)$ for various $\Gamma$. Formulae for
$u_r(\Gamma)$ are known for the genus $g+1$ orientable surface group
$\Gamma_{\!g+1}$ \cite{Me}, the free group on $s+1$ generators
$F_{s+1}$, and the genus $h+2$ non-orientable surface group
$\Lambda_{h+2}$ \cite{MP}. However, our method quickly gives the same
result in a uniform way. In fact, we prove a formula which applies to
any group $\Gamma$ to calculate $u_r(\Gamma)$.

First we recall the formula (6-$6'$) which relates the numbers
$u_r(\Gamma)$ to the number of conjugacy classes of
homomorphisms into symmetric groups: 
$$
\sum_{n\ge0}q^n|\text{Hom}(\Gamma,\frak S_n)/\frak S_n|
=\prod_{r\ge1}(1-q^r)^{-u_r(\Gamma)}.
\tag8-1
$$
We study the left hand side in detail using our formula (5-5) on the
number of homomorphisms into wreath product:
$$
\sum_{n\ge0}q^n\frac{|\text{Hom}(\Gamma,G\!\wr\!\frak S_n)|}
{|G|^nn!}
=\exp\biggl[\sum_{r\ge1}\frac{q^r}r
\Bigl\{\!\!\!\!\!\!\!\!\!
\sum\Sb H \\ \ \ \ \ |\Gamma\!/\!H|=r \endSb 
\!\!\!\!\!\!\!\!\!
\frac{|\text{Hom}(H,G)|}{|G|}\Bigr\}\biggr].
\tag8-2
$$
The main result in this section is the following formula. Let 
$\Bbb Z_r=\Bbb Z/r\Bbb Z$. 

\proclaim{Theorem 8-1}  \rom{(1)}\qua For any group $\Gamma$, the generating
function of the number of conjugacy classes of homomorphisms into
symmetric groups is described in terms of subgroups of $\Gamma$ as
follows\rom{:}
$$
\sum_{n\ge0}q^n|\text{\rm Hom}(\Gamma,\frak S_n)/\frak S_n|
=\exp\Bigl[\sum_{m\ge0}\frac{q^m}m
\Bigl\{\sum_{r|m}\!\!\!\!\!\!\!\!\!\!\!\!\!
\sum\Sb H \\ \ \ \ \ \ \ |\Gamma\!/\!H|=m/r\endSb 
\!\!\!\!\!\!\!\!\!\!\!\!\!
|\text{\rm Hom}(H,\Bbb Z_r)|\Bigr\}\Bigr].
\tag8-3
$$

\rom{(2)}\qua The number $u_r(\Gamma)$ of conjugacy classes of index $r$
subgroups of $\Gamma$ satisfies the following recursive
relation in terms of subgroups of $\Gamma$\rom{:}
$$
j_m(\Gamma\times\Bbb Z)=\sum_{r|m}r\cdot u_r(\Gamma)
=\sum_{r|m}\!\!\!\!\!\!\!\!\!\!\!\!\!
\sum\Sb H \\\ \ \ \ \ \ \ |\Gamma\!/\!H|
=m/r \endSb
\!\!\!\!\!\!\!\!\!\!\!\!\!
|\text{\rm Hom}(H_{\text{\rm ab}},\Bbb Z_r)|,
\tag8-4
$$
where $H_{\text{\rm ab}}$ denotes the abelianization of $H$. 
\endproclaim
\demo{Proof} Using the second formula in (2-6), we have 
$$
\align
|\text{Hom}(\Gamma,\frak S_n)/\frak S_n|
&=\frac{|\text{Hom}(\Gamma\times\Bbb Z,\frak S_n)|}{|\frak S_n|}
=\!\!\sum_{\sigma\in\frak S_n}\!\!
\frac{|\text{Hom}\bigl(\Gamma,C(\sigma)\bigr)|}{|\frak S_n|}\\
&=\!\!\!\!\!\sum_{[\sigma]\in{\frak S_n}_*}\!\!\!\!\!
\frac{|\text{Hom}\bigl(\Gamma,C(\sigma)\bigr)|}{|C(\sigma)|},
\endalign
$$
where ${\frak S_n}_*$ denotes the set of conjugacy classes in $\frak
S_n$. It is well known that the conjugacy class $[\sigma]$ in $\frak
S_n$ depends only on the cycle type of $\sigma$. Suppose the cycle
type of $\sigma$ is given by $\prod_{r\ge1}(r)^{m_r}$ with
$\sum_{r}rm_r=n$, then it is well known that its centralizer is a
product of wreath products given by $C(\sigma)\cong\prod_{r\ge1}(\Bbb
Z_r\!\wr\!\frak S_{m_r})$
\cite{M2}.  Then the above formula becomes
$$
|\text{Hom}(\Gamma,\frak S_n)/\frak S_n|
=\!\!\!\!\!\!\!\sum\Sb m_r\ge0 \\ \ \sum rm_r=n \endSb\!\!\!\!\!\!
\prod_{r\ge1}\frac{|\text{Hom}(\Gamma,\Bbb Z_r\!\wr\!\frak S_{m_r})|}
{|\Bbb Z_r\!\wr\frak S_{m_r}|}.
$$
Now we form a generating function. Using the above formula, we have 
$$
\align
\sum_{n\ge0}q^n|\text{Hom}(\Gamma,\frak S_n)/\frak S_n|
&=\sum_{m_r\ge0}\prod_{r\ge1}q^{rm_r}
\frac{|\text{Hom}(\Gamma,\Bbb Z_r\!\wr\!\frak S_{m_r})|}
{|\Bbb Z_r\!\wr\frak S_{m_r}|} \\
&=\prod_{r\ge1}\sum_{m\ge0}q^{rm}
\frac{|\text{Hom}(\Gamma,\Bbb Z_r\!\wr\!\frak S_{m_r})|}
{|\Bbb Z_r\!\wr\frak S_{m_r}|}  
\endalign
$$
Next, we use our crucial ingredient, namely the formula
(8-2). Continuing our calculation, we get 
$$
\align
&=\prod_{r\ge1}
\exp\Bigl[\sum_{d\ge1}\frac{q^{rd}}d
\Bigl\{\!\!\!\!\!\!\!\!\!
\sum\Sb H \\ \ \ \ \ |\Gamma\!/\!H|=d \endSb 
\!\!\!\!\!\!\!\!\!\!
\frac{|\text{Hom}(H,\Bbb Z_r)|}{r}\Bigr\}\Bigr] 
=\exp\Bigl[\sum_{r,d\ge1}\frac{q^{rd}}{rd}
\Bigl\{\!\!\!\!\!\!\!\!\!
\sum\Sb H \\ \ \ \ \ |\Gamma\!/\!H|=d \endSb 
\!\!\!\!\!\!\!\!
\bigl|\text{Hom}(H,\Bbb Z_r)\bigr|\Bigr\}\Bigr]  \\
&=\exp\Bigl[\sum_{m\ge1}\frac{q^m}{m}\Bigl\{\sum_{r|m}
\Bigl(\!\!\!\!\!\!\!\!\!
\sum\Sb H \\ \ \ \ \ |\Gamma\!/\!H|=m/r \endSb
\!\!\!\!\!\!\!\!
\bigl|\text{Hom}(H,\Bbb Z_r)\bigr|\Bigr)\Bigr\}\Bigr].
\endalign
$$
This proves (1). For (2), we simply take logarithm of the right hand
side of (8-1), and compare it with (8-3). Here we may replace
$\text{Hom}(H,\Bbb Z_r)$ by $\text{Hom}(H_{\text{ab}},\Bbb Z_r)$ since
$\Bbb Z_r$ is abelian. The equality $j_m(\Gamma\times\Bbb
Z)=\sum_{r|m}r\cdot u_r(\Gamma)$ is a well known identity. For example,
see [\Cite{St}, p.112]. This completes the proof
\qed
\enddemo

Note that once we have (8-4), we can write down the formula for
$u_r(\Gamma)$ using M\"obius inversion formula. This formula says that
if $g$ is a function defined on $\Bbb N$ and $f(m)=\sum_{r|m}g(r)$,
then $g(m)=\sum_{r|m}\mu(r)f(m/r)$, where $\mu(r)$ is the M\"obius
function given by $\mu(1)=1$, $\mu(n)=(-1)^t$ if $n$ is a product of
$t$ distinct primes, and $\mu(n)=0$ if $n$ is divisible by the square
of a prime. We leave the formula (8-4) in the above form because it is
simpler.

The numbers $u_r(\Gamma_{\!g+1})$ and $u_r(\Lambda_{h+2})$ were
calculated in \cite{Me}, \cite{MP} around mid 1980s. Our formula (8-4)
immediately proves their results which were expressed using M\"obius
inversion formula. But they are equivalent to the formulae below.

\proclaim{Corollary 8-2} The numbers $u_r(\Gamma)$ for
$\Gamma=\Gamma_{\!g+1}, F_{s+1}, \Lambda_{h+2}$ are given in the
following recursive relation\rom{:} 
$$
\gather
\sum_{r|m}r\cdot u_r(\Gamma_{\!g+1})
=\sum_{r|m}j_{\frac{m}{r}}(\Gamma_{\!g+1})r^{2(\frac{gm}r+1)}
\tag8-5 \\
\sum_{r|m}r\cdot u_r(F_{s+1})
=\sum_{r|m}j_{\frac mr}(F_{s+1})r^{\frac{sm}r +1}
\tag8-6 \endgather
$$
{\small
$$\sum_{r|m}r\cdot u_r(\Lambda_{h+2})
=\sum_{r|m}\bigl[\bigl\{j_{\frac mr}(\Lambda_{h+2})-j_{\frac mr}
(\Lambda_{h+2})^+\bigr\}(2,r)r^{\frac{mh}r+1}
+j_{\frac mr}(\Lambda_{h+2})^+r^{\frac{mh}r+2}\bigr],
\tag8-7
$$}%
where $j_r(\Lambda_{h+2})^+$ is the number of index $r$ orientable
subgroups, and $(2,r)$ denotes the greatest common divisor. 
\endproclaim

The proof of this corollary is a straightforward application of
(8-4). We simply note that the abelianizations are given by
$(\Gamma_{\!g+1})_{\text{ab}}\cong\Bbb Z^{2g+2}$,
$(F_{s+1})_{\text{ab}}\cong\Bbb Z^{s+1}$, and
$(\Lambda_{h+2})_{\text{ab}}\cong\Bbb Z^{h+1}\oplus\Bbb Z_2$, and that
among index $m/r$ subgroups of $\Lambda_{h+2}$, there are 
$j_{\frac mr}(\Lambda_{h+2})^+$ subgroups which are isomorphic to
$\Gamma_{\frac{mh}{2r}+1}$, and there are $j_{\frac mr}
(\Lambda_{h+2})-j_{\frac mr}(\Lambda_{h+2})^+$ subgroups isomorphic to
$\Lambda_{\frac{mh}r+2}$.

\Addresses\recd

\enddocument